\documentclass[11pt, reqno]{amsart}
\pdfoutput=1
\usepackage{amssymb,amsmath,latexsym}
\usepackage{epsfig}
\usepackage{mathrsfs}
\usepackage{color}
\usepackage{graphicx}
\usepackage{epic,eepic,wrapfig,color, ifthen} 
\usepackage{url}
\usepackage[colorlinks=true, citecolor=blue, filecolor=cyan, linkcolor=red, urlcolor=magenta]{hyperref}

\usepackage{upgreek}
\usepackage{amsthm}
\usepackage{amsfonts}
\usepackage{pmboxdraw}
\usepackage{verbatim}
\usepackage{color, mathtools, textcomp}
\usepackage{mathtools}
\usepackage{kotex}
\usepackage{cite}
\usepackage[normalem]{ulem}
\usepackage{enumerate}
\mathtoolsset{showonlyrefs}
\usepackage[bottom]{footmisc}

\newcommand{\cal}[1]{\mathcal{#1}}

\def\Up{{\Upsilon}}
\def\ep{{\mathsf{d}}}
\def\up{{\upsilon}}

\def\1{{\bf 1}}
\def\nn{\nonumber}

\def\sA {{\cal A}} \def\sB {{\cal B}} \def\sC {{\cal C}}

\def\sF {{\cal F}}
  \def\sI {{\cal I}}
\def\sJ {{\cal J}} 
\def\sL {{\cal L}}
\def\sM {{\cal M}}
\def\sN {{\cal N}} \def\sO {{\cal O}}
 \def\sQ {{\cal Q}}  

 \def\sU {{\cal U}}

\def\R {{\mathbb R}}  
 \def\Z {{\mathbb Z}}
\def\E {{\mathbb E}} 

\def\U{{\mathrm{US}}}
\def\L{{\mathrm{LS}}}

\def\as{Assumption \ref{ass0}}
\def\ass{Assumption \ref{assinf}}
\def\asss{Assumption \ref{assinf2}}

\def \VRDM {$\mathrm{VRD}_{R_0}(M)$}
\def \VRDo {${\mathrm {VRD}}_{R_0}(U)$}
\def \VRDi {${\mathrm {VRD}}^{R_\infty}$}
\def \Cho {$\mathrm{Ch}_{R_0}(U)$}
\def \Chi {$\mathrm{Ch}^{R_\infty}$}
\def \Tailo {${\mathrm{TJ}}_{R_0}(\psi,U)$}
\def \Tailol {${\mathrm{TJ}}_{R_0}(\psi,\le,U)$}
\def \Tailog {${\mathrm{TJ}}_{R_0}(\psi,\ge,U)$}
\def \Taili {$\mathrm{TJ}^{R_\infty}(\psi)$}
\def \Tailil {$\mathrm{TJ}^{R_\infty}(\psi, \le)$}
\def \Tailig {$\mathrm{TJ}^{R_\infty}(\psi, \ge)$}
\def \Eo {$\mathrm{E}_{R_0}(\phi,U)$}
\def \Ei {$\mathrm{E}^{R_\infty}(\phi)$}

\def \NDLo {${\mathrm{NDL}}_{R_0}(\phi,U)$}
\def \NDLi {${\mathrm{NDL}}^{R_\infty}(\phi)$}

\def\tO {{\tilde{\Omega}}}
\def\tP {{\tilde{\P}}}
\def\tE {{\tilde{\E}}}

\def\NDL {\mathrm {NDL}}

\def\PI{\mathrm{PI}}
\def\CS{\mathrm{CS}}
\def\FK{\mathrm{FK}}

\def\Tail{\mathrm{TJ}}

\numberwithin{equation}{section}
\def\qed{{\hfill $\Box$ \bigskip}}

\def\FF{{\mathcal F}}
\def\EE{{\mathcal E}}
\def\GG{{\mathcal G}}

\def\R{{\mathbb R}}

\def\P{{\mathbb P}}
\def\N{{\mathbb N}}
\def\bL{{\mathbb L}}
\def\eps{\varepsilon}

\def\wt{\widetilde}
\def\pf{\noindent{\bf Proof. }}

\def\ll{{\Lambda}}
\def\vt{{\vartheta}}

\def\vk{{\varkappa}}

\DeclareMathOperator*{\esssup}{ess\,sup}
\DeclareMathOperator*{\essinf}{ess\,inf}

\theoremstyle{plain}
\newtheorem{thm}{Theorem}[section]
\newtheorem{lem}[thm]{Lemma}
\newtheorem{cor}[thm]{Corollary}
\newtheorem{remark}[thm]{Remark}
\newtheorem{prop}[thm]{Proposition}
\newtheorem{defn}[thm]{Definition}
\newtheorem{definition}[thm]{Definition}
\newtheorem{assumption}[thm]{Assumption}
\newtheorem{example}[thm]{Example}

\theoremstyle{definition}
\newtheorem*{eg*}{Example}
\newtheorem*{egs*}{Examples}
\newtheorem*{def*}{Definition}

\theoremstyle{remark}

\begin{document}
	\title{
		General Law of iterated logarithm for Markov processes: Limsup law}	
	
	\author{Soobin Cho, Panki Kim, Jaehun Lee}\thanks{This research is  supported by the National Research Foundation of Korea(NRF) grant funded by the Korea government(MSIP) (No. 2016R1E1A1A01941893).
	}

	\address[Cho]{Department of Mathematical Sciences,
		Seoul National University,
		Seoul 08826, Republic of Korea}
	\curraddr{}
	\email{soobin15@snu.ac.kr}
	
	\address[Kim]{Department of Mathematical Sciences and Research Institute of Mathematics,
		Seoul National University,
		Seoul 08826, Republic of Korea}
	
	\curraddr{}
	\email{pkim@snu.ac.kr}
	
	\address[Lee]{
		Korea Institute for Advanced Study,  
		Seoul 02455,
		Republic of Korea}
	
	\curraddr{}
	\email{hun618@kias.re.kr}

	\maketitle
	
	\begin{abstract}
		In this paper, we discuss general criteria of limsup law  of iterated logarithm (LIL) for continuous-time Markov processes. 
		We consider minimal assumptions for LILs to hold at zero 
		(at infinity, respectively) in general metric measure spaces. 
		We establish LILs under local assumptions near zero (near infinity, respectively) on uniform bounds of the 
		 expectations of first exit times 
		 from balls in terms of a function $\phi$
		and  uniform bounds on the tails of the jumping kernel in terms of a function $\psi$. 
		The main result is that a simple ratio test in terms of the functions $\phi$ and $\psi$ completely determines whether there exists a positive non-decreasing  function $\Psi$ such that $\limsup |X_t|/\Psi(t)$ is positive and finite a.s., or not.
		Our results cover a large class of subordinate diffusions, jump processes with mixed polynomial local growths, 
		jump processes with singular jumping kernels
		and  random conductance models with long range jumps.
		
		\medskip
		
		\medskip
		
		
	\end{abstract}
	\allowdisplaybreaks

	\section{Introduction and general result}\label{s:intro}
	
	Law of the iterated logarithm (LIL) for a stochastic process describes the magnitude of the fluctuations of its sample path behaviors. 
	In addition to the law of large numbers and the central limit theorem, 
	the LIL is considered  the fundamental limit theorem in Probability theory. See \cite{Fe45} and the references therein.

	Let $Y:=(Y_t)_{t\ge0}$ be a non-trivial strictly $\beta$-stable process on ${\mathbb R}^d$ with $0<\beta \le 2$ in the sense of \cite[Definition 13.1]{Sa13}. 
	 Then $Y$ satisfies the following well known limsup LIL: If $\beta=2$, then there exist constants $ c_1,c_2\in (0,\infty)$ such that 
	\begin{align}
		\label{e:ieq1}
		\limsup_{t \to 0 \;\; (\text{resp. }  t\to \infty) } \frac{\sup_{0<s\le t} |Y_s|}{(t\log|\log t|)^{1/
				\beta}}
		=c_1, \quad \limsup_{t \to 0 \;\; (\text{resp. }  t\to \infty) } \frac{ |Y_t|}{(t\log|\log t|)^{1/\beta}}
		=c_2 \;\;\mbox{ a.s.}
	\end{align}
	Otherwise,  i.e., $\beta<2$, then for every positive non-decreasing function $\Psi$,
	\begin{equation}\label{e:ieq2}
		\begin{split}
			&\int_0\Psi(t)^{-\beta}dt < \infty \;\;\text{ or }\;\;=\infty  \quad \Big({\rm resp.}  \int^\infty\Psi(t)^{-\beta}dt < \infty \;\;\text{ or }\;\;=\infty \Big)\\
			& \Leftrightarrow \;\; \limsup_{t \to 0\;\; (\text{resp. }  t\to \infty)} \frac{\sup_{0<s \le t}|Y_s|}{\Psi(t)} = \limsup_{t \to 0\;\; (\text{resp. }  t\to \infty)} \frac{|Y_t|}{\Psi(t)} =0 \;\;\text{ or }\;\;\infty, \;\; \mbox{ a.s.}
		\end{split}
	\end{equation}
	See \cite[Chapters 47, 48]{Sa13} and the references therein.

	The limsup LIL of the second type (for  $|Y_t|$) in \eqref{e:ieq1} was obtained by  Khintchine  \cite{Kh33}, Kolmogorov \cite{Ko29},  Hartman and Wintner \cite{HW41}  
	and L\'evy \cite{
		Levy2},    for various random walks  and Brownian motions in any dimension. 
	The  limsup LIL of the first type  (for $\sup_{0<s\le t} |Y_s|$)  in \eqref{e:ieq1} was obtained by Strassen \cite{Str64} for some random walks and  Brownian motion in $\R^1$. Later, Barlow and Perkins \cite{BP88} and Barlow \cite{Ba98} showed that the limsup LIL of the second type holds with some $\beta \ge 2$ for Brownian motions on some fractals including Sierpinski gasket or carpet, and   Bass and Kumagai \cite{BK} showed that the first type holds for Brownian motions on some fractals and Riemannian manifolds. (Although these results only consider either $t$ tends to zero or infinity, one can  prove the results for the other direction by modifying their proofs.)
	For jump processes (Markov processes whose trajectories have discontinuities), the limsup LIL of the second type at infinity was done for L\'evy processes in $\R$
	with finite second moment   by Gnedenko \cite{Gn43}, symmetric stable processes in $\R$ without large jumps  by Griffin \cite{Gri85},  and  non-L\'evy processes in $\R^d$
	with finite second moment by Shiozawa and Wang \cite{SW19}, and Bae, Kang, Kim and Lee \cite{BKKL19a}. 
	Very recently, the work on 
	limsup LIL at infinity
	in  \cite{BKKL19a} is  extended by same authors to metric measure space in  \cite{BKKL19b}.

	 The limsup LIL of type \eqref{e:ieq2} was first  obtained by Khintchine \cite{Kh38}  in $\R$.  Those results were extended to subordinators and some L\'evy processes in $\R$ by Fristedt \cite{Fr67, Fr71}. 
	  Recently, the limsup LILs of type \eqref{e:ieq2} at zero were discussed for more general L\'evy processes in $\R$ by Savov \cite{Sav} and
	some L\'evy-type processes in $\R$ by Knopova and Schilling \cite{KS}. We refer to \cite{KKW17} for a multi-dimensional version.
	
	In the study of the limsup laws, the following natural and important question has been raised and partially answered. Cf. questions in \cite{EL05, He69, Ro68, Pr81b
	} for random walks.   We denote by $\sM_+$ the set of all positive non-decreasing functions defined on $(0,c_1] \cup [c_2,\infty)$ for some $c_2 \ge c_1>0$.
	
	\medskip
	
	\setlength{\leftskip}{1cm}
	\hspace{-1.2cm}$\mathbf{(Q)}$ 
	How do we determine that a limsup law (either at zero or at infinity) of a given process $Y$ is 
	of type \eqref{e:ieq1}, or else of type \eqref{e:ieq2}? More precisely, what is a necessary and sufficient condition for the existence of $\Psi \in \sM_+$ such that $\limsup$ of $|Y_t|/\Psi(t)$ converges to a positive, finite and deterministic value  as $t$ tends to zero or infinity?
	
	\setlength{\leftskip}{0cm}

	\medskip
	
	Feller \cite{Fe68} studied the above question for symmetric  random walks which belong to the domain of attraction of the normal distribution, and established an integral test  as an answer. He also found the exact form of $\Psi$ when it exists. We also refer to  Kesten's work \cite{Ke72}. For continuous time processes, Fristedt \cite{Fr71} partially answered the above question for one-dimensional symmetric L\'evy processes. Then Wee and Kim \cite{WK88} 
	established a ratio test for entire one-dimensional  L\'evy processes which determine whether $\limsup |Y_t|/\Psi(t) \in (0, \infty)$.

	\smallskip
	
	The purpose of this paper is to understand asymptotic behaviors of a given Markov process by establishing limsup law  of iterated logarithms for both near zero and near infinity under some minimal assumptions. 
	The main contribution of this paper is that we answer the above question ${\bf (Q)}$ completely for a large class of Markov processes including random conductance models with stable-like jumps (see Theorems \ref{t:limsup0-1}--\ref{t:limsupinf-1} and Section \ref{s:RCM} below).
	
	\smallskip
	Formally, our answer is 
	
	\medskip
	
	\setlength{\leftskip}{1cm}
	
	\hspace{-1.2cm}$\mathbf{(A)}$ 
	There is a dichotomous classification on continuous time Markov processes (satisfying a near diagonal lower bound estimate for the Dirichlet heat kernel):	Based on whether the tails of the jumping kernel and the mean exit times (the expectation of first exit times)
	 from balls are comparable or not,  
	we can categorize  Markov processes into two non-overlapping classes $\sA^\infty_1$ and $\sA^\infty_2$:
	$$
	\sA^\infty_1=\Big\{\text{For all }\Psi \in \sM_+, \;\, \limsup_{t \to \infty} \frac{d(x,X_t)}{\Psi(t)} 
	\text{ is either $0$ for all $x$, or $\infty$ for all $x$.}\Big\},
	$$
	$$
	\sA^\infty_2=\Big\{\text{There is $\Psi \in \sM_+$ such that } \limsup_{t \to \infty} \frac{d(x,X_t)}{\Psi(t)}
	\text{ is  deterministic  and in $(0, \infty)$ for all $x$.}\Big\}
	$$

	\noindent If a given Markov process $X$ is in the class $\sA^\infty_1$, then there is an integral test that determines whether the above $\limsup$ is zero or infinite. 
	If a given Markov process $X$ is in the class $\sA^\infty_2$, then there is 
	a natural explicit function $\Psi$ and 
	another integral test to determine when
	one can take such specific $\Psi$. An analogous   result holds for $\limsup$ LIL at zero.

	\setlength{\leftskip}{0cm}
	
	\smallskip
	
	The reader will see that processes near  the borderline of the above classes 
	are so-called Brownian-like jump processes, i.e., 
	jump processes  with high intensity of small jumps   and ones  with  low intensity of  large jumps.

	\medskip

	Assumptions in this paper are motivated by the second named author's previous paper \cite{KKW17}. In \cite{KKW17}, the authors established liminf and limsup LILs (of type \eqref{e:ieq2}) for Feller processes on a general metric measure space enjoying mixed stable-like heat kernel estimates (see Assumption 2.1 therein). Recently, 
	relationships among 
	those heat kernel estimates, and certain conditions on the jumping kernel and 
	the mean exit time from open balls are extensively studied.
	See, e.g.  \cite{BKKL19b, BKKL19a, 
		 CKW16a, CKW19, CKW20, BGK09, GHH18}. We adopt this framework and consider localized and relaxed conditions not only on the heat kernel but also on the jumping kernel and the mean exit time. 	
		We emphasize that, unlike the references mentioned above, we do not assume a weak lower scaling property of the scale function $\phi$ in several statements. (See Definition \ref{d:ws} for the notion of weak lower scaling property.) 
	However, we still need the weak lower scaling property in some  other statements, especially ones concerning long time behaviors. 
	
	Our assumptions are  weak enough so that our results cover 	a lot of Markov processes including  random conductance models with long range jumps,	jump processes with diffusion part,  jump processes with
	low intensity of small jumps,
	some  non-symmetric processes and  processes with singular jumping kernels. See the examples in Sections \ref{s:dset} and \ref{s:RCM}, and  also see \cite{CKL} for further examples.   In particular,  the class of Markov processes considered in this paper extends the corresponding results of \cite{KKW17} significantly.  Moreover, metric measure spaces in this paper can be random, disconnected and highly space-inhomogeneous (see Definition \ref{d:VD}).

	\medskip

	Let us now describe the main result of this paper   precisely and, at the same time, fix the setup and the  notation of the paper.

	Throughout this paper, 
	 we assume that $(M,d)$ is a locally compact separable metric space with a base point $o \in M$, and $\mu$ is a positive Radon measure on $M$ with full support. Denote by $B(x,r):= \{ y \in M : d(x,y) <r  \}$ and $V(x,r):= \mu(B(x,r))$ an open ball in $M$ and its volume, respectively. 
	We add a cemetery point $\partial$ to $M$ and define $M_\partial:=M \cup \{\partial\}$.   Any function $f$ on $M$ is extended to $M_\partial$ by setting $f(\partial)=0$.

	For an open set $U \subset M$ and $x \in M$, we denote by  $\updelta_U(x)$ the distance between $x$ and $M\setminus U$. 
	 { We define a map $\ep:M \to [1,\infty)$ by }
	\begin{equation*}
		\ep(x):=d(x,o)+1.
	\end{equation*}
{ When $M=\R^d$ and $o$ is the origin, $\ep(x)$ equals to $|x|+1$. Note that the map $\up \mapsto \ep(x)^\up$ is non-decreasing on $\R$ since $\ep(x) \ge 1$.}

	We now introduce 
	following versions of volume doubling  property.
	We denote $a \land b$ for $\min\{a,b\}$ and $a \vee b$ for $\max\{a,b\}$.
	
	\begin{defn}\label{d:VD}
		{\rm (i) 
			For an open set $U \subset M$ and $R_0 \in (0,\infty]$, we say that \textit{the interior volume doubling and reverse doubling property near zero}   \VRDo \ holds (with $C_V$) 
			if there exist constants $C_V \in (0,1)$, $d_2 \ge d_1>0$ and  $c>1$ such that for all $x \in U$ and $0<s\le r<R_0 \land (C_V\updelta_U(x))$,
			\begin{equation}\label{VD1_0}
			c^{-1} \Big( \frac{r}{s} \Big)^{d_1} \le \frac{V(x,r)}{V(x,s)} \le  c \Big( \frac{r}{s} \Big)^{d_2}.
			\end{equation}
			
			\noindent (ii) {	For  $R_\infty\ge1$,  we say that \textit{the weak  volume doubling and reverse doubling property near infinity} \VRDi \ holds} if there exist constants $\up \in (0,1)$, $d_2 \ge d_1>0$ and  $c>1$ such that \eqref{VD1_0} holds for all $x \in M$ and $ r\ge s> R_\infty \ep(x)^\up$.
	} \end{defn}

\begin{figure}
	\includegraphics[width=0.36\columnwidth]{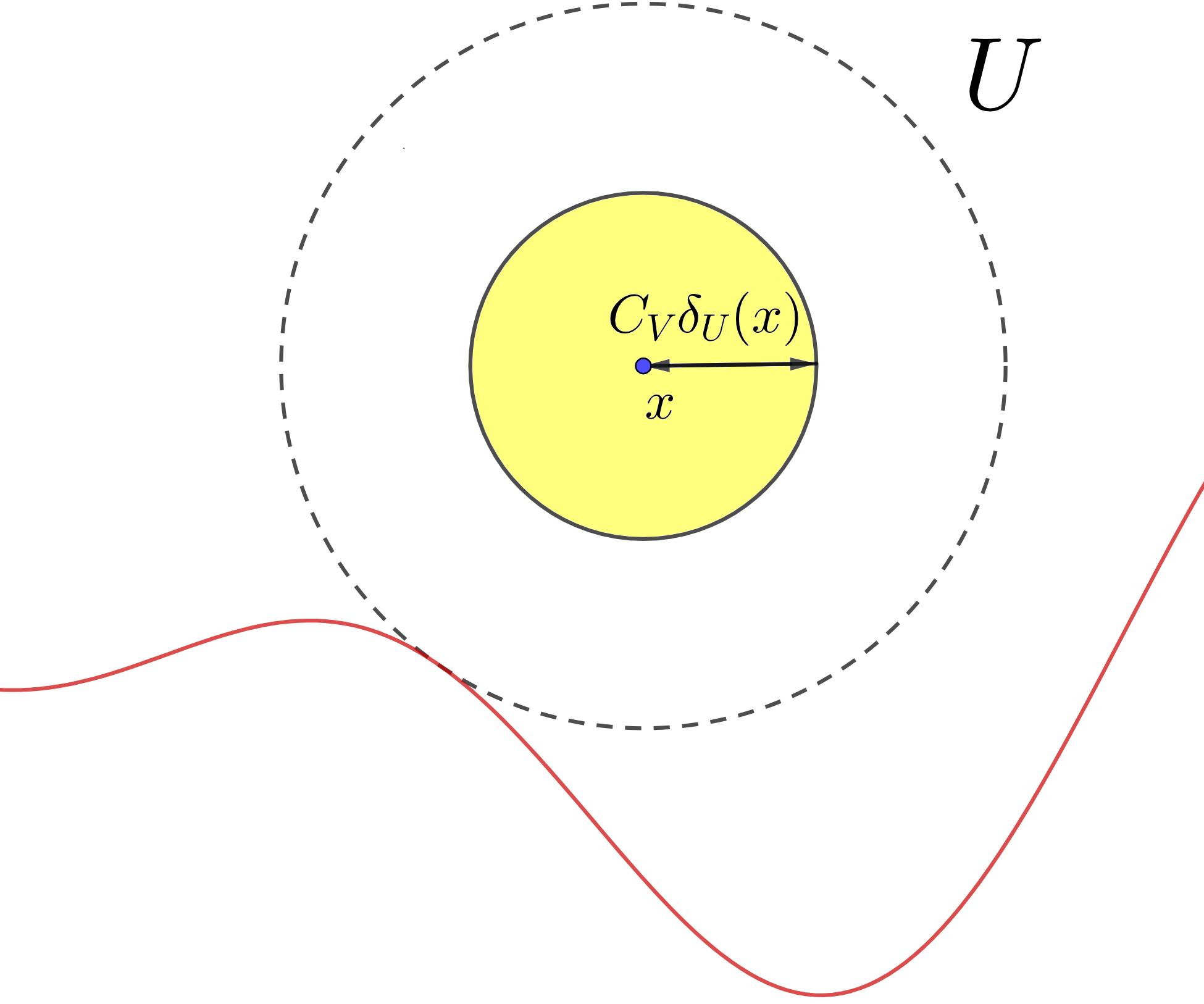}  \hspace{20mm}	\includegraphics[width=0.32
	\columnwidth]{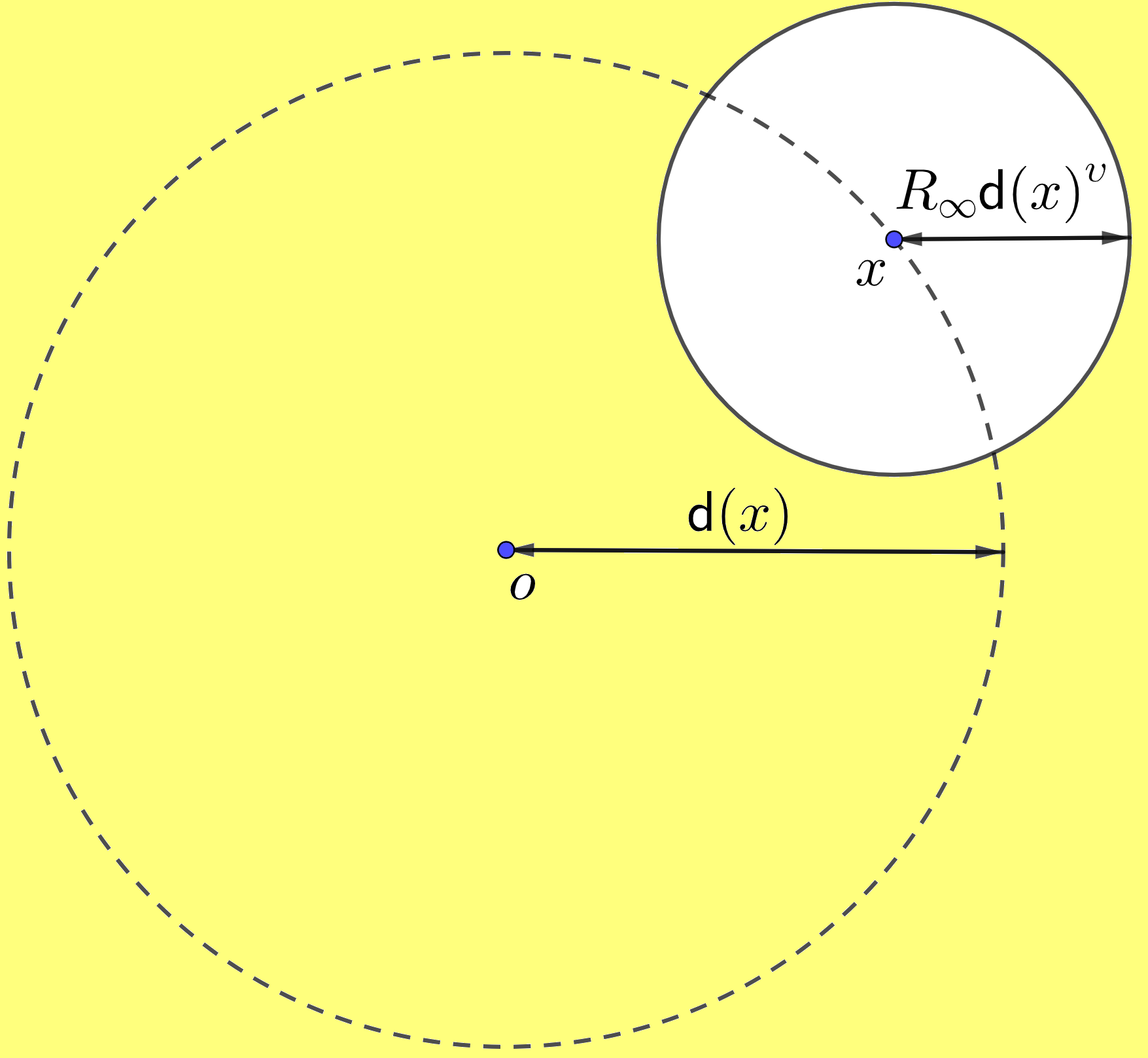} \vspace{-2mm}
		\caption{Range of $r$ in local conditions}
\end{figure}

	Note that, under \VRDo \ (resp. \VRDi), it holds that with a constant $\ell:=(2/c_\mu)^{1/d_1}$,
	\begin{equation}\label{RVD}
		\frac{V(x,  r)}{V(x,r/\ell)} \ge 2  \;\; \mbox{ for all} \; x\in U \mbox{ and } 0<r < R_0 \land (C_V\updelta_U(x))  \;\; \big(\text{resp. } x\in M \mbox{ and } r >  \ell R_\infty \ep(x)^\up \big).
	\end{equation}
	
	Next, we define local  versions of the so-called chain condition.
	\begin{defn}\label{d:Ch}
		{\rm (i) For an open set $U \subset M$ and $R_0 \in (0, \infty]$, we say that \textit{the chain condition near zero} \Cho \ holds  if there exists a constant $A \ge 1$ such that for all $x,y \in U$ with $d(x,y) < R_0$ and $n \in \N$, there is a sequence $(z_{i})_{0 \le i \le n}\subset M$ such that
			\begin{equation}\label{e:Chain}
				z_{0}=x, \;\; z_{n}=y \;\; \text{and} \;\; d(z_{i-1},z_{i}) \le  \frac{A}{n}d(x,y) \quad \text{for all} \;\; 1 \le i \le n.
			\end{equation}

			\noindent	(ii) For $R_\infty \ge 1$, we say that \textit{the weak chain condition near infinity}  \Chi \ holds   if there exist  constants  $\up\in(0,1)$ and $A \ge 1$ such that for all $x,y \in M$ and $n \in \N$ with $d(x,y)/n > R_\infty(\ep(x) \vee \ep(y))^\up$,  there is a sequence $(z_{i})_{0 \le i \le n} \subset M$ satisfying \eqref{e:Chain}.
	} \end{defn}
	
	\begin{remark}
		{\rm If $(M,d)$ is geodesic, then $\mathrm{Ch}_\infty(M)$ clearly holds.}
	\end{remark}
	
	Let $X = (\Omega, \sF_t, X_t, \theta_t, t \ge 0;  \P^x, x \in M_\partial)$ be a Borel standard Markov process on $M_\partial$ with the lifetime $\zeta:=\inf\{t>0:X_t=\partial\}$. Here $(\theta_t, t \ge 0)$ is  the shift operator with respect to  $X$, which is defined as $X_s ( \theta_t \, \omega) = X_{s+t}(\omega)$ for  $t,s>0$ and $\omega \in \Omega$. A family of $[0,\infty]$-valued random variables $\mathfrak{A}=(\mathfrak{A}_t, t \ge0)$ is called a  \textit{positive continuous additive functional} (PCAF) (in the strict sense) of $X$, if there exists $\Lambda \in \sF_\infty$ with $\P^x(\Lambda)=1$ for all $x \in M$ such that the following conditions are satisfied: (i) For each $t \ge 0$, $\mathfrak{A}_t|_\Lambda$ is $\sF_t|_\Lambda$-measurable and (ii) For any $\omega \in \Lambda$, $\mathfrak{A}_0(\omega)=0$, $|\mathfrak{A}_t(\omega)|<\infty$ for $t<\zeta(\omega)$, $\mathfrak{A}_t(\omega)=\mathfrak{A}_{\zeta(\omega)}(\omega)$ for $t \ge \zeta(\omega)$, $t \mapsto \mathfrak{A}_t(\omega)$ is a continuous function on $[0,\infty]$ and  $\mathfrak{A}_{t+s}(\omega)=\mathfrak{A}_t(\omega)+\mathfrak{A}_s(\theta_t\,\omega)$ for all $s,t\ge 0$. See \cite[Chapter 5]{FOT11}.

	According to \cite{BJ73}, 
	since $X$ is a Borel standard process on $M$, it has a L\'evy system $(N, H)$. Here $N(x,dy)$ is a kernel on $(M_\partial, \sB(M_\partial))$ and $H$ is a PCAF of $X$ with bounded $1$-potential. 
	{\it We assume that $X$ has a L\'evy system $(N,H)$ such that
	the Revuz measure of $H$ is given by $\nu_H(x)\mu(dx)$ for a measurable function $\nu_H(x)$.}
 Let $J(x, dy)=N(x,dy)\nu_H(x)$. Then, by the L\'evy system formula, for any non-negative Borel function $F$ on $M \times M_\partial$ vanishing on the diagonal, it holds that
	\begin{equation}\label{e:system}
		\E^z \bigg[ \sum_{s \le t} F(X_{s-}, X_s)  \bigg] = \E^z \left[ \int_0^t \int_{M_\partial} F(X_s,y)J(X_s,dy)ds \right], \quad z\in M, \; t>0.  
	\end{equation}
	The  measure $J(x,dy)$ on $M_\partial$ is called the \textit{L\'evy measure} of $X$. See \cite{Wat64}. Here, we note that the killing term  $J(x, \partial)$ is  included in the  L\'evy measure. 
	Also, we emphasize that 
	$J(x,dy)$ can be identically zero and $J(x,dy)$ may not be absolutely continuous with respect to  $\mu$. 
	
	\smallskip
	
	For an open set $D \subset M$, denote by $\tau_D:= \inf \{ t>0 : X_t \in M_\partial \setminus D \}$ \textit{the first exit time} of $X$ from $D$.
	The subprocess $X^D$, defined by $X_t^D := X_t \, \1_{\{\tau_D > t\}} + \partial \, \1_{\{\tau_D \le t\}}$,  is called the \textit{killed process} of $X$  upon  leaving $D$. 
	We call  a measurable function $p^D: (0,\infty) \times D \times D \to [0,\infty]$ the \textit{heat kernel} (or the transition density) of $X^D$ if the followings hold: 
	
	\smallskip
	
	(1) $\E^x [f(X^D_t)] = \int_{D} p^D(t,x,y)f(y) \mu(dy)$ for all $t>0$, $x \in D$ and $f \in L^\infty(D;\mu)$;
	
	\smallskip
	
		(2) $p^D(t+s,x,y) = \int_{D} p^D(t,x,z) p^D(t,z,y) \mu(dz)$ for all $t,s>0$ and $x,y \in D$.
	
	\smallskip
	
\noindent When the process $X=X^M$ has a heat kernel, 	we simply write $p(t,x,y)$ instead of  $p^M(t,x,y)$.

	Hereinafter, we say that a real-valued function $f$ on an interval is increasing (resp. decreasing) if $f(s)<f(r)$ ($f(s)>f(r)$) for all $s<r$ in the domain of $f$. If $f(s)\le f(r)$ ($f(s)\ge f(r)$) for all $s<r$, we say that $f$ is non-decreasing (resp. non-increasing).
	
	{\it Throughout this paper, we consider increasing and continuous functions  $\phi,\psi:(0,\infty) \to (0,\infty)$  such that $ \lim_{r \to 0}\phi(r)=0$, $\lim_{r \to \infty}\phi(r)=\infty$ and
\begin{align}\label{e:phipsi}
\phi(r) \le \psi(r) \quad \text{for all $r>0$.}
\end{align}
}
Using $\phi$ and $\psi$, we now introduce three types of local conditions: Tail estimates on the L\'evy measure $J(x, dy)$, estimates on the mean exit times from balls and near diagonal lower estimates of heat kernels.  (See, e.g. \cite{CKW16b} for their global versions.)

	\begin{definition}
		{\rm  Let  $R_0 \in (0,\infty]$ be a constant and $U \subset M$ be an open set.
			
\noindent	(i)  We say that  \Tailo \ holds 
(with  $C_1$) 
if there exist  constants  $C_1\in (0,1)$ and $c>1$  such that for all $x \in U$ and  $0<r<R_0 \wedge (C_1\updelta_U(x))$,
\begin{equation}\label{e:Tail_0}
	\frac{c^{-1}}{\psi(r)} \le J(x,M_\partial \setminus B(x,r)) \le \frac{c}{\psi(r)}.
\end{equation}
We say that \Tailol \ (resp. \Tailog) holds  if the upper bound (resp. lower bound) in \eqref{e:Tail_0} holds for all $x \in U$ and $0<r<R_0 \wedge (C_1\updelta_U(x))$.

\smallskip

\noindent (ii) We say that \Eo \ holds
 (with  $C_1$) 
if there exist constants $C_1 \in (0,1)$, $A_0>0$ and $c> 1$ such that for all $x \in U$ and $0<r<R_0 \land (C_1\updelta_U(x))$,
\begin{equation}\label{e:Eo}
	c^{-1}\phi(A_0 r) \le \E^x[\tau_{B(x,r)}]  \le c \phi(A_0 r).
\end{equation}

\noindent (iii)  We say that \NDLo \ holds 
(with  $C_2$)
 if there exist constants $C_2, \eta \in (0,1)$ and  $c>0$ such that for all $x \in U$ and $0<r<R_0 \land (C_2\updelta_U(x))$, the heat kernel $p^{B(x,r)}(t,y,z)$ of $X^{B(x,r)}$ exists and 
\begin{equation}\label{e:NDL_inf}
	p^{B(x,r)}(\phi(\eta r),y,z) \ge \frac{c}{V(x, r)} \quad \text{for all}\;\, y,z \in B(x, \eta^2r).
\end{equation} 
}
	\end{definition}

\begin{definition}\label{d:inf}
	{\rm 
		Let $R_\infty \ge 1$ be a constant.
		
		\noindent  (i)  We say that  \Taili \ holds  if there exist constants $\up\in(0,1)$ and $c>1$ such that \eqref{e:Tail_0} holds for all $x \in M$ and $r >R_\infty \ep(x)^\up$. We say that \Tailil \ (resp. \Tailig) holds if the upper bound (resp. lower bound) in \eqref{e:Tail_0} holds for all $x \in M$ and $r>R_\infty \ep(x)^\up$.

\smallskip

\noindent (ii)  We say that \Ei \ holds if there exist constants  $\up \in (0,1)$, $A_0>0$ and  $c> 1$ such that \eqref{e:Eo} holds for all  $x \in M$ and $r>R_\infty \ep(x)^\up$.

\smallskip

\noindent (iii)  We say that \NDLi  \ holds  if there exist
constants $\up,\eta \in (0,1)$ and $c>0$ such that for all $x \in M$ and  $r >R_\infty \ep(x)^\up $,  the heat kernel $p^{B(x,r)}(t,y,z)$ of $X^{B(x,r)}$ exists and satisfies \eqref{e:NDL_inf}.}
\end{definition}

	\begin{remark}\label{r:inf}
	{\rm $\mathrm{VRD}_\infty(M)$ implies \VRDi \ for any $R_\infty \ge 1$ since \eqref{VD1_0} holds for all $x \in M$ and $0<s \le r$.
		Similar results hold concerning  other conditions $\mathrm{Ch}$, $\mathrm{TJ}$, $\mathrm{E}$ and $\mathrm{NDL}$. 
	}
\end{remark}

	\begin{remark}\label{r:tail}
		{\rm (i) 
			The inequality $\phi(r) \le \psi(r)$ in \eqref{e:phipsi} is quite natural
	under the assumptions \VRDo,  \Tailog \  and 
	\Eo. See \cite[Remark 2.7]{BKKL19b} and the paragraphs below it.

			\noindent (ii)	{ By \eqref{e:phipsi},   \Tailol \ implies $\mathrm{TJ}_{R_0}(\phi,  \le,U)$ and \Tailil \ implies $\mathrm{TJ}^{R_\infty}(\phi, \le)$.}

			\noindent (iii) If  $X$ is a conservative diffusion process, then   $\Tail_\infty(\phi,\le,M)$    holds. 
			Here, the conservativeness is required since the killing term $J(x, \partial)$ is included in \eqref{e:Tail_0}.
			
			\noindent (iv) For $\alpha \in (0,1)$ and  any open set $D \subset \R^d$,   a killed isotropic $2\alpha$-stable process in $D$ (whose  infinitesimal generator is the Dirichlet fractional Laplacian $-(-\Delta)^{\alpha}\vert_{D}$) satisfies $\Tail_\infty(r^{2\alpha},D)$.
			  Moreover, by \cite[(2.22),(2.23) and (3.7)]{KSV20}, for $\beta, \delta \in (0,1)$ and  a bounded open set $D \subset \R^d$, a subordinate killed stable process in $D$ with the  infinitesimal  generator $-\big( (-\Delta)^\delta\vert_{D}\big)^\beta$ satisfies $\Tail_\infty(r^{2\beta \delta },D)$. We emphasize again that the killing term  $J(x, \partial)$ is included in \eqref{e:Tail_0}.
		}
	\end{remark}

	As you see from the above Definitions
	\ref{d:VD}(ii), 
	\ref{d:Ch}(ii) 
	and
	\ref{d:inf}, 
	our conditions at infinity are weaker by adding the restriction $r >R_\infty \ep(x)^\up $
	with $\up \in (0,1)$. 
	Thanks to such weak assumptions at infinity, 
	we can cover  LILs for random conductance models  at infinity. See Section \ref{s:RCM} below.

	We now introduce (local) weak lower and upper scaling properties for positive functions.
	
	\begin{defn}\label{d:ws}
		{\rm  Let $g$ be a given positive function defined on an interval, and  $\beta_1,\beta_2,c_1,c_2$ be  constants. For $a \in (0,\infty]$, 		 	we say that $\L_a(g,\beta_1, c_1)$ (resp. $\L^a(g,\beta_1, c_1)$) holds if	
			$$ \frac{g(r)}{g(s)} \geq c_1 \Big(\frac{r}{s}\Big)^{\beta_1} \quad \text{for all} \;\; s\leq r< a\;\;(\text{resp.}\,\;a < s\leq r),$$
			and	we say that $\U_a(g,\beta_2, c_2)$ (resp. $\U^a(g,\beta_2, c_2)$) holds if
			$$ \frac{g(r)}{g(s)} \leq c_2 \Big(\frac{r}{s}\Big)^{\beta_2} \quad \text{for all} \quad s\leq r< a\;\;(\text{resp.}\;a < s\leq r).$$
			We say that  $\L(g,\beta_1,c_1)$ holds if  $\L_\infty(g,\beta_1,c_1)$ holds, and that $\U(g,\beta_2,c_2)$ holds if  $\U_\infty(g,\beta_2,c_2)$ holds.
	}\end{defn}

	Finally, we are ready to present our results in full generality which give precise answers for the question ${\bf (Q)}$.  Recall that we have assumed \eqref{e:phipsi} and that $\sM_+$ is the set of all positive non-decreasing functions defined on $(0,c_1] \cup [c_2,\infty)$ for some $c_2 \ge c_1>0$. 
	
	\begin{thm}\label{t:limsup0-1}
	{	Let $U\subset M$ be an open set. Suppose that $\U_{R_0}(\phi,\beta_2,C_U)$ and  $\U_{R_0}(\psi,\beta_3,C_U')$  hold for some $R_0>0$.  The following limsup laws at zero hold.}
	
			\smallskip
		
		\noindent (i) 
		Assume that \Eo \ and \Tailo \  hold.
		If 
		\begin{equation}\label{limsup0-case1}
			\limsup_{r \to 0} \frac{\psi(r)}{\phi(r)}<\infty,
		\end{equation}
		then for any $\Psi \in \sM_+$,
		\begin{equation}\label{3}
			\begin{split}
				&\int_0 \frac{dt}{\phi(\Psi(t))} = \infty \,\,({\rm resp.} < \infty) \\
				& \Leftrightarrow \;\; \limsup_{t \to 0} \frac{\sup_{0<s \le t}d(x,X_s)}{\Psi(t)} = \limsup_{t \to 0} \frac{d(x,X_t)}{\Psi(t)} = \infty \,\,({\rm resp.} =0), \quad  \P^x\text{-a.s.} \;\; \forall x \in U.
			\end{split}
		\end{equation}

		\noindent (ii)  Assume that \VRDo, \Cho, \NDLo \ and \Tailo \  hold. 
		If  
		$$\limsup_{r \to 0} \frac{\psi(r)}{\phi(r)} = \infty,$$ then  there exist  $\Psi \in \sM_+$  and constants $0<a_1 \le a_2<\infty$  such that for all $x \in U$, there exist $a_{x,1}, a_{x,2}\in [a_1, a_2]$ satisfying
		\begin{equation}\label{e:limsuplaw1}
			\limsup_{t \to 0} \frac{\sup_{0<s \le t}d(x,X_s)}{\Psi(t)} = a_{x,1} \quad \text{and} \quad 	\limsup_{t \to 0} \frac{d(x,X_t)}{\Psi(t)} = a_{x,2},\qquad
			\P^x\mbox{-a.s.}
		\end{equation}
	\end{thm}

	\begin{thm}\label{t:limsupinf-1}
		{ Suppose that $\L^{R_\infty}(\phi,\beta_1,C_L)$, $\U^{R_\infty}(\phi,\beta_2,C_U)$ and $\U^{R_\infty}(\psi,\beta_3,C_U')$ hold for some $R_\infty \ge 1$. The following limsup laws at infinity hold.}
		
		\smallskip
		
		\noindent (i) 
		Assume that \Ei \ and \Taili \ hold.
		If 
			\begin{equation}\label{limsupinf-case1}
				\limsup_{r \to \infty} \frac{\psi(r)}{\phi(r)}<\infty,
			\end{equation}
		then for any $\Psi\in \sM_+$, 
		\begin{equation}\label{4}
			\begin{split}
				&\int^\infty \frac{dt}{\phi(\Psi(t))} = \infty \,\,({\rm resp.} < \infty) \\
				& \Leftrightarrow \quad \limsup_{t \to \infty} \frac{\sup_{0<s \le t}d(x,X_s)}{\Psi(t)} = \limsup_{t \to \infty} \frac{d(x,X_t)}{\Psi(t)} = \infty \,\,({\rm resp.} =0), \quad  \P^y\text{-a.s.} \;\; \forall x,y \in M.
			\end{split}
		\end{equation}

		\noindent (ii) 
		Assume  that \VRDi, \Chi, \NDLi \ and \Taili \ hold. If 
				$$\limsup_{r \to \infty} \frac{\psi(r)}{\phi(r)}=\infty,$$  then there exist  $\Psi\in \sM_+$  and   constants $b_1, b_2 \in (0,\infty)$  such that
		\begin{equation}\label{e:limsuplaw2}
			\limsup_{t \to \infty} \frac{\sup_{0<s\le t}d(x,X_s)}{\Psi(t)} =b_1 \quad \text{and} \quad \limsup_{t \to \infty} \frac{d(x,X_t)}{\Psi(t)} = b_2, \quad
			\P^y\mbox{-a.s.} \;\; \forall x,y \in M.
		\end{equation}
	\end{thm}
	
	\smallskip

{l It is natural to seek an explicit form of the rate function $\Psi$ in Theorems \ref{t:limsup0-1}(ii) and \ref{t:limsupinf-1}(ii).
	When $X$ is a Brownian motion on certain Riemannian manifolds or  fractals, it has been proven in \cite{BP88, Ba98, BK} that \eqref{e:limsuplaw1} and \eqref{e:limsuplaw2} hold with $\Psi(t)= (t/\log|\log t|)^{1/d_w}\log |\log t|$ where $d_w \ge 2$ is the \textit{walk dimension} of the underlying space.
	See also \cite[Theorem 5.4]{BKKL19b} for Brownian-like jump processes in metric measure spaces. 
Given these cases, a natural choice for the function $\Psi$ in  \eqref{e:limsuplaw1} and \eqref{e:limsuplaw2} is  $\Psi_1(t):=\phi^{-1}(t/\log|\log t|) \log |\log t|$.
	In the following theorems, we provide integral tests that fully determine whether  $\Psi=\Psi_1$ can be used in equations \eqref{e:limsuplaw1} and \eqref{e:limsuplaw2},  respectively. If these integral tests fail, we observe that for any  $\Psi \in \sM_+$ satisfying \eqref{3} (or \eqref{4}), it holds that $\liminf \Psi(t)/\Psi_1(t)<\infty$  as $t$ tends to zero (or infinity), while  $\limsup \Psi(t)/\Psi_1(t)=\infty$. See \eqref{e:suppre1} and \eqref{e:suppre2} below.}
	
Note that when $\limsup_{r \to 0} \psi(r)/\phi(r)<\infty$, 
	 the integral \eqref{e:supint} below is always infinite. In fact, in such cases, one can see that the integrand in \eqref{e:supint} is greater than some constant multiple of $r^{-1}(\log |\log r|)^{-\beta_2}$. Similarly,  
	 if  $\limsup_{r \to \infty} \psi(r)/\phi(r)<\infty$, then the integral \eqref{e:supint2} below is always infinite.

	\begin{thm}\label{t:supprecise1}
{ Let $U\subset M$ be an open set.		Suppose that  \VRDo, \Cho,  \NDLo, \Tailol ,  $\U_{R_0}(\phi,\beta_2,C_U)$ and  $\U_{R_0}(\psi,\beta_3,C_U')$ hold for some $R_0>0$. Suppose also that   $\L_{R_0}(\phi, \beta_1, C_L)$ holds with $\beta_1 \ge 1$. Then the following statements hold.

}

			\smallskip
			
		\noindent(i) If the integral
		\begin{equation}\label{e:supint}
			\int_{0} \frac{\phi(r) \log |\log r|}{r \psi(r \log |\log r|)} dr 
		\end{equation}
		is finite, then there exist constants $0<a_1 \le a_2<\infty$ such that for all $x \in U$, there exists $a_x \in [a_1, a_2]$ satisfying
		\begin{equation}\label{e:1.14}
			\limsup_{t \to 0} \frac{d(x,X_t)}{\phi^{-1} ( t/\log |\log t| )\log |\log t|} = a_x,~\qquad\,
			\P^x\mbox{-a.s.}
		\end{equation}
		
		\noindent(ii) Assume further that \Tailog \ holds. If the integral \eqref{e:supint} is infinite, then
		for all $x\in U$, the left hand side of \eqref{e:1.14} is infinite, $\P^x$-a.s. Moreover, if we also have $\limsup_{r \to 0}\psi(r)/\phi(r)=\infty$, then any  $\Psi \in \sM_+$ satisfying \eqref{e:limsuplaw1} is factorized as 
		\begin{align}\label{e:suppre1}
			&\Psi(t) = f(t) \phi^{-1} ( t/\log |\log t|)\log |\log t|, \\[2pt]
			&\text{for a function $f$ such that }   \liminf_{t \to 0} f(t)< \infty, \; \limsup_{t \to 0}f(t)=\infty.
		\end{align}
	\end{thm}
	
	\begin{thm}\label{t:supprecise2}
	{	Suppose that \VRDi, \Chi, \NDLi, \Tailil, $\U^{R_\infty}(\phi,\beta_2,C_U)$ and \linebreak $\U^{R_\infty}(\psi,\beta_3,C_U')$ hold for some $R_\infty \ge 1$. Suppose also that  $\L^{R_\infty}(\phi,\beta_1,C_L)$ holds with  $\beta_1 \ge 1$. Then the following statements hold.}

		\smallskip

		\noindent(i)	If the integral
		\begin{equation}\label{e:supint2}
			\int^\infty \frac{\phi(r) \log \log r}{r \psi(r \log \log r)} dr 
		\end{equation}
		is finite, then there exists a constant $b_\infty \in (0, \infty)$ such that for all $x,y \in M$,
		\begin{equation}\label{e:1.15}
			\limsup_{t \to \infty} \frac{d(x,X_t)}{\phi^{-1} ( t/\log \log t)\log \log t} =  b_\infty,~\qquad\,
			\P^y\mbox{-a.s.,}
		\end{equation}
		
		\noindent (ii) Assume further that \Tailig \ holds. If the integral  \eqref{e:supint2} is infinite, then for all $x,y \in M$, the left hand side of \eqref{e:1.15} is infinite, $\P^y$-a.s.  Moreover,  if we also have $\limsup_{r \to \infty}\psi(r)/\phi(r)=\infty$, then any $\Psi \in \sM_+$ satisfying \eqref{e:limsuplaw2} is factorized as
		\begin{align}\label{e:suppre2}
			&\Psi(t) = f(t) \phi^{-1} ( t/\log \log t)\log \log t, \\[2pt]
			&\text{for a function $f$ such that }   \liminf_{t \to \infty} f(t)< \infty, \; \limsup_{t \to \infty}f(t)=\infty.
		\end{align}
	\end{thm}

	\begin{remark}\label{r:suppre}
		{\rm (i) \eqref{e:1.14} and \eqref{e:1.15} are also valid (with possibly different constants $a_x$ and $b_\infty$) with the numerator $\sup_{0<s \le t}d(x,X_s)$ instead of  $d(x,X_t)$.
			
			\noindent (ii)  The function $\Psi_1(t)=\phi^{-1}(t/\log|\log t|)\log|\log t|$ may not be non-decreasing near zero so that it may not belong to $\sM_+$. Nevertheless, under $\U_{R_0}(\phi, \beta_2, C_U)$, since $t/\log|\log t| \downarrow 0$ as $t \downarrow 0$, one can verify that $\Psi_1(t) \asymp \sup_{s \in (0, \, t \wedge (1/16)]} \Psi_1(s) =:\Psi_2(t)\in \sM_+$ for $t \in (0, 1/16]$. Hence, by the Blumenthal's zero-one law, \eqref{e:1.14} is still valid even if we use the function $\Psi_2$ as the denominator instead of $\Psi_1$. However, for brevity, we simply used the function $\Psi_1$ in \eqref{e:1.14}, instead of $\Psi_2$. 
		}
	\end{remark}

	To prove Theorems  \ref{t:limsupinf-1}(ii) and \ref{t:supprecise2}, we need  a proper zero-one law. Note  that   there is no assumption near zero in
these theorems. Hence, under the setting of  Theorem  \ref{t:limsupinf-1}(ii) or \ref{t:supprecise2},  it is not possible to prove the continuity for parabolic functions in $M$. However, under that setting, we can establish an oscillation result of  parabolic functions for large distances in Proposition \ref{p:phr}. Then using this result, we show that a zero-one law holds for shift-invariant events without assuming that $M$ is connected. 
	Cf. \cite[Theorem 8.4]{BB99}, \cite[Proposition 2.3]{BK} and \cite[Theorem 2.10]{KKW17}.

	\vspace{1mm}
	
	The rest of the paper is organized as follows. In Section \ref{s:dset}, we apply our main theorems to two classes of Markov processes: subordinate processes and symmetric Hunt processes. 
	The proofs of assertions in Section \ref{s:dset} are given in Appendices A and B.
	In Section \ref{s:RCM},  we apply our results to random conductance models.	
	
	In  Sections \ref{s:pre} and  \ref{s:limsup},  we prove our limsup LILs: 
	We first introduce auxiliary functions and use these to obtain precise bounds on  tail probability on the first exit times from balls. 
	After establishing the zero-one law in Section \ref{s:tp}, 
	we present the proofs of our main results in Section \ref{s:limsup}. 
	Proofs of Propositions \ref{p:EP}, \ref{p:EPL}, Theorems \ref{t:limsup0-1}(ii) and \ref{t:limsupinf-1}(ii) are the most delicate part of this paper.

	\medskip
	
	{\bf Notations}: 
We use same fixed  positive real constants $d_1$, $d_2$, $\beta_1$, $\beta_2$, $\beta_3$,  $\eta$,  $C_L$, $C_U$, $C_U'$ and $C_i$, $i=1,2,...$ on conditions and statements both at zero and at infinity. On the other hand,  lower case letters $\eps$, $\delta$, $c$,  $b_\infty$, $a_x$, $a_{x,i}$, $a_i$, $b_i$ and $c_i$, $i=0,1,2,...$ denote positive real constants  and are fixed in each statement and proof, and the labeling of these constants starts anew in each proof. 
	
	We use the symbol ``$:=$'' to denote a definition, 
	which is read as ``is defined to be.''  
	Recall that $a\wedge b:=\min\{a,b\}$ and  $a\vee b:=\max\{a,b\}$. We set $\lceil a \rceil := \min \{ n \in \Z : n \ge a \}$ and $\text{diag}:=\{(x,x):x \in M\}$, and denote  $\overline{A}$ for the closure of $A$.  The notation $f(x) \asymp g(x)$ means that there exist constants $c_2 \ge c_1>0$ such that $c_1g(x)\leq f (x)\leq c_2 g(x)$ for a specified range of $x$. We denote by $C_c(D)$ the space of all continuous functions  with compact support in $D$.

	\section{LILs for subordinate processes and symmetric jump processes}\label{s:dset}

	Recall that $(M,d)$ is a locally compact separable metric space, and $\mu$ is a positive Radon measure on $M$ with full support. 
	Set $\bar R:=\sup_{y,z \in M}d(y,z)$. 
	Let $F$ be an increasing continuous function on $(0,\infty)$ such that 
	\begin{equation}\label{e:Fscale}	\text{$\L(F, \gamma_1, c_L)$ and  $\U(F, \gamma_2, c_U)$   hold for some constants $\gamma_2 \ge \gamma_1>1$ and $c_L,c_U>0$.}
	\end{equation} 
{  We set $F(\infty):=\infty$.}
	
	\medskip
	{\it Throughout this section and Appendices \ref{s:A} and \ref{s:B}, we assume that
		$\mathrm{VRD}_{\bar R}(M)$ and  $\mathrm{Ch}_{\bar R}(M)$  hold.		 We  also  assume that
		 there is a  conservative Hunt  process  $Z = (\Omega, \FF, \FF_t, Z_t,  t \ge 0;\, 
		\P^x, x \in M)$ on $M$ which has a heat kernel $q(t,x,y)$ (with respect to $\mu$) enjoying the following estimates:	There exist constants $R_1 \in (0, \infty]$ and  $c_1,c_2,c_3>0$ such that for all $t \in (0, F(R_1))$ and $x,y \in M$,
		\begin{equation}\label{e:diffusion}
			\frac{c_1}{V(x, F^{-1}(t))}\1_{\{F(d(x,y)) \le t\}}\le  q(t,x,y) \le \frac{c_2}{V(x, F^{-1}(t))} \exp \big(- c_3  F_1(d(x,y),t) \big),
		\end{equation}
		where the function $F_1$ is defined as
		\begin{equation}\label{e:F1}
			F_1(r,t):=\sup_{s>0} \bigg( \frac{r}{s}-\frac{t}{F(s)}\bigg).
		\end{equation}
	}
	
	The above function $F_1$ is widely used in heat kernel estimates for diffusions on metric measure spaces including the Sierpinski gasket or carpet, nested fractals and affine nested fractals. See  \cite{BP88, 
		 BB99, GT12,  HK99
	 }. Note that  if $F(r)=r^\gamma$ for $\gamma>1$, then $F_1(r,t) = c\gamma (r^{\gamma}/t)^{1/(\gamma-1)}$ for some $c_\gamma>0$.

	We mention that the process $Z$  may not be $\mu$-symmetric. For example, $Z$ can be a Brownian motion with drift on $\R^d$, which has $\Delta + \mathsf{p} \cdot \nabla$ as the infinitesimal generator where the function $\mathsf{p}$ belongs to some suitable Kato class  (see, e.g. \cite{Zh97, KS06}).   In this case, $\mu$ is the Lebesgue measure so $V(x,r)=\omega_dr^d$ for a constant $\omega_d>0$ and $F(r)=r^2$.

	\subsection{General subordinate processes}

{	A non-negative L\'evy prcoess  is called a \textit{subordinator}.}
	Let $S=(S_t)_{t\ge0}$ be a subordinator  on the probability space $(\Omega, \FF, \P)$ which is independent of  $Z$.  
	{ The Laplace transform of $S$ can be expressed as
\begin{equation*}
	 \E\big[e^{-\lambda S_t}\big]= e^{-t\phi_1(\lambda)}, \quad  t,\lambda \ge 0,
\end{equation*}
where the function $\phi_1$ is referred to as \textit{the Laplace exponent} of $S$. It is known that there exist a unique non-negative constant $b$ and a unique measure $\nu$ on $(0,\infty)$  satisfying $\int_{(0,\infty)}(1 \land u) \nu(du)<\infty$ such that 
	\begin{equation*}
	\phi_1(\lambda)=b\lambda + \int_{(0,\infty)}(1-e^{-\lambda u})\nu(du), \quad \lambda \ge 0.
	\end{equation*}
The measure $\nu$ is called \textit{the L\'evy measure} of $S$. In this subsection and Appendix \ref{s:A}, we always assume that $\lim_{\lambda \to \infty} \phi_1(\lambda)=\infty$, or equivalently, }
	\begin{equation}\label{e:not-compound}
		\textit{either } \;\; b \neq 0 \;  \text{ or } \; \nu((0,\infty)) = \infty.
	\end{equation}

Let $X_t:=Z_{S_t}$, which is called a subordinate process.  The subordinate process $X$ is a {Hunt} process on $M$ and it has the heat kernel $p(t,x,y)$ and the jumping kernel $J(x,y)\mu(dy)$ given by {(see \cite[p. 67, 73--75]{Bo84})}
	\begin{equation}\label{e:PJ}
		p(t,x,y) = \int_0^\infty q(s,x,y) \P(S_t \in ds) \quad\;\; \text{and} \quad\;\; J(x,y) = \int_0^\infty q(s,x,y) \nu(ds).
	\end{equation} In order to analyze the limsup LILs for the process $X$, we introduce a set of auxiliary functions. (See \cite{JP87}). Let $\phi_0$ be a pure jump part of $\phi_1$, namely,
	$$
	\phi_0(\lambda):=\int_{(0,\infty)}(1-e^{-\lambda u})\nu(du),
	$$
and let $H_0(\lambda):=\phi_0(\lambda)-\lambda \phi_0'(\lambda)$. Note that  $H_0$ is non-negative and increasing. Denote by $w(r)=\nu((r,\infty))$ the tail of the L\'evy measure of $S$, and define
	\begin{equation}\label{e:subaux1}
		\Phi(r):= \frac{1}{\phi_1(F(r)^{-1})}, \quad \;\; \Theta(r):=\frac{1}{H_0(F(r)^{-1})} \quad   \text{ and } \quad \Pi(r):= \frac{2e}{w(F(r))} \quad \text{for } \; r>0.
	\end{equation} 
	Note that $\Phi, \Theta$ and $\Pi$ are non-decreasing. Moreover,  $\lim_{r \to 0}\Phi(r)  =0$ by \eqref{e:not-compound}. For all $r>0$ and $\lambda \ge 1$, we have $\phi_1(\lambda r) \le \lambda \phi_1(r)$ and $H_0(\lambda r) \le \lambda^2 H_0(r)$.
	See, e.g. \cite[Lemma 2.1(a)]{Mi16}. 
	Thus, by \eqref{e:Fscale},  we see that 
	\begin{equation}\label{e:aux}
	\text{$\U(\Phi, \gamma_2, c_U)\;$  and  $\;\U(\Theta, 2\gamma_2, c_U)\;$ hold.}
	\end{equation}
	Further,  by \cite[Lemma 2.6]{Mi16},  we also see that
	\begin{equation}\label{e:aux2}
		\Phi(r)  \le \phi_0(F(r)^{-1})^{-1}  \le \Theta(r)  \le \Pi(r) \quad\;\; \text{for all} \;\; r>0.
	\end{equation}
	
	Now, we are ready to state the limsup LILs for $X$. The proofs of the next two theorems are given in Appendix \ref{s:A}.

	\begin{thm}\label{t:SDsup0} {\bf  (Limsup LILs at zero)}  (i) Suppose that 
		$$\limsup_{r \to \infty}\frac{\phi_1(r)}{H_0(r)}<\infty \quad (\text{or, equivalently,  $\limsup_{r \to \infty}\frac{\phi_1(r)}{w(r^{-1})}<\infty$}).$$
		Then for every $\Psi \in \sM_+$ and $x \in M$, \eqref{3} holds with $\phi=\Phi$.
		
		\smallskip
		
		\noindent (ii) Suppose that 
			$$\limsup_{r \to \infty}\frac{\phi_1(r)}{H_0(r)}=\infty \quad (\text{or, equivalently,  $\limsup_{r \to \infty}\frac{\phi_1(r)}{w(r^{-1})}=\infty$}).$$ Then the following two statements are true.

		(a) Assume that $\L^1(\phi_1, \alpha_1, c_1)$  holds with $\alpha_1 \ge 1/\gamma_1$. If the integral \eqref{e:supint} is finite with $\phi=\Phi$ and $\psi=\Theta$, then there exist constants $a_2 \ge a_1>0$ such that for all $x \in M$, there exists $a_x \in [a_1,a_2]$ satisfying \eqref{e:1.14} with $\phi=\Phi$.

		(b) Assume that $\L^1(\phi_1, \alpha_1, c_1)$ and  $\U_{F(R_2)}(w^{-1}, \alpha_2, c_2)$  hold with some $R_2>0$. 
		Then we can find a function $\Psi \in \sM_+$ and  constants $a_2 \ge a_1>0$  such that
		\eqref{e:limsuplaw1} holds  for all $x \in M$. 
		Moreover, if  $\alpha_1  \ge 1/\gamma_1$, then Theorem \ref{t:supprecise1} holds with $\phi=\Phi$ and $\psi=\Pi$.
		
	\end{thm}
	
	\begin{thm}\label{t:SDsup2} {\bf  (Limsup LILs at infinity)}
		Assume that $M$ is unbounded and $R_1 = \infty$.
		
		\noindent (i) Suppose that 	$$\limsup_{r \to 0}\frac{\phi_1(r)}{H_0(r)}<\infty \quad (\text{or, equivalently,  $\limsup_{r \to 0}\frac{\phi_1(r)}{w(r^{-1})}<\infty$}).$$
		Then for every  $\Psi \in \sM_+$, \eqref{4} holds with $\phi=\Phi$.
		
		\smallskip
		
		\noindent (ii)
	Suppose that	$$\limsup_{r \to 0}\frac{\phi_1(r)}{H_0(r)}=\infty \quad (\text{or, equivalently,  $\limsup_{r \to 0}\frac{\phi_1(r)}{w(r^{-1})}=\infty$}).$$
		Then the following two statements are true.

		(a) Assume that $\L_1(\phi_1, \alpha_3, c_1)$ holds with $\alpha_3 \ge 1/\gamma_1$. If the integral \eqref{e:supint2} is finite with $\phi=\Phi$ and $\psi=\Theta$, then there exists a constant $b_\infty \in (0,\infty)$ such that  \eqref{e:1.15} holds with $\phi=\Phi$ for all $x,y \in M$.

		(b) Assume that $\L_1(\phi_1, \alpha_3, c_1)$  and  $\U^{F(R_3)}(w^{-1}, \alpha_4, c_2)$ hold with some $R_3>0$. 
		Then we can find a function $\Psi \in \sM_+$  and   constants $b_1, b_2\in (0,\infty)$  such that  \eqref{e:limsuplaw2} holds for all $x,y \in M$.
		Moreover, if  $\alpha_3 \ge 1/\gamma_1$, then Theorem \ref{t:supprecise2} holds with $\phi=\Phi$ and $\psi=\Pi$.
	\end{thm}

	\begin{remark}
		{\rm  (i) According to \cite[Remark 1.3(1)]{CK20}, for all $\alpha>0$, we have
			\begin{equation}\label{e:assonw}
				\text{$\U_{F(R_2)}(w^{-1}, \alpha, c_1)$ (resp. $\U^{F(R_3)}(w^{-1}, \alpha, c_1)$) } \text{ $\Rightarrow$  \, $\U^1(\phi_0, \alpha \land 1, c_2)$ (resp. $\U_{1}(\phi_0, \alpha \land 1, c_2)$).}
			\end{equation}
			Moreover,  if $\alpha<1$, then there exists $r>0$ such that 
			\begin{equation*}
			 \text{ $\U^1(\phi_0, \alpha, c_4)$ (resp. $\U_{1}(\phi_0, \alpha, c_4)$) } 	\text{ $\Rightarrow$  \, $\U_{F(r)}(w^{-1}, \alpha, c_3)$ (resp. $\U^{F(r)}(w^{-1}, \alpha, c_3)$).}
			\end{equation*}
			We note that the constant $\alpha$ in  \eqref{e:assonw} can be larger than $1$. By imposing a weak scaling property on $w$ instead of $\phi_0$, our results cover 
			not only mixed $\alpha$-stable-like subordinators but also
			 a large class of subordinators whose Laplace exponent is regularly varying at infinity of index $1$.
		}
	\end{remark}

	Below, we give some concrete examples for the	subordinate process $X_t=Z_{S_t}$. In the remainder of this subsection, we assume that $F(r)=r^\gamma$ for some $\gamma>1$.

	\begin{example}\label{e:mixedtype} {\bf (Nonzero diffusion term)}
		{\rm Assume that $\phi_1(\lambda)=b\lambda + \phi_0(\lambda)$ for $b>0$. Then   $\phi_1(r^{-1}) \asymp br^{-1}$ for $r \in (0,1)$ and  $\limsup_{r \to 0}\phi_1(r^{-1})/w(r) \ge b(2e)^{-1}\lim_{r \to 0}r^{-1}/\phi_0(r^{-1}) = \infty$.	Moreover, by the change of the variables $u = r^{-\gamma} (\log |\log r|)^{-\gamma}$ and Tonelli's theorem, we get
			\begin{align}
				&\int_{0} \frac{\Phi(r) \log|\log r|}{r \Theta(r \log |\log r|)}dr = \int_{0} \frac{ H_0(r^{-\gamma}(\log|\log r|)^{-\gamma})\log|\log r|}{r\phi_1(r^{-\gamma})}dr \\
				&\le c_1 \int_{0} r^{\gamma-1}H_0(r^{-\gamma}(\log|\log r|)^{-\gamma})\log|\log r|dr \le c_2\int^{\infty} u^{-2}H_0(u) (\log \log u)^{1-\gamma}  du \label{e:cal2.5}\\
				&= c_2\int_0^\infty \int^{\infty}u^{-2}(1-e^{-us}-use^{-us}) (\log \log u)^{1-\gamma} du\,\nu(ds)  \\
				& \le c_2\int_0^\infty \int^{\infty}u^{-2}((us)^2 \land 1)  du\,\nu(ds) \le c_3\int_0^\infty (s \wedge 1)\nu(ds) <\infty.
			\end{align}
			In the third inequality above, we used the fact that  $1-e^{-\lambda} -\lambda e^{-\lambda} \le \lambda^2 \land 1$ for all $\lambda \ge 0$. Therefore,  according to Theorem  \ref{t:SDsup0}(ii-b), 
			there exist constants $a_2 \ge a_1>0$ such that for all $x \in M$, there is a constant $a_{x} \in [a_1,a_2]$ satisfying
			\begin{equation}\label{e:ex2.5}
				\limsup_{t \to 0} \frac{d(x, X_t)}{b^{1/\gamma}t^{1/\gamma}(\log |\log t|)^{1-1/\gamma}}=a_{x},\quad \P^x\text{-a.s.}
			\end{equation}	
			In particular, by taking $b=1$ and $\phi_0\equiv 0$, we see that the process $Z_t$ satisfies  \eqref{e:ex2.5}.
			We remark here  that limsup LIL of $d(x,Z_t)$ at zero is  studied in \cite[Theorem  4.7]{BP88} when $Z$ is $\mu$-symmetric.	\qed
		}
	\end{example}

	\begin{example}\label{e:geo} {\bf (Low intensity of small jumps)}
		{\rm Let $\phi_1(\lambda)=(\log (1+ \lambda^\alpha))^{\delta}$ for $\alpha, \delta \in (0,1]$, $\alpha\delta<1$. Here, $\phi_1$ is indeed the Laplace exponent of a subordinator by \cite[Theorem 5.2, Proposition 7.13 and Example 16.4.26]{SRV12}. A prototype of such $X_t=Z_{S_t}$ is a geometric $2\alpha$-stable process on $\R^d$ ($0<\alpha<1$), that is, a L\'evy process on $\R^d$ with the characteristic exponent $\log (1+ |\xi|^{2\alpha})$.
			
			In view of \cite[Lemma 2.1(iii) and (iv)]{CK20}, we have
			\begin{equation*}
				\big(\log(1+r^{-\alpha})\big)^{\delta} = 	\phi_1(r^{-1}) \asymp w(r) \asymp \begin{cases}
					|\log r|^{\delta}, \; &\mbox{for} \;\; r\in (0, e^{-1}),\\
					r^{-\alpha \delta}, \; &\mbox{for} \;\; r\in (e, \infty).
				\end{cases}
			\end{equation*}
			Thus, using Theorems \ref{t:SDsup0}(i) and \ref{t:SDsup2}(i), we get that for all $x \in M$, the limsup LIL at zero \eqref{3} holds with $\phi(r) = |\log r|^{-\delta}$ and the limsup LIL at infinity \eqref{4} holds with $\phi(r)=r^{\gamma \alpha \delta}$.
			
			\medskip
			
			Let $T_t$ be a $\alpha\delta$-stable subordinator with the Laplace exponent $\lambda^{\alpha \delta}$.
			By Theorem \ref{t:SDsup2}(i), the subordinate process $Z_{T_t}$ also enjoys the  limsup LIL at infinity \eqref{4} with $\phi(r)=r^{\gamma \alpha \delta}$. 
			However, we show in the below that if $\delta<1$, then estimates on the heat kernel $p(t,x,y)$ of $X_t$ are different from estimates on the heat kernel $\wt p(t,x,y)$ of $Z_{T_t}$, even for large $t$. Precisely, we will see that $p(t,x,x) = \infty$ for all $t>0$ and $x \in M$, but $\wt p(t,x,y) \le C/V(x, t^{1/(\gamma \alpha \delta)})$ for all $t>0$ and $x,y \in M$ for a constant $C>0$.

			Indeed, we can find the above estimates on $\wt p(t,x,y)$ from \cite[Example 6.1]{CKW16a}. On the other hand, observe that $\phi_1'(\lambda) 
			\asymp \lambda^{-1}(\log \lambda)^{\delta-1}$ for all $\lambda>e^2$ so that $(\phi_1')^{-1}(r) \asymp r^{-1} |\log r|^{\delta-1}$ and  $(H_0 \circ (\phi_1')^{-1})(r) \asymp |\log r|^\delta$ for $r \in (0, e^{-3})$. Thus, by  \eqref{e:ndlq} in Appendix,  the integration by parts and \cite[Lemma 5.2(ii)]{JP87}, since $\mathrm{VRD}_{\bar R}(M)$ holds, we get that for all $t>0$ and $x\in M$, 
			\begin{equation*}
				\begin{split}
					p(t,x,x) &\ge \int_0^\infty q^{B(x,1)}(s,x,x) \P(S_t \in ds) \ge  c_6\limsup_{\eps \to 0} \int_{0}^\eps V(x, s^{1/\gamma})^{-1} \P(S_t \in ds) \\
					&\ge c_7\limsup_{\eps \to 0}  \eps^{-c_8} \P(S_t \le \eps) \ge c_9 \limsup_{\eps \to 0}  \exp\big( c_8|\log \eps|-c_{10} t |\log (\eps/t)|^\delta  \big) = \infty.
				\end{split}
			\end{equation*}
			
			The above example shows that the condition $\mathrm{NDL}^{R}(r^{\gamma \alpha \delta})$ for $X$ may not be  sharp (which holds for some $R\ge1$ according to Lemma \ref{l:subNDL}(ii) in Appendix). However, it is sufficient for applying our theorems to obtain the LIL for $X$.

			We mention that our results also cover the cases when $\phi_1(\lambda)=\phi_0(\lambda) = ({{Log}} \circ...\circ {Log})(\lambda)$, where $Log(\lambda):=\log(1+\lambda)$, for arbitrary finite number of compositions. (These functions are the Laplace exponent of a subordinator according to \cite[Theorem 5.2, Corollary 7.9 and Example 16.4.26]{SRV12}.) \qed
			
		}
	\end{example}

	\begin{example}\label{e:polydecay} {\bf (Jump processes with mixed polynomial growths)}
		{\rm In this example, we work with pure jump processes  enjoying the limsup LIL of type \eqref{e:limsuplaw1} or \eqref{e:limsuplaw2}. Assume $b=0$ so that $\phi_1=\phi_0$. Then by \cite[Lemma 2.1]{CK20}, 		 we have	\begin{equation}\label{e:phiw}
				\Phi(r)= \phi_0(r^{-\gamma})^{-1} \asymp \frac{2er^\gamma}{ \int_0^{r^\gamma}w(s)ds} = \frac{r^\gamma}{\gamma \int_0^{r}s^{\gamma-1}\Pi(s)^{-1}ds}  \quad\;\; \text{for} \;\; r>0,
			\end{equation}
			where comparison constants are independent of the subordinator $S$.
			
			\smallskip
			
			\noindent (i) \textit{limsup LIL at zero.}
			Here, we give examples of $w$ which make the integral \eqref{e:supint}  be finite or infinite with $\phi=\Phi$ and $\psi=\Pi$. Based on that, we study a small time limsup LIL for $X$.
			
			\smallskip

			(a) Let $p>1$ and assume that 
			$$
			w(r) \asymp r^{-1}|\log r|^{-p} \quad \text{for} \;\, r \in (0,e^{-1}).
			$$
			 Here, $p$ must be larger than $1$ because it should holds that $\int_0^1 w(s)ds<\infty$.  By \eqref{e:phiw}, we  have that for $r, t \in (0, e^{-4})$,
			\begin{equation}\label{eg1}
				\Pi(r) \asymp r^\gamma |\log r|^{p} \quad  \text{and} \quad    \Phi(r) \asymp r^\gamma|\log r|^{p-1}
			\end{equation}
			and
			\begin{equation}\label{eg11}
				\Phi^{-1}(t/\log|\log t|)\log|\log t| \asymp t^{1/\gamma} |\log t|^{(1-p)/\gamma} (\log|\log t|)^{(\gamma-1)/\gamma}=:\Psi_1(t).
			\end{equation}
			Since $\lim_{r \to 0} \Pi(r)/\Phi(r)=\infty$, by Theorem \ref{t:SDsup0}(ii-a) and \eqref{e:subaux1},  there exists $\Psi \in \sM_+$ satisfying \eqref{e:limsuplaw1}. Moreover, one can see from \eqref{eg1} that the integral  \eqref{e:supint} is finite  if $\gamma>2$, and is infinite if $\gamma \le 2$. Hence, by Theorem \ref{t:SDsup0}(ii-a) and \eqref{eg11},  if $\gamma >2$, then  there exist constants $a_2 \ge a_1>0$ such that for all $x \in M$, there is $a_x \in [a_1, a_2]$ satisfying
			\begin{equation}\label{e:mixedex1}
				\limsup_{t \to 0}\frac{d(x,X_t)}{t^{1/\gamma} |\log t|^{(1-p)/\gamma} (\log|\log t|)^{(\gamma-1)/\gamma}}= a_{x},\qquad
				\P^x\mbox{-a.s.},
			\end{equation}
			and if $\gamma \le 2$, then for all $x \in M$, $\P^x$-a.s., the left hand side of \eqref{e:mixedex1} is infinite. Thus, for example, if $M$ is a Sierpinski gasket or carpet, then the left hand side of \eqref{e:mixedex1} is a deterministic constant a.s. while it is infinite a.s. if $M$ is a Euclidean space $\R^d$, $d \ge 1$. 
			
			Now, we give an example of a proper rate function $\Psi$ satisfying \eqref{e:limsuplaw1} in case of $\gamma \le 2$ (see \eqref{e:defsn} and \eqref{e:defPsi} in Section \ref{s:limsup}): Choose any $\delta>0$ and define
			\begin{align*}
				&\Psi(t):= \sum_{n=4}^\infty (\log n) \exp\big(-\gamma^{-1}n (\log n)^{2-\gamma + \delta}\big) \cdot \1_{(t_{n+1}, t_n]}(t) \; \asymp \; \sum_{n=4}^\infty \Psi_1(t_n)\cdot \1_{(t_{n+1}, t_n]}(t) \\
				&\text{where} \quad  t_n:=n^{p-1}(\log n)^{(2-\gamma + \delta)p +\gamma-\delta-1}\exp\big(- n (\log n)^{2-\gamma + \delta}\big).
			\end{align*}
			For the above $\Psi$ and $t_n$,  one can see that  since $2-\gamma + \delta>0$, 
			$$\lim_{n \to \infty} \Psi(2t_{n+1})/\Psi_{1}(2t_{n+1}) \asymp \lim_{n \to \infty} \Psi_1(t_{n})/\Psi_1(t_{n+1}) \ge c\lim_{n \to \infty} (t_n/t_{n+1})^{1/(2\gamma)}= \infty.$$ Thus, the above $\Psi$ can be factorized as \eqref{e:suppre1} with $\phi = \Phi$.

			\smallskip

			(b) Let $p>1$ and assume that 
			$$w(r) \asymp r^{-1}|\log r|^{-1}(\log|\log r|)^{-p} \quad \text{for} \;\,r \in (0, e^{-4}).$$ (As in (a), $p$ must be larger than $1$ since $\int_0^1 w(s)ds<\infty$.) In this case, by \eqref{e:phiw}, we see that
			\begin{equation}\label{eg2}
				\Pi(r) \asymp r^\gamma |\log r| (\log |\log r|)^{p} \quad \text{and} \quad  \Phi(r) \asymp r^{\gamma} (\log|\log r|)^{p-1}  \;\; \text{for} \;\; r \in (0, e^{-4}).
			\end{equation}
			Since $\lim_{r \to 0} \Pi(r)/\Phi(r)$, as in (a), there exists $\Psi\in \sM_+$ which satisfies \eqref{e:limsuplaw1}.
			Moreover, we see from \eqref{eg2} that the integral \eqref{e:supint} is always finite.  According to Theorem \ref{t:SDsup0}(ii-a), it follows that  there exist constants $a_2 \ge a_1>0$ such that for all $x \in M$, there is $a_x \in [a_1, a_2]$ satisfying
			\begin{equation}\label{e:mixedex11}
				\limsup_{t \to 0}\frac{d(x,X_t)}{ t^{1/\gamma} (\log|\log t|)^{(\gamma-p)/\gamma}}= a_{x},\qquad
				\P^x\mbox{-a.s.}
			\end{equation}

			\smallskip
			
			\noindent (ii) \textit{limsup LIL at infinity.}  Let $R_1 = \infty$, $k,q>0$ and  assume that 
			$$w(r) \asymp r^{-1}(\log r)^{-1}\exp\big(k (\log \log r)^{q}\big)\quad\text{for}\;\, r \in (e^4, \infty).$$ Then  for $r,t$ large enough,
			\begin{equation}\label{eg3}
				\Pi(r) \asymp r^{\gamma}(\log r)\exp\big(-k (\log \log r^\gamma)^{q}\big) \;\; \text{and} \;\;    \Phi(r) \asymp r^{\gamma} (\log\log r)^{q-1}\exp\big(-k (\log \log r^\gamma )^{q}\big)
			\end{equation}
			and
			\begin{equation}\label{eg31}
				\Phi^{-1}(t/\log\log t)\log \log t \asymp t^{1/\gamma} (\log \log t)^{(\gamma-q)/\gamma} \exp\big(\gamma^{-1}k (\log \log t)^{q}\big)=:\Psi_\infty(t).
			\end{equation}
			Since $\lim_{r \to \infty} \Pi(r)/\Phi(r) = \infty$, by Theorem \ref{t:SDsup2}(ii-a) and \eqref{e:subaux1},  there exists $\Psi\in \sM_+$ satisfying \eqref{e:limsuplaw2}. Besides, in view of \eqref{eg3}, one can check that  the integral \eqref{e:supint2} is finite if $q \in (0, \gamma-1)$, and is infinite if $q \ge \gamma-1$. By Theorem \ref{t:SDsup2}(ii-a) and \eqref{eg31}, it follows that if $q \in (0, \gamma-1)$, then there exists a constant $b_\infty \in (0,\infty)$ such that
			\begin{equation*}
				\limsup_{t \to \infty}\frac{d(x,X_t)}{t^{1/\gamma} (\log \log t)^{(\gamma-q)/\gamma} \exp\big(\gamma^{-1}k (\log \log t)^{q}\big)}= b_\infty \quad \text{for all} \; x,y \in M, \;
				\P^y\mbox{-a.s.,}
			\end{equation*}
			and if $q \ge \gamma-1$, then  the left hand side of the above equality is infinite for all $x,y \in M$, $\P^y$-a.s. Moreover, in case of $q \ge \gamma-1$,  one can find an example of a proper rate function $\Psi$ satisfying \eqref{e:limsuplaw2}: Choose any $\delta>0$ and define 
			\begin{align}\label{e:largeex}
				&\Psi(t):= \sum_{n=5}^\infty (\log n) \exp\big(\, \gamma^{-1}n (\log n)^{q-\gamma +1 +  \delta}\,\big) \cdot \1_{(t_{n-1}, t_n]}(t) \; \asymp \;  \sum_{n=5}^\infty \Psi_\infty(t_n) \cdot \1_{(t_{n-1}, t_n]}(t)\quad \\
				&\text{where} \quad  t_n:=(\log n)^{q}\exp\big(\, n (\log n)^{q-\gamma +1+ \delta} -k \big(\log n + (q-\gamma + 1 +\delta)\log \log n \big)^q\,\big).
			\end{align}
			With
			the above $\Psi$ and $t_n$, one can see that  $\lim_{n \to \infty} \Psi(2t_{n-1})/\Psi_\infty(2t_{n-1})=\infty$ so that $\Psi$ admits a factorization \eqref{e:suppre2} with $\phi=\Phi$.
			
			By putting $q=1$ and $k=1-p$, we  obtain the following analogous result to \eqref{e:mixedex1}: Let $p<1$ and suppose that $w(r) \asymp r^{-1}(\log r)^{-p}$ for $r \in (e^4, \infty)$. If $\gamma > 2$, then there exists a constant $b_\infty' \in (0,\infty)$ such that 
			\begin{equation}\label{e:mixedex2}
				\limsup_{t \to \infty}\frac{d(x,X_t)}{t^{1/\gamma} (\log t)^{(1-p)/\gamma} (\log\log t)^{(\gamma-1)/\gamma}}= b_{\infty}' \quad \text{for all} \; x,y \in M, \;
				\P^y\mbox{-a.s.,}
			\end{equation}
			and if $\gamma \le 2$, then the  left hand side in \eqref{e:mixedex2} is infinite for all $x,y \in M$, $\P^y$-a.s. and an example of  $\Psi$ satisfying \eqref{e:limsuplaw2} is given as \eqref{e:largeex} (with $q=1$, $k=1-p$ and any $\delta>0$).\qed
		}
	\end{example}

	\subsection{Symmetric Hunt processes}\label{s:Hunt}
	
	Recall that we have assumed that the conditions $\mathrm{VRD}_{\bar R}(M)$  and $\mathrm{Ch}_{\bar R}(M)$ hold, and there exists a conservative Hunt process $Z$ in $(M,d,\mu)$ having a heat kernel $q(t,x,y)$ which satisfies \eqref{e:diffusion}. 
	In this subsection, we also assume that  $Z$ is $\mu$-symmetric  and associated with a \textit{regular strongly local Dirichlet form}  $(\EE^Z, \FF^Z)$ on $L^2(M;\mu)$. See \cite{FOT11} for the definitions  of a regular Dirichlet form and the strongly local property. Then by \cite[Section 3]{CKS87}, for any bounded $f \in \FF^Z$, there exists a unique positive Radon measure $\Gamma^Z(f,f)$ on $M$ such that $\int_M g d \Gamma^Z(f,f) = \EE^Z(f,fg) - \frac{1}{2}\EE^Z(f^2,g)$ for every $g \in \FF^Z \cap C_c(M)$.  One can see that the measure $\Gamma^Z(f,f)$ can be uniquely extended to any $f \in \FF^Z$ as the non-decreasing limit of $\Gamma^Z(f_n,f_n)$, where $f_n:=((-n) \vee f) \wedge n$. The measure $\Gamma^Z(f,f)$ is called the \textit{energy measure}  (which is also called the \textit{carr\'e du champ}) of $f$ for  $\EE^Z$. 
	
	\smallskip
	
	Let $(\EE^X, \FF^X)$ be a regular  Dirichlet form on $L^2(M;\mu)$  having the following expression: There exist a constant $\ll \ge 0$ and a 
	symmetric
	Radon measure of the form
	$\bar  J_X(dx, dy)=
	J_X(x,dy)\mu(dx)$ on $M \times M \setminus \text{diag}$ such that
	\begin{align}
		\EE^X(f,g)&=\ll\EE^{X,(c)}(f,g) + \int_{M \times M \setminus \text{diag}} (f(x)-f(y))(g(x)-g(y))J_X(x,dy)\mu(dx), \;\; f,g \in \FF^X \hspace{-1.5mm}, \nn\\
		\FF^X&=\overline{\big\{f \in C_c(M): \EE^X(f,f)<\infty\big\}}^{\EE^X_1},\label{e:Form}
	\end{align}
	where $\EE^{X,(c)}(f,g)$ is the strongly local part of $(\EE^X,\FF^X)$ (that is, $\EE^{X,(c)}(f,g)=0$ for all $f, g \in \FF^X$ such that $(f-c)g=0$ $\mu$-a.e. on $M$ for some $c \in \R$) and $\EE^X_1(f,f):=\EE^X(f,f) + \lVert f\rVert_{L^2(M;\mu)}^2$. See \cite[Theorem 3.2.1]{FOT11} for a general representation theorem on regular Dirichlet forms.  
	
	Since $(\EE^X,\FF^X)$ is regular, there exists an associated $\mu$-symmetric Hunt process $X=(X_t, t \ge 0; \P^x, x \in M \setminus \sN)$ where $\sN$ is a properly exceptional set in the sense that $\sN$ is nearly Borel,  $\mu(\sN)=0$ and $M_\partial \setminus \sN$ is $X$-invariant. This Hunt process is unique up to a properly exceptional set. See, e.g. \cite[Chapter 7 and Theorem 4.2.8]{FOT11}. We fix $X$ and $\sN$, and write $M_0:=M \setminus \sN$.

Our assumption on $X$ is the following. We use a convention $B(x,\infty)=M$ and  write $\Gamma^{X,(c)}(f,f)$ for the energy measure of $f \in \FF^X$ with respect to $\EE^{X,(c)}$. 
	
	\medskip

	\noindent {\bf Assumption L.} There exist an open set $\sU \subset M$,  constants $R_4 \in (0, \infty]$, $\kappa_0 \in (0,1)$ and an increasing function $\psi_0$ on $(0,R_4)$ such that the followings hold.
	
	\smallskip
	
 (L1)  $\L_{R_4}(\psi_0, \beta_1, c_1)$ and $\U_{R_4}(\psi_0,  \beta_2, c_2)$ hold with $\beta_2 \ge \beta_1>0$ and $\int_0 \psi_0(s)^{-1} dF(s)<\infty$.

(L2) There exist a kernel $J_X(\cdot, \cdot)$ on $\sU \times M$ and $c_3 \ge 1$   such that for all $x \in \sU$,
	$
	J_X(x,dy)=J_X(x,y)\mu(dy)$ in $B\big(x, R_4 \land (\kappa_0 \updelta_\sU(x))\big)$ and	
	\begin{equation*}
		\frac{c_3^{-1}}{V(x, d(x,y)) \psi_0(d(x,y))} \le J_X(x,y) \le \frac{c_3}{V(x, d(x,y)) \psi_0(d(x,y))}, \quad\;  y \in B\big(x, R_4 \land (\kappa_0 \updelta_\sU(x))\big).
	\end{equation*}

 (L3) There exists  $c_4>0$ such that for all $x \in \sU$, 
	$$J_X\left(x,M_\partial \setminus B\big(x,R_4 \land (\kappa_0\updelta_\sU(x))\big)\right) \le \frac{c_4}{\psi_0(R_4 \land \updelta_\sU(x))}.$$

 (L4) If $\ll>0$, then $\FF^X\cap C_c(\, \sU) = \FF^Z \cap C_c(\,\sU)$	and $d\Gamma^{X,(c)}(f,f)|_\sU \asymp d\Gamma^{Z}(f,f)|_\sU$ for  $f \in \FF^X$.

	\bigskip

	Under Assumption L, we define (cf. \eqref{e:phiw} and  \cite[(2.20)]{BKKL19b})
	\begin{equation}\label{d:HuntPhi}
		\Phi_1(r):=\frac{F(r)}{\int_0^r \psi_0(s)^{-1}dF(s) + \ll}, \quad \;\;  r \in (0, R_4),
	\end{equation}
	which is well-defined due to (L1). Note that $\Phi_1(r) \le  \psi_0(r)F(r)/\int_0^r dF(s)=\psi_0(r)$ for $r \in (0, R_4)$. 
	
		Now, we give the  limsup LIL for  $X$. The proof is given in Appendix \ref{s:B}.
	
	\begin{thm}\label{t:Hunt0} Suppose that  Assumption L  holds.

		\noindent (i) For all $x \in \sU \cap M_0$,  the limsup LILs at zero given in Theorems \ref{t:limsup0-1} and \ref{t:supprecise1} hold  with $\phi = \Phi_1$ and $\psi=\psi_0$.
		
		\noindent (ii) Assume  that $\sU=M$ and $R_1=R_4 =\infty$. Then for all $x \in M_0$,  the limsup LILs at infinity given in Theorems \ref{t:limsupinf-1} and \ref{t:supprecise2} hold  with $\phi = \Phi_1$ and $\psi=\psi_0$.
	\end{thm}

   We give two explicit examples, which are similar to  Examples \ref{e:mixedtype} and \ref{e:polydecay}(i).  In the following two examples, we assume that $F(r)=r^\gamma$ for some $\gamma>1$ (where  $F$ is the function in \eqref{e:diffusion}),  $X$ is a  symmetric Hunt process associated with the Dirichlet form \eqref{e:Form}, and Assumption L holds.
	
	\begin{example}
		{\rm 	Suppose that $\Lambda>0$ so that the Dirichlet form $(\EE^X,\FF^X)$ has non-zero strongly local part. Then  $\Phi_1(r) \asymp \Lambda^{-1}r^\gamma$ for $r \in (0,1]$ so that  by the change of the variables and (L1) (cf. \eqref{e:cal2.5}),
			$$\int_0  \frac{\Phi_1(r)\log|\log r|}{r \psi_0(r \log|\log r|)}dr \le c\int_0 \frac{s^{\gamma-1}(\log |\log s|)^{1-\gamma}}{\psi_0(s)}ds \le c \gamma^{-1}\int_0 \frac{d F(s)}{\psi_0(s)}<\infty.$$
			Thus, by Theorem \ref{t:Hunt0}, the limsup LIL at zero \eqref{e:ex2.5} holds for all $x \in \sU \cap M_0$ with $b=\Lambda$.  \qed
		}
	\end{example}
	
	\begin{example}
		{\rm  Suppose that $\ll =0$ and $\psi_0(r)=r^\gamma |\log r|^p$ (resp. $\psi_0(r)=r^\gamma |\log r|(\log |\log r|)^p$) for some $p \in \R$. Then by Theorem \ref{t:Hunt0}(i),  the limsup LIL at zero given in Example \ref{e:polydecay}(i-a) (resp. (i-b)) holds for all $x \in \sU \cap M_0$. \qed
		}
	\end{example}

\vspace{-2mm}

	Analogously, by  Theorem \ref{t:Hunt0}(ii), one can construct concrete examples of the LILs at infinity for symmetric pure-jump process from Example \ref{e:polydecay}(ii).

	\section{LILs for random conductance models}\label{s:RCM}

	Let $G=(\bL, E_\bL)$ be a locally finite connected infinite undirected graph, where $\bL$ is the set of vertices, and $E_\bL$  the set of edges. For $x,y \in \bL$, we write $d(x,y)$ for the graph distance, namely, the length of the shortest path joining $x$ and $y$. Let $\mu_c$ be the counting measure on $\bL$. Recall that $B(x,r):=\{y \in \bL : d(x,y)<r\}$ denotes the open ball of radius $r$ centred at $x$. In  this section, we always assume that there exists a constant $d>0$ such that
	\begin{equation}\label{e:graph-d-set}
		\mu_c (B(x,r)) \asymp r^d \quad \text{for}\;\; x \in \bL, \; r>10.
	\end{equation}
Obviously, if $G$ is the lattice graph $\Z^{d_1} \times \Z_+^{d_2}$,  then  \eqref{e:graph-d-set} holds with $d=d_1 + d_2$. Also, infinite graphs based on fractals such as graphical generalized Sierpinski gaskets and carpets, and the Vicsek tree satisfy \eqref{e:graph-d-set}.  For the precise definitions of these graphs, we refer to \cite{BB99graph, GT02}.

	 Let $(\eta_{xy}: x,y \in \bL)$  be a family of non-negative random variables  defined on a probability space $(\varOmega,\mathbf{F},\mathbf{P})$  such that   $\eta_{xx}=0$ and $\eta_{xy}=\eta_{yx}$ for all $x,y \in \bL$. The family $(\eta_{xy} : x,y \in \bL)$	is called a \textit{random conductance} on $\bL$.  
	We write $\nu_x= \sum_{y \in \bL} \eta_{xy}$ for $x \in \bL$. 
	 
 For each $\omega \in \varOmega$, let $X^\omega=(X^\omega_t, t \ge 0 ; \P^x_\omega, x \in \bL)$ be the  \textit{variable speed random walk} (VSRW)  associated with the random conductance $(\eta_{xy}:x,y \in \bL)$, that is, $X^\omega$ is the  process on $\bL$ that waits at each vertex $x$ for an exponential time with mean $1/\nu_x(\omega)$ and jumps according to the transitions $\eta_{xy}(\omega)/\nu_x(\omega)$. Then $X^\omega$ is the symmetric Markov process with $L^2(\bL;\mu_c)$-generator
  \begin{equation}\label{e:LV}
 	\sL^\omega_Vf(x)=\sum_{y}\eta_{xy}(\omega)(f(y)-f(x)), \quad x \in \bL.
 \end{equation}
	There is another  natural continuous time random walk $Y^\omega=(Y^\omega_t, t \ge 0 ; \P^x_\omega, x \in \bL)$, which is called the \textit{constant speed random walk} (CSRW), associated with  $(\eta_{xy}:x,y \in \bL)$. This process $Y^\omega$ also jumps according to the transitions $\eta_{xy}(\omega)/\nu_x(\omega)$, however, it waits at each vertex $x$ for an exponential time with mean $1$. Note that the CSRW $Y^\omega$ is the symmetric Markov process  with $L^2(\bL;\nu)$-generator
	\begin{align*}
		&\sL^\omega_C f(x) = \nu_x^{-1}(\omega) \sum_{y}\eta_{xy}(\omega)(f(y)-f(x)), \quad x \in \bL.
	\end{align*}

One of the most	important examples of random conductance model is the  \textit{bond percolation} in  $\Z^d$. We recall the definition of this model:
	Let $G$ be the lattice graph $\Z^d$, $d \ge 1$ and $p \in [0,1]$. For edges $e \in E_{\Z^d}$, let $\eta_e$ be i.i.d. Bernoulli random variables in $(\varOmega,\mathbf{F},\mathbf{P})$ with $\mathbf{P}(\eta_e = 1) = p $. Then for $x,y \in \Z^d$, we assign $\eta_{xy}=\eta_e$ if $\{x,y\}=e \in E_{\Z^d}$ and $\eta_{xy} = 0$ otherwise. In each configuration $\omega \in \varOmega$, edges $e \in  E_{\Z^d}$  with $\eta_e(\omega)=1$ are called  \textit{open} and a subset $\sC \subset \Z^d$ is called an \textit{open cluster} if every $x, y \in \sC$ are connected by an open path. 
	It is known that there exists a critical probability $0<p_c<1$ depends on  $d$ such that if $p>p_c$, then for $\mathbf P$-a.s.  $\omega$, there exists  a unique infinite open cluster, which we denote $\sC_\infty = \sC_\infty(\omega)$.

	In the bond percolation with parameter $p>p_c$,  the following limsup LIL for CSRW is obtained in \cite{Dum08} using the weak Gaussian bounds derived in \cite{Ba04}: For $\mathbf{P}$-a.s. $\omega$ and all $x \in \sC_\infty(\omega)$, the CSRW $Y^\omega$ satisfies that   
	$ \limsup_{t \to \infty} {|Y_t^\omega|}/{(t\log \log t)^{1/2}} = c(p,d) \in (0, \infty)$, $\P^x_\omega$-a.s. This result was extended in \cite{CK16, Ch17} to cover more general models such that enjoys weak sub-Gaussian heat kernel estimates (see \cite[Assumption 1.1]{Ch17}). 	In this section, we study limsup LIL for  random conductance models  with long jumps that do not enjoy weak sub-Gaussian heat kernel estimates.

	\medskip

	Let $\alpha \in (0,2)$ and $(w_{xy}:x,y \in \bL)$ be a family of non-negative independent random variables such that $w_{xy}=w_{yx}$ for all $x \neq y \in \bL$. Define
\begin{equation}\label{e:eta}
\eta_{xx}=0 \quad \text{and} \quad \eta_{xy}(\omega) = \frac{w_{xy}(\omega)}{|x-y|^{d+\alpha}}, \quad  \, x\neq y, \;\; x,y \in \bL.
	\end{equation}
As an application of \cite{CKW18, CKW20-1} and our Theorem \ref{t:limsupinf-1}, we have the following limsup LIL for random conductance models  with long jumps.

Let $\mathbf E$ be the expectation with respect to $\mathbf{P}$.

\begin{thm}\label{t:LIL-RCM}
Suppose  $X^\omega$ (resp.   $Y^\omega$) is the VSRW (resp.  CSRW)
associated with  $(\eta_{xy}:x,y \in \bL)$ defined as  \eqref{e:eta}.  Suppose that $d>2(2-\alpha)$,
\begin{equation}\label{e:pq0} \sup_{x,y \in \bL, \, x \neq y} \mathbf{P}(w_{xy} = 0) < 1/2 \qquad \mbox{and} \quad \sup_{x,y \in \bL, \, x \neq y} \big( \mathbf E[w_{xy}^{\,p}]+ \mathbf E[w_{xy}^{-q} \, \1_{\{w_{xy} > 0\}}] \big) < \infty 	\end{equation} 
for some constants
\begin{equation}\label{e:pq}
	p > \frac{d+2}{d} \lor \frac{d+1}{2(2-\alpha)} \quad \mbox{and} \quad q> \frac{d+2}{d}.
\end{equation}
When we consider the CSRW $Y^\omega$,  we also assume  that there are constants $m_2 \ge m_1>0$ such that $\mathbf P$-a.s. $\omega$,
\begin{equation}\label{e:CSRW-condi}
	w_{xy}(\omega)>0  \;\; \text{for all} \;\, x,y \in \bL, \, x \neq y \quad\; \text{and} \;\quad  m_1 \le   \sum_{y \in \bL, \, y \neq x} \frac{w_{xy}(\omega)}{|x-y|^{d+\alpha}} \le m_2 \;\; \text{for all} \;\, x \in \bL.
\end{equation}
Then $X^\omega$ and $Y^\omega$  enjoy the limsup law at infinity \eqref{4} with $\phi(r)=r^\alpha$. Moreover, when $\alpha \in (0,1)$, the coclusion is still valid if $d>2(1-\alpha)$ and \eqref{e:pq0} holds with
\begin{equation}\label{e:pq'}
	p > \frac{d+2}{d} \lor \frac{d+1}{2(1-\alpha)} \quad \mbox{and} \quad q> \frac{d+2}{d}.
\end{equation}
\end{thm}
\pf  $\mathrm{VRD}^{r_1(\omega)}$  
with the reference measure $\mu=\mu_c$ holds by \eqref{e:graph-d-set}.
	By exactly same the Borel-Cantelli arguments in the proof of  \cite[Proposition 5.6]{CKW18}, we can deduce from \eqref{e:pq0}-\eqref{e:pq} that there is a constant  $\up_0\in (0,1)$ such that for  $\mathbf P$-a.s. $\omega$,	there is $r_1(\omega) \ge 1$ such that  the assumption (Exi.$(\theta)$) 	 in \cite{CKW18}  holds with the reference measure $\mu=\mu_c$, and  the associated constants $R_0\ge 1$ and $\theta \in (0,1)$ being replaced by $r_1(\omega)$ and $\up_0$ respectively.    Then by  \cite[Lemma 3.4]{CKW18}, for  $\mathbf P$-a.s. $\omega$, $\mathrm{TJ}^{r_1(\omega)}(r^\alpha)$ for the VSRW $X^\omega$ holds.

	In the following, we show that  there exist constants $\up, \theta' \in (\up_0,1)$ and $c_i>0$, $i=1,2,3$  such that for $\mathbf P$-a.s. $\omega$, there is $r_2(\omega) \ge 1$ so  that for all $x_0 \in \bL$, $R \ge r_2(\omega) (1+|x_0|)^\up$  and $R^{\theta' \alpha} \le t \le c_1 R^\alpha$,
	\begin{equation}\label{e:RCM-NDL1}
		p^{B(x_0,R)}(t,y,z) \ge c_2 t^{-d/\alpha}, \quad  y,z \in B(x_0,R/4) \;\; \text{with} \;\; |y-z| \le 4c_3t^{1/\alpha},
	\end{equation}
Using \cite[Theorem 3.3]{CKW18}, we see that  there exist constants $c_4>0$ and  $\up_0<\delta<\theta'<1$ such that  for $\mathbf P$-a.s. $\omega$, there is $r_3(\omega) \ge 1$  such that for all $s \ge r_3(\omega)$, $s^\delta \le r \le s$ and $t \ge r^{\theta' \alpha}$,
	\begin{equation}\label{e:RCM1}
		\sup_{x \in B(0,2s)} \P^x_\omega(\tau_{B(x,r)} \le t ) \le c_4 \bigg( \frac{t}{r^\alpha}\bigg)^{1/3}.
	\end{equation}
Let $\up := \delta/\theta' \in (\delta,1)$. Fix $x_0 \in \bL$.
 Then  for  $\mathbf P$-a.s. $\omega$ and all $R \ge 2r_3(\omega) (1+|x_0|)^\up$, $x \in B(x_0,R)$, $R^{\theta'} \le r \le R/2$
  and $t \ge 2R^{\theta' \alpha}$, by letting $s=(R+ |x_0|)/2$ in \eqref{e:RCM1}, we get that
  	\begin{equation}\label{e:RCM.2.30}
  	\sum_{y \in B(x,r)^c} p^{B(x_0,R)}(t/2,x,y) \le \P_\omega^x(\tau_{B(x,r)} \le t/2)  \le    \sup_{y \in B(0,|x_0|+R)}  \P_\omega^y (\tau_{B(y,r)} \le t) \le  c_4 \bigg( \frac{t}{r^\alpha}\bigg)^{1/3}. 
  	\end{equation}
Indeed, in the above situation, the constant $s=(R+ |x_0|)/2$ satisfies that $s \ge R/2 \ge r_3(\omega)$, $s \ge r$ and $s^\delta \le (R \vee |x_0|)^\delta = R^{\theta'\up} \vee |x_0|^{\theta' \up} \le R^{\theta'}$.
Since the result of  \cite[Theorem 2.10]{CKW20-1} holds with the constant $\delta \in (\up_0,1)$ above by \cite[Theorem 3.6]{CKW18},  by following \textbf{Step (1)} in the proof of \cite[Proposition 2.11]{CKW20-1}, using \eqref{e:RCM.2.30} instead of  \cite[(2.30)]{CKW20-1}, we can deduce that  \eqref{e:RCM-NDL1} holds.

In the end, since  $\mathrm{NDL}^{r_2(\omega)}$ for $X^\omega$ holds by \eqref{e:RCM-NDL1}, we deduce from Proposition \ref{p:E}(ii) and Theorem \ref{t:limsupinf-1}(i) that the desired limsup law for the VSRW $X^\omega$ holds. Since $\nu_x$ is comparable with the counting measure $\mu_c$ under the assumption \eqref{e:CSRW-condi}, one can also conclude that the desired limsup for the CSRW $Y^\omega$ holds by the same way. 
 
Now we suppose that $\alpha \in (0,1)$ and \eqref{e:pq0} holds with some $p,q$ satisfying  \eqref{e:pq'}. Then using the Borel-Cantelli arguments again, we get from \eqref{e:pq0} and \eqref{e:pq'} that there is $\up_0 \in (0,1)$ such that for  $\mathbf P$-a.s. $\omega$,	there is $r_1(\omega) \ge 1$ such that  the assumption (Exi.$(\theta)'$) 	 in \cite{CKW18} and \cite[(3.10)]{CKW18} hold with $\mu=\mu_c$,  $R_0=r_1(\omega)$ and $\theta=\up_0$.  Hence,   $\mathrm{TJ}^{r_1(\omega)}(r^\alpha)$ still holds true by  \cite[Lemma 3.4]{CKW18}. Moreover, according to \cite[Remark 3.10(2)]{CKW18}, the results of \cite[Theorems 3.3 and 3.6]{CKW18} are still valid under (Exi.$(\theta)'$). Therefore, by following the argument above, we get the desired result from our Theorem  \ref{t:limsupinf-1}(i). We have finished the proof.   \qed

\begin{remark}
	{\rm (i) The assumption that $w_{xy}>0$ for all $x,y \in \bL, \, x \neq y$ in \eqref{e:CSRW-condi} can be weaken to that $\sup_{x,y \in \bL, \, x \neq y} \mathbf{P}(w_{xy} = 0)<c_1$ for some  $c_1 \in (0,1)$ depends on the constants $m_1,m_2$. See \cite[Remrak 4.4(1)]{CKW20-1}.

		\noindent (ii)  We emphasize that the assumption \eqref{e:pq} is weaker than the assumptions of \cite[Theorems 4.1 and 4.3]{CKW20-1}. Hence,  under \eqref{e:pq}, even for 
 $\bL=\Z^{d_1} \times \Z_+^{d_2}$, we do not know whether $X^\omega$ and $Y^\omega$ enjoy two-sided $\alpha$-stable like heat kernel estimates or not  (see \cite[(4.4)]{CKW20-1}). However, since our LIL works without any information on off-diagonal heat kernel estimates, we succeeded in obtaining LILs for $X^\omega$ and $Y^\omega$. 
}	
\end{remark}

	\section{Estimates on the first exit times  and zero one Law
	}\label{s:pre}

	\subsection{Auxiliary functions $\Up_1, \Up_2, \vt_1$ and  $\vt_2$}\label{s:aux}
	
\begin{figure}[!h]
		\centering \hspace{-0.9cm}
		\includegraphics[width=0.54\columnwidth]{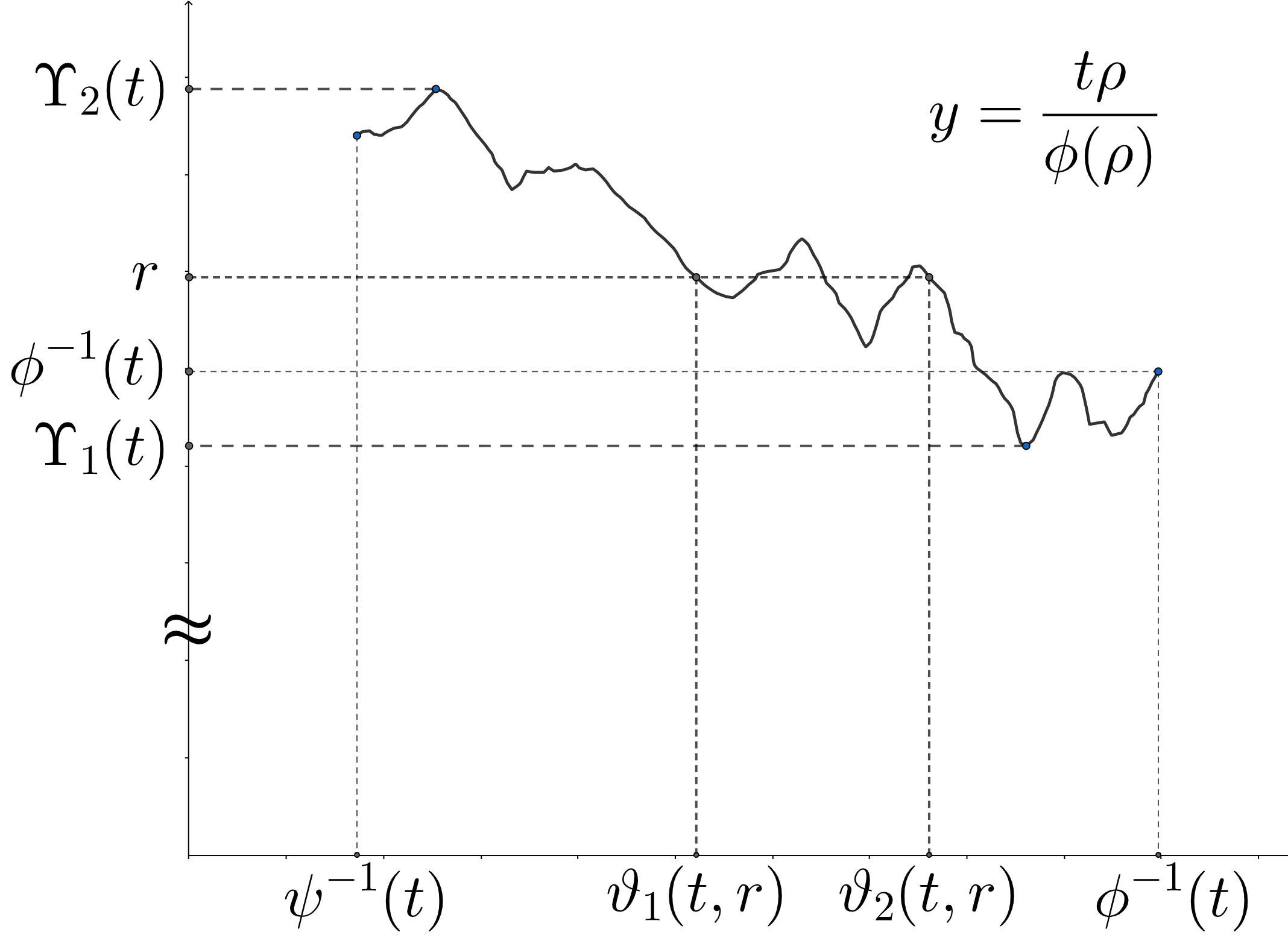}
		\vspace{-0.2cm}
		\caption{Auxiliary functions}\label{f:1}
	\end{figure}

	In this subsection, we introduce auxiliary functions which will be used in tail probability estimates on the first exit times from balls (see Propositions \ref{p:EP} and \ref{p:EPL}).
	Recall that we always assume \eqref{e:phipsi}. Thus, $\phi^{-1}(t) \ge \psi^{-1}(t)$ for all $t > 0$. We define 
	\begin{equation}
		\Up_1(t):=\min\left\{\frac{t \rho}{\phi(\rho)}: \psi^{-1}(t)  \le \rho \le \phi^{-1}(t) \right\} \quad \text{and} \quad \Up_2(t):=\max\left\{\frac{t \rho}{\phi(\rho)}: \psi^{-1}(t)  \le \rho \le \phi^{-1}(t) \right\}.
	\end{equation}
By the definition, it holds that
	\begin{equation}\label{e:Up}
		\Up_1(t) \le \frac{t \phi^{-1}(t)}{\phi(\phi^{-1}(t))} = \phi^{-1}(t) \le \Up_2(t) \quad\;\; \text{for all} \;\; t>0.
	\end{equation}
	We also define functions $\vt_1, \vt_2:(0, \infty) \times [0, \infty) \to (0,\infty)$ as (cf. \cite[(1.13)]{CK202})
	\begin{equation}\label{e:theta1}
		\vt_1(t,r):=\begin{cases}
			\phi^{-1}(t), &\mbox{if} \;\; r \in [0, \Up_1(t)),\\
			\displaystyle\min\Big\{\rho \in [\psi^{-1}(t),  \; \phi^{-1}(t)]: \frac{t \rho}{\phi(\rho)} \le r \Big\}, \quad  &\mbox{if} \;\; r \in [\Up_1(t), \Up_2(t)],\\
			\psi^{-1}(t), &\mbox{if} \;\; r \in (\Up_2(t), \infty)
		\end{cases}
	\end{equation}
	and
	\begin{equation}\label{e:theta2}
		\vt_2(t,r):=\begin{cases}
			\phi^{-1}(t), &\mbox{if} \;\; r \in [0, \Up_1(t)),\\
			\displaystyle\max\Big\{\rho \in [\psi^{-1}(t),  \; \phi^{-1}(t)]: \frac{t \rho}{\phi(\rho)} \ge r \Big\}, \quad &\mbox{if} \;\; r \in [\Up_1(t), \Up_2(t)],\\
			\psi^{-1}(t), &\mbox{if} \;\; r \in (\Up_2(t), \infty).
		\end{cases}
	\end{equation}
	See Figure \ref{f:1}. Since $\phi$ and $\psi$ are continuous, the above $\vt_1$ and $\vt_2$ are well-defined. For each fixed $t>0$, $r \mapsto \vt_1(t,r)$ and $r \mapsto \vt_2(t,r)$ are non-increasing.
	Intuitively, $\Up_2(t)$ represents the maximal distance reachable by a standard chaining argument at time $t$, and $\vt_1(t,r)$ and $\vt_2(t,r)$ represent the number of minimal and maximal number of chains to reach the distance $r$, respectively.
	
	\vspace{-0.4cm}
	
	\begin{lem}\label{l:theta1}
	For all $\eps \in (0,1)$ and $t,r>0$ satisfying  $\Up_1(t) \le \eps r$, it holds that
		\begin{equation*}
			\frac{\eps r}{\vt_1(t, \eps r)} \ge \frac{t}{\phi(\vt_1(t, \eps r))}.
		\end{equation*}
	\end{lem}
	\pf If $\eps r \le \Up_2(t)$, then $t \vt_1(t,\eps r)/\phi(\vt_1(t,\eps r)) \le \eps r $ and hence 
	\begin{equation*}
		\frac{r}{\vt_1(t, \eps r)} - \frac{t}{\phi(\vt_1(t, \eps r))} = \frac{1}{\vt_1(t, \eps r)}\Big(r - \frac{t \vt_1(t, \eps r)}{\phi(\vt_1(t, \eps r))}\Big) \ge  \frac{(1-\eps)r}{\vt_1(t, \eps r)}.
	\end{equation*}
	Otherwise, if $\eps r>\Up_2(t)$, then since $t \vt_1(t,\eps r)/\phi(\vt_1(t,\eps r)) \le \Up_2(t) $, we also obtain
	\begin{equation*}
		\frac{r}{\vt_1(t, \eps r)} - \frac{t}{\phi(\vt_1(t, \eps r))} = \frac{1}{\vt_1(t, \eps r)}\Big(r - \frac{t \vt_1(t, \eps r)}{\phi(\vt_1(t, \eps r))}\Big) \ge \frac{r - \Up_2(t)}{\vt_1(t, \eps r)}>\frac{(1-\eps)r}{\vt_1(t, \eps r)}.
	\end{equation*}
	\qed
	
	The next lemma follows from the inequality $e^{-x} \le \beta^{\beta} x^{-\beta}$, $x, \beta>0$.
	\begin{lem}\label{l:exp1}
		(i) If $\U_{R_0}(\psi, \beta_3, C_U')$ holds, then for all $k>0$ and $0< \psi^{-1}(t) \le r<R_0$,
		\begin{equation}\label{e:exp1}
			\exp\left(-{k r}/{\psi^{-1}(t)}\right) \le {C_U'\beta_3^{\beta_3} k^{-\beta_3}t}/{\psi(r)}.
		\end{equation}
		
		\noindent 	(ii) If $\U^{R_\infty}(\psi, \beta_3, C_U')$ holds, then  \eqref{e:exp1} holds for all $k>0$ and  $R_\infty < \psi^{-1}(t) \le r$.
	\end{lem}

	\subsection{Estimates on the first exit times from open balls}
	
	In this subsection, we study the first exit time from open balls. The main results  are Propositions \ref{p:EP} and \ref{p:EPL}. Note that a  similar result appears in \cite[(3.25)]{BKKL19b} (see also \cite[(3.3)]{CKW20}). However, in this paper, we not only give  local versions (in the variable $r$) of that result, but also remove an extra assumption that the local lower scaling index of $\phi$ is strictly larger than $1$ therein, by using auxiliary functions defined in subsection \ref{s:aux} and adopting the ideas from the first and second named authors' paper \cite[Theorem 3.10]{CK202}.

	In the following proposition,  we let $\ell>1$ be the constant in \eqref{RVD},
	\begin{equation}\label{e:defr0}
		r_0:=\eta^2R_0/(6\ell) \;\; \text{with the constant $\eta \in (0,1)$  from ${\mathrm {NDL}}_{R_0}(\phi, U)$}
	\end{equation} 
	in the first statement, and
	\begin{equation}\label{e:defrinf}
		r_\infty:= (4\ell R_\infty)^{1/(1-\up)} \vee (2R_\infty/\eta^2) \;\; \text{with the constants $\up, \eta \in (0,1)$ from \NDLi}
	\end{equation}
	in the latter statement.  Note that $r_0=\infty$ if $R_0=\infty$.

	\begin{prop}\label{p:E}
		(i) Suppose that \VRDo, ${\mathrm {NDL}}_{R_0}(\phi, U)$ and  $\U_{R_0}(\phi,\beta_2,C_U)$  hold. 
		Then there exist constants $C_1',A_0',a_1,a_2>0$ such that for all  $n \ge 1$, $x \in U$ and  $0<r<r_0 \land (C_1' \updelta_U(x))$, 
		\begin{equation}\label{e:SP}
		e^{-a_1n} \le	\P^x(\tau_{B(x,r)} \ge n \phi(A_0' r) ) \le  e^{-a_2n}.
		\end{equation}
		Moreover,  $\mathrm{E}_{r_0}(\phi, U)$ holds.

		\noindent	(ii) Suppose that  \VRDi, \NDLi \ and  $\U^{R_\infty}(\phi,\beta_2,C_U)$ hold. 		 Then there exist constants $A_0',a_1,a_2>0$ such that \eqref{e:SP} holds for all $n\ge1$, $x \in M$ and $r>r_\infty \ep(x)^\up$, and 	$\mathrm{E}^{r_\infty}(\phi)$ holds.
	\end{prop}
	\pf We adopt the idea from the proof for \cite[Proposition 3.5(ii)]{CKW16b}.
	
 (i)   Let  $C_1'=\eta^2(C_V \land C_2)/(5\ell)$ and $A_0'=4\ell/\eta$, where $C_V,C_2,\eta\in(0,1)$ are the constants in \VRDo, \NDLo \ and \NDLo, respetively, and $\ell>1$ is the constant in \eqref{RVD}. 
 Choose $x \in U$ and $0<r<r_0 \land (C_1' \updelta_U(x))$. Set $t_r:= \phi(3 \ell \eta^{-1} r)$. Then $\phi^{-1}(t_r)/\eta = 3\ell r/\eta^{2} <  R_0 \land ((C_V \land C_2)\updelta_U(x))$. Thus, by ${\mathrm {NDL}}_{R_0}(\phi, U)$ and \VRDo, it hold that
	\begin{equation}\label{e:E_0-1}
		p^{B(x,\phi^{-1}(t_r)/\eta)}(t_r,y,z) \ge \frac{c_1}{V(x,\phi^{-1}(t_r)/\eta)} \ge \frac{c_2}{V(x,r)}, \quad\;\; \forall \, y,z \in B(x,3\ell r).
	\end{equation}
 Meanwhile, by \eqref{RVD},  one can pick $w \in B(x, 2\ell r) \setminus B(x,2r)$. 
	Observe that $B(w, r) \cap B(x,r) = \emptyset$,    $B(x, r) \subset B(w,3\ell r)$, $B(w, r) \subset B(x,3\ell r)$ and $\updelta_U(w) \ge \updelta_U(x)-2\ell r > (C_1'^{-1}-2\ell)r \ge 3C_V^{-1}\ell r$. By  \VRDo, it follows that $V(w,r) \ge c_3V(w, 3 \ell r) \ge c_3 V(x, r)$. Then  by  \eqref{e:E_0-1}, it holds that 
	\begin{equation}\label{e:E_0-2}
		\P^y\big(X_{t_r} \notin B(x,r)\big) \ge \P^y\big(X^{B(x, \phi^{-1}(t_r)/\eta)}_{t_r} \in B(w, r) \big) \ge \frac{c_2V(w,r)}{V(x,r)} 
		\ge c_2c_3,  \quad
		\forall \, y \in B(x,r).
	\end{equation}
Using the Markov property, we get from \eqref{e:E_0-2} that for all $n \ge 1$,
	\begin{align}\label{e:E_0-4}
		\P^x\big(\tau_{B(x,r)} \ge n \phi(A_0' r)\big) &\le \E^x \left[ \1_{\{ \tau_{B(x,r)} \ge  (n-1) \phi(A_0' r) \}} \P^{X_{(n-1) \phi(A_0' r)}}( \tau_{B(x,r)} \ge \phi(A_0'r)  )\right] \nn\\[3pt]
		& \le \P^x\big(\tau_{B(x,r)} \ge  (n-1) \phi(A_0' r)\big) \sup_{y \in B(x,r)} \P^y \big(\tau_{B(x,r)} \ge  \phi(A_0' r)\big) \nn\\
		&\le  (1-c_2c_3) \, \P^x\big(\tau_{B(x,r)} \ge (n-1) \phi(A_0' r)\big) \le  \dots \le  (1-c_2c_3)^{n}.
	\end{align}

	Next, we prove the  lower bound in \eqref{e:SP}. Set $k:=\lceil C_U(A_0'/\eta)^{\beta_2} \rceil$.
	Using  $\U_{R_0}(\phi, \beta_2, C_U)$, the semigroup property and \VRDo, we get that for all $n \ge 1$,
	\begin{align}\label{e:E_0-6}
		&\P^x\big(\tau_{B(x,r)} \ge  n\phi(A_0' r) \big) \ge \P^x\big(\tau_{B(x,r)} \ge kn\phi(  \eta r) \big)  \ge \int_{B(x,\eta^2 r)} p^{B(x,r)}(kn\phi(\eta r),x,z_{kn})\mu(dz_{kn}) \nn \\
		&\ge \int_{B(x,\eta^2 r)} \dots \int_{B(x,\eta^2 r)}  p^{B(x,r)}(\phi(\eta r),x,z_1)  \dots  p^{B(x,r)}(\phi(\eta r),z_{kn-1},z_{kn})   \mu(dz_1) \dots \mu(dz_{kn})  \nn \\
		&\ge \bigg(\frac{c_4}{V(x,r)}\bigg)^{kn} \int_{B(x,\eta^2 r)} \dots \int_{B(x,\eta^2 r)}    \mu(dz_1) \dots \mu(dz_{kn}) =	 \bigg(\frac{c_4V(x, \eta^2 r)}{V(x,r)}\bigg)^{kn}  \ge e^{-c_5n}.
	\end{align}
	We have established  \eqref{e:SP}. Now using  Markov inequality, we get from \eqref{e:SP} that 
	\begin{equation}\label{e:SP-E-1}
		\E^x[\tau_{B(x,r)}] \ge  \phi(A_0' r)\P^x(\tau_{B(x,r)} \ge  \phi(A_0' r) ) \ge e^{-a_1} \phi(A_0' r)
	\end{equation}
	and
	\begin{align}\label{e:SP-E-2}
		\E^x [\tau_{B(x,r)}] &\le \sum_{n=1}^\infty  n \phi(A_0' r) \P^x\big(\tau_{B(x,r)} \ge (n-1) \phi(A_0' r)  \big)\nn\\
		&\le  \phi(A_0' r)\sum_{n=1}^\infty  ne^{-a_2(n-1)}  = c_6\phi(A_0' r).
	\end{align} 
Therefore, $\mathrm{E}_{r_0}(\phi, U)$ holds.

 (ii) We follow the proof of (i). Choose any  $x \in M$, $r>r_\infty \ep(x)^\up$ and let $t_r:= \phi(3\ell \eta^{-1} r)$. Since  $
	\phi^{-1}(t_r)/\eta = 3\ell\eta^{-2} r  > R_\infty\ep(x)^\up$,  by \NDLi \ and \VRDi,  \eqref{e:E_0-1} is still valid. Also, we can pick  $w \in B(x, 2\ell r) \setminus B(x,2r)$ by \eqref{RVD}. Using the inequality $(a+b)^\up \le a^\up +b^\up$ for $a,b \ge 0$, we get
	\begin{equation*}
		R_\infty \ep(w)^\up \le R_\infty (\ep(x)+2\ell r)^\up \le R_\infty (\ep(x)^\up + (2\ell r)^\up) < (2^{-1} + 2 \ell R_\infty r_\infty^{\up-1})r \le r.
	\end{equation*}
	Hence, we see from \VRDi \ that $V(w,r) \ge c_1V(w,3\ell r) \ge c_1V(x,r)$ and \eqref{e:E_0-2} is true. Using the Markov property again, we deduce that \eqref{e:E_0-4}  holds.  On the other hand, we can still get \eqref{e:E_0-6} from $\U^{R_\infty}(\phi, \beta_2, C_U)$ and \NDLi \ since $\eta^2 r_\infty>R_\infty$. Therefore, \eqref{e:SP} is true. Now, we conclude that 	$\mathrm{E}^{r_\infty}(\phi)$ holds by \eqref{e:SP-E-1} and \eqref{e:SP-E-2}.	 \qed

	For sake of brevity, we give three families of conditions.

	\begin{assumption}\label{ass0}
		{\rm  \Eo,  ${\mathrm {TJ}}_{R_0}(\phi,  \le,U)$  and  $\U_{R_0}(\phi,\beta_2,C_U)$ hold  (with $C_1$) for  some $R_0 \in (0,\infty]$ and an open set $U \subset M$.
	}
	\end{assumption}

	\begin{assumption}\label{assinf}	{\rm     \Ei,  ${\mathrm {TJ}}^{R_\infty}(\phi, \le)$  and $\U^{R_\infty}(\phi,\beta_2,C_U)$   hold  for some  $R_\infty \ge1$ and $\up \in (0,1)$.}
	\end{assumption}

	\begin{assumption}\label{assinf2}
		{\rm \VRDi, \NDLi,  ${\mathrm {TJ}}^{R_\infty}(\phi, \le)$,  $\L^{R_\infty}(\phi,\beta_1,C_L)$ and  $\U^{R_\infty}(\phi,\beta_2,C_U)$   hold for some $R_\infty\ge1$ and $\up \in (0,1)$.}
	\end{assumption}

By  Proposition \ref{p:E}(ii), \asss \ is stronger than \ass. \asss \ will be used in the next subsection to establish estimates on oscillation of bounded parabolic functions and obtain a zero-one law.

	For $\rho>0$, let $X^{(\rho)}$ be a Markov process on $M$ obtained by eliminating all jumps with absolute jump size bigger than $\rho$ from $X$. Then $X^{(\rho)}$ is a Borel standard Markov process on $M$ with the L\'evy measure $J^{(\rho)}(x,dy):=\1_{B(x,\rho)}(y)J(x,dy)$.
	By the Meyer's construction  (see \cite{Me75}, and also \cite[Section 3]{BGK09}),  $X$ can be constructed from $X^{(\rho)}$ by attaching large jumps whose sizes are bigger than $\rho$. 
	For an open set $D \subset M$, we denote $\tau^{(\rho)}_D := \inf\{ t>0 : X^{(\rho)}_t \in M_\partial \setminus D\}$.

	Our first goal is obtaining Proposition \ref{p:EP} below. For this, we prepare two lemmas. 
	The main strategies  of the proofs for the next two lemmas are similar to  those for \cite[Lemmas 4.20 and 4.21]{CKW16a}.  
	Under \ass, for each $\up_1 \in (\up,1)$, we define
	\begin{equation}\label{e:infty'}
		R_\infty'=R_\infty'(\up_1):=(7R_\infty)^{\up_1/(\up-\up\up_1)}.
	\end{equation}

	\begin{lem}\label{l:4.20}
		(i) Suppose that \as \ holds.  Then, there exist constants $\delta \in (0,1)$ and $C_3 \ge 1$ such that for all $x \in U$,  $0<\rho \le r<3^{-1}\big(R_0 \land (C_1 \updelta_U(x))\big)$ and  $\lambda \ge C_3/\phi(\rho)$,
		\begin{equation}\label{e:4.20}
			\E^x \left[ \exp\big(-\lambda \tau^{(\rho)}_{B(x,r)}\big)\right] \le 1-\delta.
		\end{equation}
		(ii) Suppose that \ass \ holds. Then,  there exist constants  $\delta \in (0,1)$ and $C_3\ge 1$ such that \eqref{e:4.20} holds for any $\up_1 \in (\up,1)$, $x_0\in M$, $r>R_\infty' \ep(x_0)^{\up_1}$, $x \in B(x_0,6r)$, $2r^{\up/\up_1} \le \rho \le r$  and $\lambda \ge C_3/\phi(\rho)$, where $R_\infty'=R_\infty'(\up_1)$ is defined by \eqref{e:infty'}.
	\end{lem}
	\pf (i) Choose any  $x \in U$ and $0<r< 3^{-1}\big(R_0 \land (C_1 \updelta_U(x))\big)$. Set $B := B(x,r)$. Then we have
	\begin{equation}\label{e:3.8.1}
		C_1\updelta_U(z) \ge C_1(\updelta_U(x)-r) > 2r \quad \text{for all} \;\; z \in B.
	\end{equation}
	Hence,  since  \Eo \ holds, by using the inequality $\tau_B \le t + (\tau_B-t)  1_{\{\tau_B > t\}}$ and the Markov property, we get that  for all $t>0$,
	\begin{align*}
		c_1 \phi(r) &\le \E^x \tau_B \le t + \E^x \left[ \1_{\{  \tau_B>t \} } \E^{X_t} \tau_B  \right] \le t+  \P^x(\tau_B >t ) \sup_{z \in B} \E^z \tau_B \\[2pt]
		& \le t + \P^x(\tau_B >t ) \sup_{z \in B} \E^z \tau_{B(z,2r)} \le t + c_2 \P^x (\tau_B >t ) \phi(2A_0r).
	\end{align*}
	Since $\U_{R_0}(\phi,\beta_2,C_U)$ holds and $\phi$ is increasing, it follows that 
	\begin{equation}\label{e:4.16}
		\P^x( \tau_B>t) \ge  \frac{c_1\phi(r)}{c_2\phi(2A_0r)} - \frac{t}{c_2\phi(2A_0r)} \ge c_3  - \frac{t}{c_3 \phi(r)} \quad \text{for all} \;\; t>0,
	\end{equation}
for some constant $c_3\in(0,1)$.	We claim that there exists a  constant $c_4>0$ independent of $x$ and $r$ such that
	\begin{equation}\label{e:7.8}
		\Big|\P^x\big(\tau_B > t\big) - \P^x\big( \tau_B^{(\rho)} > t\big)\Big| \le  \frac{c_4t}{\phi(\rho)} \quad\;\; \text{ for all} \;\; \rho \in (0, r], \;\; t>0.
	\end{equation}
	To prove \eqref{e:7.8}, we need some preparations. Choose any $\rho \in (0,r]$.  Let $\xi$ be an exponential random variable with rate parameter $1$ independent of  $X^{(\rho)}$. Define $\sJ_t := \int_0^t   J\big(X_s^{(\rho)}, M_\partial \setminus B(X_s^{(\rho)},\rho)\big)  ds$ and $T_0:=\inf \{ t \ge 0: \sJ_t \ge \xi\}$. 
	In view of the Meyer's construction, we may identify $T_0$ with $T_1:= \inf\{t>0: X_t \neq X^{(\rho)}_t\}$, the first attached jump time for $X^{(\rho)}$ (see \cite[Section 3.1]{BGK09}).  Since $\Tail_{R_0}(\phi, \le,U )$ holds and  $X_s^{(\rho)} \in B$ for all $0 \le s < t \land \tau_B^{(\rho)}$, by \eqref{e:3.8.1}, we have that 
	\begin{equation}\label{e:7.8.1}
		\P^x\big(T_1 \le t \land \tau_B^{(\rho)}\big) = \P^x\big( \sJ_{t \land \tau_B^{(\rho)}} \ge \xi\big) \le \P^x\big(\frac{c_5 t}{\phi(\rho)} \ge \xi\big) \le \frac{c_5 t}{\phi(\rho)}.
	\end{equation}
	Besides, since $X_s = X^{(\rho)}_s$ for $s < T_1$, we have  $\{t < T_1 \land \tau_B  \} = \{ t < T_1 \land \tau^{(\rho)}_B \}$ and $ \{ T_1 \le t < \tau_B  \} \subset \{ T_1 \le t \land \tau^{(\rho)}_B \}$  for all $t>0$. Thus, we can see from \eqref{e:7.8.1} that for all $t>0$,
	\begin{align*}
		&\Big| \P^x\big(t < \tau_B \big) - \P^x\big( t<\tau_B^{(\rho)} \big) \Big| = \Big| \P^x\big(T_1 \le t < \tau_B \big) - \P^x\big(T_1 \le  t<\tau_B^{(\rho)} \big) \Big|   \\  
		&\le \P^x\big( T_1 \le t < \tau_B \big)  + \P^x\big(T_1 \le t < \tau_B^{(\rho)} \big) \le 2 \, \P^x\big( T_1 \le t \land \tau_B^{(\rho)}\big) \le \frac{2c_5t}{\phi(\rho)}.
	\end{align*}
	This yields \eqref{e:7.8} with $c_4 = 2c_5$.

	Now, set $\delta = c_3/3 \in (0, 1/3)$ and $C_3=(c_3^{-1}+c_4)\delta^{-1}\log(\delta^{-1})$. In view of \eqref{e:4.16} and \eqref{e:7.8}, by taking $t_\rho = \delta \phi(\rho)/(c_3^{-1}+c_4)$, we conclude that for all $\rho \in (0,r]$ and $\lambda \ge C_3/\phi(\rho)$,
	\begin{align}\label{e:l3.8.}
		\E^x \left[ \exp\big(-\lambda \tau_B^{(\rho)}\big)\right] &\le \E^x \left[\exp\big(-\lambda \tau_B^{(\rho)}\big);\tau_B^{(\rho)} \le t_\rho \right] + 
		e^{-\lambda t_\rho}
		\le \P^x(\tau_B^{(\rho)} \le t_\rho) + e^{-\lambda t_\rho} \nn\\[2pt]
		&\le 1- c_3 +\frac{(c_3^{-1}+c_4) t_\rho}{\phi(\rho)}  + \exp \Big( -\frac{C_3 t_\rho}{\phi( \rho)}\Big) \le 1- c_3 + 2\delta \le 1- \delta.
	\end{align}
	Since $\delta$ and $C_3$ are independent of $x$ and $r$, this completes the proof.

 (ii) Choose any $x_0 \in M$, $r> R_\infty' \ep(x_0)^{\up_1}$, $x \in B(x_0,6r)$ and $2r^{\up/\up_1} \le \rho \le r$. By \eqref{e:infty'}, it holds that for any $z \in B(x,r) \subset B(x_0,7r)$,
	\begin{equation}\label{e:ballstable}
		R_\infty\ep(z)^\up \le R_\infty\ep(x_0)^\up +(7r)^\up R_\infty < (R_\infty R_\infty'^{-\up/\up_1}+  7^\up R_\infty R_\infty'^{\,\up-\up/\up_1})r^{\up/\up_1}  < 2r^{\up/\up_1} \le \rho.
	\end{equation}
	Hence, since \Ei, $\mathrm{TJ}^{R_\infty}(\phi, \le)$   and $\U^{R_\infty}(\phi, \beta_2,C_U)$ hold, we see that \eqref{e:4.16}  and  \eqref{e:7.8} hold by following the proof for (i). By repeating the calculations \eqref{e:l3.8.}, we obtain the result.
	\qed

	\begin{lem}\label{l:4.21}	(i) Suppose that \as \ holds. Then  there exist constants $C_4>0$, $C_5\in (0,1]$ such that for all $x \in U$ and  $0<\rho \le r<3^{-1}\big(R_0 \land (C_1\updelta_U(x))\big)$,
		\begin{equation}\label{e:4.21-1}
			\E^x \bigg[ \exp\Big(-\frac{C_3}{\phi(\rho)}\tau_{B(x,r)}^{(\rho)}  \Big)\bigg] \le C_4 \exp\big(-C_5 \frac{r}{\rho}\big),
		\end{equation}
	where $C_3$ is the constant in Lemma \ref{l:4.20}(i).
		
		\noindent (ii) Suppose that \ass \ holds. Then there exist constants $C_4>0$, $C_5\in (0,1]$ such that \eqref{e:4.21-1} holds for any $\up_1 \in (\up,1)$, $x_0 \in M$, $r>R_\infty'\ep(x_0)^{\up_1}$, $x \in B(x_0,5r)$ and  $2r^{\up/\up_1} \le \rho \le r$ where $C_3$ is the constant in Lemma \ref{l:4.20}(ii) and  $R_\infty'$ is defined by \eqref{e:infty'}.
	\end{lem}
	\pf  (i) Choose  any $x \in U$ and $0<\rho \le r<3^{-1}\big(R_0 \land (C_1\updelta_U(x))\big)$. If $6\rho>r$, then by taking $C_4$ larger than $e^6$, we are done. Hence, we assume that $6\rho \le r$.
	
	Set $\lambda_\rho:=C_3/\phi(\rho)$, $n_\rho:=\lceil r/(6\rho) \rceil$, $\tau_0:=\tau^{(\rho)}_{B(x,r)}$, $u(y) := \E^y[ e^{-\lambda_\rho \tau_0}]$ and $m_k := \Vert u \Vert_{L^\infty(\overline{B(x,k\rho)};\mu)}$. Let $\delta \in (0,1)$ be the constant in  Lemma \ref{l:4.20}(i). Fix any  $\delta' \in (0,\delta)$. Then for each $1 \le k \le 2n_\rho$, choose any $x_k \in \overline{B(x,k\rho)}$ such that $(1-\delta')m_k \le u(x_k) \le m_k$. Since
	$(2n_\rho + 2) \rho  < r/3 + 4\rho \le r$, we see that  $B(x_k,2\rho) \subset B(x,(k+2)\rho)  \subset B(x,r)$ for all $1 \le k \le  2n_\rho$.
	
	Set $\tau_k := \tau^{(\rho)}_{B(x_k,\rho)}$ and $v_k(y) := \E^y[e^{-\lambda_\rho \tau_k}]$. Note that $3^{-1}C_1\updelta_U(y) > 3^{-1}C_1(\updelta_U(x)-r) > 3^{-1}r \ge  \rho$ for all $y \in B(x,r)$. Hence, by Lemma \ref{l:4.20}(i), we obtain $v_k(x_k) \le 1- \delta$. 
	Then by the strong Markov property, it holds that  for all $1 \le k \le 2n_\rho-2$,
	\begin{align}\label{e:trunc}
		(1-\delta')m_k \le u(x_k) &= \E^{x_k}\big[e^{-\lambda_\rho \tau_0} \,  ; \, \tau_k \le  \tau_0 < \infty\big] = \E^{x_k}  \big[ e^{-\lambda_\rho \tau_k} e^{-\lambda_\rho (\tau_0 - \tau_k)} \,  ; \, \tau_k \le  \tau_0 < \infty\big] \nn\\
		&= \E^{x_k}\big[e^{-\lambda_\rho \tau_k} \E^{X^{(\rho)}_{\tau_k}}(e^{-\lambda_\rho \tau_0}) \big] = \E^{x_k}\big[e^{-\lambda_\rho \tau_k} u(X_{\tau_k}^{(\rho)})\big] \nn\\[3pt]
		&\le  v_k(x_k) \Vert u \Vert_{L^\infty(\overline{B(x_k,2\rho)};\mu)}  \le  v_k(x_k)  m_{k+2}\le (1-\delta)m_{k+2}.
	\end{align}
	In the first inequality above, we used the fact that $X^{(\rho)}_{\tau_k} \in \overline{B(x_k,2\rho)}$ since the jump size of $X^{(\rho)}$ can not be larger than $\rho$. The second inequality above holds since $B(x_k, 2\rho) \subset B(x, (k+2)\rho)$.
	
	In the end, by \eqref{e:trunc}, we conclude  that
	\begin{align*}
		\E^x \Big[ \exp\Big(-\frac{C_3}{\phi(\rho)}\tau_{B(x,r)}^{(\rho)}  \Big)\Big] =u(x) \le m_1 \le \Big( \frac{1-\delta}{1-\delta'} \Big)^{n_\rho -1} m_{2n_\rho -1} \le \frac{1-\delta'}{1-\delta}  \exp \Big( - \frac{r}{6\rho} \log \frac{1-\delta'}{1-\delta}  \Big).
	\end{align*}
	In the last inequality above, we used the fact that $m_{2n_\rho -1} \le \Vert u \Vert_{L^\infty(M;\mu)} \le 1$. This proves \eqref{e:4.21-1}.

 (ii) By following the proof for (i), using Lemma \ref{l:4.20}(ii), one can obtain \eqref{e:4.21-1}. Note that since $x \in B(x_0,5r)$, the point $x_k$ in \eqref{e:trunc} satisfies $x_k \in B(x,k\rho) \subset B(x_0,6r)$ so that we can apply Lemma \ref{l:4.20}(ii) in the counterpart of \eqref{e:trunc}.  \qed
	
	Recall that $\Up_1, \Up_2, \vt_1$ and  $\vt_2$ are auxiliary functions discussed in Subsection \ref{s:aux}.
	\begin{prop}\label{p:EP}
		(i) Suppose that \as \ holds.  Then, there exists a constant  $C_6\ge 1$ such that for all $x \in U$, $0< r<3^{-1}\big(R_0 \land (C_1\updelta_U(x))\big)$ and  $t>0$,
		\begin{equation}\label{e:EP_0}
			\P^x(\tau_{B(x,r)} \le t) \le \frac{C_6 t}{\phi(r)}. 
		\end{equation}
		If  \Tailol \ and $\U_{R_0}(\psi, \beta_3, C_U')$  also  hold, then there exist constants $C_7, C_8, C_9>0$ such that for all $x \in U$, $0< r<3^{-1}(R_0 \land (C_1\updelta_U(x)))$ and $t>0$, 
		\begin{equation}\label{e:EP_1}
			\P^x(\tau_{B(x,r)} \le t) \le C_7 \bigg( \frac{t}{\psi(r)} + \exp \Big(- \frac{C_8 r}{\vt_1(t,C_9r)} \Big) \bigg)
		\end{equation}
		
		\noindent (ii) Suppose that \ass \ holds. Then,  there exists a constant  $C_6\ge 1$ such that \eqref{e:EP_0} holds for any $\up_1 \in (\up,1)$, $x \in M$,  $r>R_\infty'\ep(x)^{\up_1}$  and $t\ge \phi(2r^{\up/\up_1})$ where  $R_\infty'$ is defined by \eqref{e:infty'}. 
		If  \Tailil \ and $\U^{R_\infty}(\psi, \beta_3, C_U')$ also hold, then there exist constants $C_7, C_8,C_9>0$ such that \eqref{e:EP_1} holds for all $x \in M$, $r>R_\infty'\ep(x)^{\up_1}$ and $t\ge \psi(2r^{\up/\up_1})$.
	\end{prop}
	\pf (i) First, assuming that  \eqref{e:EP_1} holds under the further assumptions \Tailol \ and $\U_{R_0}(\psi, \beta_3, C_U')$  for the moment. Then
	since $\mathrm{TJ}_{R_0}(\phi,\le,U)$ and  $\U_{R_0}(\phi, \beta_2, C_U)$  are contained in \as, we deduce that under \as \ only, for all $x \in U$, $0< r<3^{-1}\big(R_0 \land (C_1\updelta_U(x))\big)$ and $t\in(0,\phi(r)]$, since $\vt_1(t, C_9r) \le \phi^{-1}(t)$  and $\U_{R_0}(\phi, \beta_2, C_U)$ holds, by \eqref{e:exp1},
	\begin{equation*}
		\P^x(\tau_B \le t) \le C_7 \bigg( \frac{t}{\phi(r)} + \exp \Big(- \frac{C_8 r}{\phi^{-1}(t)} \Big) \bigg)  \le \frac{(C_7 + \beta_2^{\beta_2} C_U C_8^{-\beta_2}) t}{\phi(r)}.
	\end{equation*}
	Hence, we obtain \eqref{e:EP_0}  (since it holds trivially in case $t>\phi(r)$).  Therefore, it suffices to prove \eqref{e:EP_1} under the further assumptions  \Tailol \ and $\U_{R_0}(\psi, \beta_3, C_U')$.

	Now, we prove \eqref{e:EP_1}. Choose any $x \in U$, $0< r<3^{-1}\big(R_0 \land (C_1\updelta_U(x))\big)$ and $t>0$. Set $B:=B(x,r)$  and let $C_3,C_4,C_5$ be the constants in Lemmas \ref{l:4.20}(i) and \ref{l:4.21}(i).  Let $C_9=C_5/(8C_3)$. If $\vt_1(t, C_9r) \ge r/4$, then by taking $C_7$ larger than $e^{4C_8}$, we obtain \eqref{e:EP_1}. Besides, if $C_9r < \Up_1(t)$,  then by \eqref{e:Up} and \eqref{e:theta1}, it holds $\exp \big(- C_8 r/\vt_1(t,C_9r) \big) \ge \exp \big( - C_8C_9^{-1}\big)$. Thus, by taking $C_7$ larger than $e^{C_8/C_9}$, we get \eqref{e:EP_1}. Therefore, we suppose that $\vt_1(t, C_9r) < r/4$ and $C_9r  \ge \Up_1(t)$.

	Let $\rho \in [\psi^{-1}(t), r/4)$ be a constant chosen later.
	Let $Y^1:= X^{(\rho)}$ be a $\rho$-truncated process. 
	Recall that the L\'evy measure of $Y^1$ is given by $J_0(x,dy):=J(x,dy)\1_{\{|x-y| < \rho \}}$. 
	Set 
	\begin{equation*}
		J_1(x,dy):=J(x,dy)\1_{\{\rho \le |x-y| < r/4 \}} \quad \text{and} \quad J_2(x,dy)\1_{\{|x-y|  \ge r/4 \}}.
	\end{equation*}
	Let $(\xi^1_n)_{n \ge 1}$ be i.i.d. exponential random variables with rate parameter $1$ independent of $Y^1$, and define an additive functional of $Y^1$ as
	\begin{equation}\label{e:H1}
		H_{1,u} := \int_0^u J_1(Y^1_s,M_\partial)ds.
	\end{equation}
	Set $S_{1,1}= \inf\{u >0 : H_{1,u}  \ge \xi_1^1 \}$. Let $Y^2_s=Y^1_s$ for $0\le s<S_{1,1}$, and then define $Y^2_{S_{1,1}}$ with law $J_1( Y^2_{S_{1,1}-}, \cdot)/J_1( Y^2_{S_{1,1}-}, M_\partial)$. Repeat this procedure with the new starting point $(S_{1,1}, Y^2_{S_{1,1}})$  to determine the next random time $S_{1,2}:=\inf \{u>0: H_{1,u} \circ \theta^{Y^2}_{S_{1,1}} \ge \xi_2^1\}$, where $\theta^{Y^2}$ denotes the shift operator with respect to $Y^2$, and behaviors of $Y^2_u$ for $ u \in (S_{1,1}, S_{1,1} + S_{1,2}]$. Let 
	$
	T_{1,0}=0$ and $T_{1,n}:= \sum_{k=1}^n S_{1,k}$ for $n \ge 1$.
	By iterating this procedure,  one can construct a process $Y^2_u$ for  $u \in [0, \,\lim_{n \to \infty} T_{1,n}] \cap [0, \infty)$.  	Next, by setting $Y^2$ as the new underlying process, using the additive functional (of $Y^2$)
	\begin{equation}\label{e:H2}
		H_{2,u} := \int_0^u J_2(Y^2_s,M_\partial)ds
	\end{equation}
	and i.i.d. exponential random variables $(\xi^2_n)_{n \ge 1}$ with rate parameter $1$ independent with all of $Y^1, Y^2$ and $(\xi^1_n)_{n \ge 1}$, one can construct a process $Y^3$, $S_{2,n}$ and $T_{2,n}$ by the same way. 
	According to \cite{Me75}  (see also \cite[Section 3.1]{BGK09}), the joint law of $\big((Y^2_u, Y^3_u) : u \in [0, \, \lim_{n \to \infty} (T_{1,n} \wedge T_{2,n})\big)$ is the same as the one for $\big((X^{(r/4)}_u, X_u) :  u \in [0, \, \lim_{n \to \infty} (T_{1,n} \wedge T_{2,n}) \big)$.
	
	We note that $C_1\updelta_U(z) \ge C_1(\updelta_U(x)-r) \ge 2r$ for all $z \in B$. Thus, since \Tailol \ and $\U_{R_0}(\psi, \beta_3, C_U')$ hold, we have that  
	\begin{equation}\label{e:J1}
		\sup_{z \in B} J_1(z, M_\partial)  \le \sup_{z \in B} J\big(z, M_\partial \setminus B(z, \rho)\big) \le \frac{c_1}{\psi(\rho)},
	\end{equation}
	\begin{equation}\label{e:J2}
		\sup_{z \in B} J_2(z, M_\partial) \le \sup_{z \in U} J\big(z, M_\partial \setminus B(z, r/4)\big) \le \frac{c_1}{\psi(r/4)} \le \frac{4^{\beta_3}c_1C_U'}{\psi(r)}=:\frac{c_2}{\psi(r)}.
	\end{equation}
	Let $\wt \tau_B:= \inf\{u>0: Y^3_u \in M_\partial \setminus B \}$. In  $[0, \wt \tau_B)$,  almost surely, there can be  only finitely many extra jumps added in the construction of $Y^3$ from $Y^1$. Thus, for all $x \in U$,  $\lim_{n \to \infty} (T_{1,n} \wedge T_{2,n}) \ge \wt \tau_B$ (a.s.) under $\P^x$, and hence $\P^x(\tau_B \le t) =\P^x(\wt \tau_B \le t) $. Therefore, we have
	\begin{align}\label{e:EPdecomp}
		\P^{x}(\tau_B \le t)  &= \P^x(T_{2,1} \le \wt\tau_B \le t)+ \P^x(T_{1,2} \le \wt\tau_B \le t, \, \wt\tau_B<T_{2,1})+ \P^x(\wt\tau_B \le t, \, \wt\tau_B <  T_{1,2} \land T_{2,1} ) \nn\\
		&=: P_1 + P_2 + P_3.
	\end{align}
	First, by  \eqref{e:H2} and \eqref{e:J2}, since $Y^2_s = Y^3_s \in B$ for all $s \in [0, T_{2,1} \wedge \wt \tau_B)$, we have that
	\begin{equation}\label{e:I1}
		P_1 = \P^x \big(H_{2,S_{2,1}} \ge \xi^2_1, \, S_{2,1} \le \wt\tau_B \le t\big) \le  \P^x \big(\frac{c_2 S_{2,1}}{\psi(r)}\ge \xi^2_1, \, S_{2,1}  \le t\big)\le  \P^x \big( \frac{c_2t}{\psi(r)} \ge \xi^2_1\big) \le \frac{c_2 t}{\psi(r)}.
	\end{equation}
	Next, by a similar way, we get from \eqref{e:H1} and \eqref{e:J1} that 
	\begin{equation}\begin{split}\label{e:EP2-1}
			P_2 &\le \P^x\big(H_{1, S_{1,1}} + H_{1, S_{1,2}} \circ \theta^{Y^2}_{S_{1,1}}  \ge \sum_{k=1}^2 \xi_k^1,  \;\; S_{1,1} + S_{1,2}\le \wt\tau_B \le t, \, \wt\tau_B<S_{2,1}\big) \\
			&\le  \P^x\big( \frac{c_1t}{\psi(\rho)} \ge \sum_{k=1}^2 \xi_k^1 \big)= \int_0^{c_1t/\psi(\rho)} l e^{-l} dl \le \Big(\frac{c_1t}{\psi(\rho)}\Big)^2.
		\end{split}
	\end{equation}
	Lastly, we bound $P_3$. To do this, we set 
	\begin{equation*}
		g_k: = \sup \Big\{  d(Y^3_{\wt \tau_B\land T_{1,k-1}}, Y^3_{u-}) :u \in (\wt \tau_B \land  T_{1,k-1}, \, \wt \tau_B \land   T_{1,k} ]\, \Big\}, \quad k=1,2.
	\end{equation*}
	Note that on the event $\{\wt\tau_B \le t, \, \wt\tau_B <  T_{1,2} \wedge T_{2,1} \}$, by the triangle inequality, it holds  that
	\begin{align*}
		r \le d(x, Y^3_{\wt \tau_B}) &\le d(Y^3_{T_{1,0}}, Y^3_{\wt \tau_B\land T_{1,1}-}) + d(  Y^3_{\wt \tau_B \land  T_{1,1}-},  Y^3_{\wt \tau_B \land T_{1,1}}) + d(  Y^3_{\wt \tau_B\land T_{1,1}},  Y^3_{\wt \tau_B})\le g_1 + r/4  + g_2,
	\end{align*}
	where in the third inequality above, we used the fact  that the jump size at time $T_{1,1}-$ can not be larger than $r/4$ by the definition of $J_1$, and  $\wt \tau_B = \wt \tau_B \land  T_{1,2}$ on the event $\{\wt\tau_B \le t, \, \wt\tau_B <  T_{1,2} \land  T_{2,1} \}$. Thus, we see that on the event $\{\wt\tau_B \le t, \, \wt\tau_B <  T_{1,2} \wedge T_{2,1} \}$, 
	\begin{equation}\label{e:EP3}
		\max \{g_1, g_2\} >r/4, \quad \;\; \P^x\text{-a.s.}
	\end{equation}
	Besides, since $3^{-1}C_1 \updelta_U(z) \ge 3^{-1}C_1 (\updelta_U(x)-5r/4) > r/4$ for all $z \in B(x, 5r/4)$, according to Markov inequality and  Lemma \ref{l:4.21}(i) we have that,
	\begin{align}\label{e:EP0}
		&\sup_{z \in B(x,5r/4)}\P^z(\tau_{B(z,r/4)}^{(\rho)} \le t) = \sup_{z \in B(x,5r/4)}\P^z \Big( \exp\big(-C_3\tau_{B(z,r/4)}^{(\rho)}/\phi(\rho)\big) \ge \exp \big(- C_3 t/ \phi(\rho)\big) \Big)  \nn\\
		&\le \exp \big(  C_3t/\phi(\rho) \big) \sup_{z \in B(x,5r/4)}\E^z\Big[ \exp\big(-C_3\tau_{B(z,r/4)}^{(\rho)}/\phi(\rho)\big) \Big] \le C_4 \exp \big(- \frac{C_5r}{4\rho } + \frac{C_3  t}{\phi(\rho)} \big).
	\end{align}
	
	In view of  the Meyer's construction, on the event $\{\wt \tau_B \le t, \,  \wt\tau_B < T_{1,2} \land T_{2,1}\}$,  for each $k=1,2$, the law of $\big(Y^3_{u}: u \in [\wt \tau_B\land T_{1,k-1} , \wt \tau_B \land T_{1,k})\big)$ is the same as the one for $\big(Y^1_u: u \in [0,  \wt \tau_B \land T_{1,k}- \wt\tau_B \land T_{1,k-1} )\big)$ starting from $Y^3_{\wt \tau_B \land T_{1,k-1}}$.
	Therefore, by   \eqref{e:EP3}, \eqref{e:EP0} and  the strong Markov property, we obtain
	\begin{align}\label{e:I3}
		P_3 &= \P^x\big(  \max\{g_1, g_2\} >r/4, \; \wt\tau_B \le t, \;  \wt\tau_B < T_{1,2} \land T_{2,1}\big) \nn\\[3pt]
		& \le \P^x\big(\min_{1 \le k \le 2} \tau^{(\rho)}_{B(Y^3_{\wt \tau_B \land T_{1,k-1}}, r/4)}  \circ \theta^{Y^3}_{\wt \tau_B\land T_{1,k-1}} \le t, \; \wt\tau_B \le t, \;  \wt\tau_B < T_{1,2} \land T_{2,1} \big) \nn\\
		& \le 2 \sup_{z \in B(x,5r/4)} \P^z \big(\tau^{(\rho)}_{B(z, r/4)}  \le t\big) \le 2C_4 \exp \Big(-\frac{C_5 r}{4\rho} + \frac{C_3 t}{\phi(\rho)} \Big),
	\end{align}
	where $\theta^{Y^3}$ denotes the shift operator with respect to $Y^3$.
	In the second inequality above, we used the 
	the fact that $Y^3_{\wt \tau_B\land T_{1,k-1}} \in B(x, r + r/4)$ for each $k =1,2$. 
	
	Eventually, by combining all of \eqref{e:EPdecomp}, \eqref{e:I1}, \eqref{e:EP2-1} and \eqref{e:I3}, we conclude that
	\begin{equation}\label{e:EPdecomp2}
		\P^x(\tau_B \le t) \le \frac{c_2t}{\psi(r)} + \Big(\frac{c_1t}{\psi(\rho)}\Big)^2 + 2C_4 \exp \Big(-\frac{C_5 r}{4\rho} + \frac{C_3 t}{\phi(\rho)} \Big).
	\end{equation}
	For the last step of the proof, we first assume $t\psi(r/4) \le \psi(\phi^{-1}(t))^2$. Then  by Lemma \ref{l:theta1}, we see
	\begin{equation*}
		\frac{C_3t}{\phi(\vt_1(t, C_9r))} \le \frac{C_3C_9r}{\vt_1(t, C_9r)} = \frac{C_5r}{8\vt_1(t, C_9r)}.
	\end{equation*}
	Further, since  $\U_{R_0}(\psi, \beta_3, C_U')$ holds, 
	\begin{equation*}
		\Big(\frac{t}{\psi(\vt_1(t, C_9r))}\Big)^2 \le \Big(\frac{t}{\psi(\phi^{-1}(t))}\Big)^2 \le \frac{t}{\psi(r/4)} \le  \frac{C_U'4^{\beta_3} t}{\psi(r)}.
	\end{equation*}
	Thus, we get \eqref{e:EP_1} by substituting $\rho = \vt_1(t, C_9r) \in [\psi^{-1}(t), r/4)$ in \eqref{e:EPdecomp2}.
	
	Otherwise, if $t\psi(r/4) > \psi(\phi^{-1}(t))^2$, then we substitute $\rho = \psi^{-1}(t^{1/2} \psi(r/4)^{1/2}) \in (\phi^{-1}(t), r/4)$ in \eqref{e:EPdecomp2}. Since $\U_{R_0}(\psi, \beta_3, C_U')$ holds, we have $(t/\psi(\rho))^2 = t/\psi(r/4) \le C_U'4^{\beta_3}t/\psi(r)$ and
	\begin{align*}
		\exp \big(-\frac{C_5 r}{4\rho} + \frac{C_3 t}{\phi(\rho)} \big) \le e^{C_3}(2\beta_3/C_5)^{2\beta_3} \Big(\frac{\rho}{r/4}\Big)^{2\beta_3} \le c_3\Big( \frac{ \psi(\rho)}{\psi(r/4)}  \Big)^2= \frac{c_3t}{\psi(r/4)} \le   \frac{c_4t}{\psi(r)}.
	\end{align*}
	In the first inequality above, we used the fact that $e^{-x} \le (2\beta_3)^{2\beta_3}x^{-2\beta_3}$ for all $x>0$.
	This completes the proof for (i).

 (ii) We follow the proof of (i) using the same notations. It still suffices to prove \eqref{e:EP_1} under the further assumptions. Choose any  $x \in M$, $r>R_\infty'\ep(x)^{\up_1}$ and $t > \psi(2r^{\up/\up_1})$. Let $C_3,C_4,C_5$ be the constants chosen by Lemmas \ref{l:4.20}(ii) and \ref{l:4.21}(ii). We can still assume that $\vt_1(t, C_9r)<r/4$ and $C_9r \ge \Up_1(t)$ where $C_9=C_5/(8C_3)$.
	For each $\rho \in [\psi^{-1}(t), r/4) \subset [2r^{\up/\up_1}, r/4)$, one can construct two processes $Y^2$ and $Y^3$ from a $\rho$-truncated process $Y^1:=X^{(\rho)}$ so that  the  joint law $\big((Y^2_u, Y^3_u): u \in [0, \infty)\big)$ is the same as $\big((X^{(r/4)}_u, X_u): u \in [0, \lim_{n \to \infty}(T_{1,n} \land T_{2,n}))\big)$ by the similar way to the one given in (i). In view of \eqref{e:ballstable}, we obtain \eqref{e:J1} and \eqref{e:J2}. Hence,  \eqref{e:I1} and \eqref{e:EP2-1}  are still valid. Moreover, we get \eqref{e:EP0} from Lemma \ref{l:4.21}(ii). In the end, we obtain \eqref{e:EP_1} by the same arguments. \qed

	Recall that $\zeta$ denotes the lifetime of $X$.

	\begin{cor}\label{c:conserv}
		Suppose that \ass \ and $\L^{R_\infty}(\phi, \beta_1, C_L)$ hold. Then $X$ is conservative in the sense that 
		 $\P^x(\zeta=\infty)=1$ for all $x \in M$.
	\end{cor}
	\pf  Fix $\up_1\in (\up,1)$. Choose any  $x \in M$ and $T>0$. By Proposition \ref{p:EP}(ii) and $\L^{R_\infty}(\phi, \beta_1, C_L)$, we have that for all $\phi^{-1}(t) \ge 2R_\infty'^{\,\up/\up_1}\ep(x)^\up + \phi^{-1}(T)$, 
	\begin{equation*}
		\P^x(\zeta \le T) \le \P^x(\tau_{B(x,\, (\phi^{-1}(t)/2)^{\up_1/\up})} \le t) \le c_1(\phi^{-1}(t))^{-\beta_1(\up_1-\up)/\up}.
	\end{equation*}
	By taking $t \to \infty$, we get $\P^x(\zeta \le T)=0$. Since this holds with arbitrary $x \in M$ and  $T>0$, we get the result. \qed

	We turn to establish lower estimates on $\P^x(\tau_{B(x,r)} \le t)$. 
	\begin{prop}\label{p:EPL0}
		(i) Suppose that $\U_{R_0}(\psi, \beta_3, C_U')$ and   \Tailog \  hold with $C_1$. Then, for every $a \ge 1$, there exists $c>0$ such that for all $x \in U$ and $0<r<3^{-1}\big((a^{-1}R_0) \land (C_1\updelta_U(x))\big)$,
		\begin{equation}\label{e:EPdown0}
			\P^x(\tau_{B(x,r)}\le t) \ge \frac{c t}{\psi(r)}, \quad\;\; 0<t \le \psi(ar).
		\end{equation}
		
		\noindent (ii) Suppose that $\U^{R_\infty}(\psi, \beta_3, C_U')$ and  \Tailig \  hold. Then, for every $a\ge 1$, there exists $c>0$ such that  \eqref{e:EPdown0} holds for all  $x \in M$ and $r>(2R_\infty)^{1/(1-\up)} \ep(x)^\up$.
	\end{prop}
	\pf (i) Choose any $x \in U$, $0<r<3^{-1}\big((a^{-1}R_0) \land (C_1\updelta_U(x))\big)$ and denote $B:=B(x,r)$. Set $T_r:=\inf\{s>0: d(X_{s-},X_s) \ge 2r\}$. Then  $T_r \ge \tau_B$, $\P^x$-a.s. Hence, by the Meyer's construction and \Tailog,  we have that,  for all $t>0$  (cf. \eqref{e:H2}, \eqref{e:J2} and \eqref{e:I1}),
	\begin{equation*}
		\P^x(\tau_B > t)  =  \P^x(T_r \land \tau_B > t) \le \P^x(c_1t/\psi(2r)< \xi_1) = e^{-c_1t/\psi(2r)},
	\end{equation*}
	where $\xi_1$ is an exponential random variable with rate parameter $1$ independent of $X$. Hence, since $\U_{R_0}(\psi, \beta_3, C_U')$ holds, we get that for all $0<t \le \psi(ar)$,
	\begin{equation*}
		\P^x(\tau_B \le t) \ge 1- e^{-c_1t/\psi(2r)} \ge \frac{c_1t}{\psi(2r)}e^{-c_1\psi(ar)/\psi(2r)} \ge \frac{c_3t}{\psi(r)}.
	\end{equation*}

(ii) Using	a similar calculation to  \eqref{e:ballstable}, we get the result by the same way as the one for (i). We omit details here. \qed 
	
	Under the local chain condition (see Definition \ref{d:Ch}), we  obtain lower estimates on $\P^x(\tau_{B(x,r)} \le t)$ which are sharp in view of \eqref{e:EP_1}.
	\begin{prop}\label{p:EPL}
		(i) Suppose that \VRDo, \Cho,  \NDLo \ and $\U_{R_0}(\phi, \beta_2, C_U)$ hold with $C_V$ and $C_2$. 
		Then, there exists a constant $c>0$ such that for each $a>0$, there are constants  $C_{11}, C_{12}>0$ such that for all $x \in U$, $0<r< c\big((a^{-1}R_0) \land ((C_V \land C_2)\updelta_U(x))\big)$ and $0<t \le \phi(a r)$ satisfying $C_{12}r \le \Up_2(t)$,
		\begin{equation}\label{e:EPdown}
			\P^x(\tau_{B(x,r)}\le t) \ge  \exp\Big(-\frac{C_{11} r}{\vt_2(t, C_{12}r)}\Big),
		\end{equation}
	where $C_V, C_2 \in (0,1)$ are the constants in \VRDo \ and \NDLo.		If   \Tailog \ and $\U_{R_0}(\psi, \beta_3, C_U')$ also hold, then there exists a constant $C_{10}>0$ such that for all  $x \in U$, $0<r<c\big((a^{-1}R_0) \land ((C_V \land C_1 \land  C_2)\updelta_U(x))\big)$  and $0<t \le \phi(a r)$,
		\begin{equation}\label{e:EPdown2}
			\P^x(\tau_{B(x,r)}\le t) \ge C_{10} \Big(\frac{ t}{\psi(r)} +  \exp\Big(-\frac{C_{11} r}{\vt_2(t, C_{12}r)}\Big)\Big).
		\end{equation}

		\noindent (ii) Suppose that \VRDi, \Chi, \NDLi \ and  $\U^{R_\infty}(\phi, \beta_2, C_U)$  hold with $\up\in(0,1)$. Let $\up_1 \in (\up,1)$. Then for every $a>0$, there exist constants $R_\infty'' \in [R_\infty, \infty)$, $C_{11}, C_{12}>0$ such that  \eqref{e:EPdown} holds for all $x \in M$, $r> R_\infty'' \ep(x)^{\up_1}$ and $\psi(R_\infty'' r^{\up/\up_1}) < t \le \phi(a r)$ satisfying $C_{12}r \le \Up_2(t)$.
			If  \Tailig \ and  $\U^{R_\infty}(\psi, \beta_3, C_U')$ also hold, then there exists a constant $C_{10}>0$ such that \eqref{e:EPdown2} holds for  all $x \in M$, $r > R_\infty''\ep(x)^{\up_1}$ and $\psi(R_\infty''r^{\up/\up_1}) < t \le \phi(a r)$.
	\end{prop}
	\pf  (i) Let $\ell>1$, $A \ge 1$ and $\eta \in (0,1)$ be the constants in  \eqref{RVD}, \Cho \ and \NDLo, respectively. Then we set $c':=3^{-1}\ell^{-1}A^{-1}\eta  \in (0,1/3)$ and $C_{12}:=8 C_U (1+a)^{1+\beta_2}  \ell  A\eta^{-1}$.  Choose any $x \in U$, $0<r< c'\big(R_0 \land ((C_V\land C_2)\updelta_U(x))\big)$   and $t \in (0, \phi(ar)]$. Write $B:=B(x,r)$.

	We first prove \eqref{e:EPdown}.  Set $\vt:=\vt_2(t,C_{12}r) \in [\psi^{-1}(t), \phi^{-1}(t)]$ and  $n_{0}:=\lceil C_{12}r/(2\vt) \rceil$. Then we see that $n_{0} \ge \lceil C_{12}r/(2\phi^{-1}(t)) \rceil  \ge \lceil C_{12} / (2a) \rceil \ge 4$. 
	Since $\Up_{2}(t) \ge C_{12}r \ge ar \ge  \phi^{-1}(t) \ge \Up_1(t)$, by the definition of $\vt_2$, we have  $t\vt /\phi(\vt) \ge  C_{12}r$. Thus,  since $C_{12}r/(2\vt) \le n_0 \le C_{12}r/\vt$, we get
	\begin{equation}\label{e:EPchain}
		\phi^{-1}(t/n_0) \ge \phi^{-1}(t\vt/(C_{12} r)) \ge \vt \ge 2^{-1}n_0^{-1}C_{12} r > 4 n_0^{-1}\ell  A \eta^{-1} r.
	\end{equation}
	On the other hand, since $\U_{R_0}(\phi, \beta_2, C_U)$ holds, 
	\begin{equation*}
		n_0\phi(r) \ge \frac{C_{12}r\phi(r)}{2\phi^{-1}(t)} \ge  \frac{C_{12}  \phi(r)}{2a} \ge C_U (1+a)^{\beta_2}\phi(r) \ge \phi(ar) \ge t \quad \text{so that} \;\; \phi^{-1}(t/n_0) \le r.
	\end{equation*}
	
	By \eqref{RVD}, we can choose $y \in B(x, 2\ell r) \setminus B(x,2r)$. By \Cho,
since $2\ell r< R_0$,   there exists a sequence $(z_i)_{0 \le i \le n_{0}} \subset M$ such that $z_0=x$, $z_{n_0}=y$ and $d(z_{i-1},z_i) \le 2 \ell A r/n_0$ for all $1 \le i \le n_0$. Note that for all $1 \le i \le n_0$, we have $d(x, z_i) \le \sum_{j=1}^i d(z_{j-1}, z_j) \le 2 \ell Ar$ and hence 
	\begin{equation*}
		(C_V \land C_2)\updelta_U(z_i) \ge (C_V \land C_2)(\updelta_U(x)-2\ell Ar) \ge \eta^{-1}\ell A r  >\eta^{-1}r \ge  \eta^{-1} \phi^{-1}(t/n_0).
	\end{equation*}
	Eventually, by the semigroup property, since  \NDLo \ and \VRDo \ hold, we obtain
	\begin{align}\label{e:EPLOW}
		&\P^x(\tau_B \le t)  \ge  \P^x\big(X_t \in B(y,r)\big) \ge  \P^x\big(X_t \in B(y,\phi^{-1}(t/n_0))\big)\nn\\
		&\ge \int_{B(z_1, 2^{-1}\eta\phi^{-1}(t/n_0))}...\int_{B(z_{n_0}, 2^{-1}\eta\phi^{-1}(t/n_0))} p^{B(x, \eta^{-1}\phi^{-1}(t/n_0))}(t/n_0, x, u_1) \nn\\
		& \hspace{1.8cm} \times p^{B(z_1, \eta^{-1}\phi^{-1}(t/n_0))}(t/n_0,u_1,u_2) \dots p^{B(z_{n_0-1}, \eta^{-1}\phi^{-1}(t/n_0))}(t/n_0,u_{n_0-1},u_{n_0}) \, du_{n_0}...du_{1}\nn\\
		& \ge \prod_{i=0}^{n_0-1}\left(\frac{c_1V(z_i, 2^{-1}\eta \phi^{-1}(t/n_0))}{V(z_i, \eta^{-1}\phi^{-1}(t/n_0))}\right) \ge e^{-c_2 n_0} \ge e^{-c_2C_{12}r/\vt}.
	\end{align}
	In the third inequality above, we used the fact that since \eqref{e:EPchain} holds, for all $1 \le j \le n_0$ and $u_j \in B(z_j, {2^{-1}\eta \phi^{-1}(t/n_0)})$,  it holds  that
	$$d(u_{j}, z_{j-1}) \vee d(u_{j}, z_{j}) \le d(u_j, z_j) +  d(z_{j-1}, z_{j}) \le 2^{-1}\eta \phi^{-1}(t/n_0) + 2 \ell Ar/n_0  < \eta \phi^{-1}(t/n_0).$$ The proof of \eqref{e:EPdown} is complete.

	Now, we also assume \Tailog \ and  $\U_{R_0}(\psi, \beta_3, C_U')$, and then prove \eqref{e:EPdown2}.
	 By Proposition \ref{p:EPL0} and \eqref{e:EPdown}, since $\psi \ge \phi$,
	 \eqref{e:EPdown2} holds when $C_{12}r \le \Up_2(t)$. Let $C_{12}r>\Up_2(t)$. Then $\vt_2(t, C_{12}r) = \psi^{-1}(t)$  by the definition \eqref{e:theta2}. Thus, in view of  Lemma \ref{l:exp1}(i) and Proposition \ref{p:EPL0}, since $\psi^{-1} \le \phi^{-1}$, we can deduce  that \eqref{e:EPdown2} holds. 
	
	\smallskip
	
 (ii) We follow the proof of (i). Set  $C_{12}:=8C_U(1+a)^{1+\beta_2}\ell A\eta^{-1}$ and $R_\infty'':=(2\ell)^\up C_{12}R_\infty$ with the constants $\ell, A, \eta$ in \eqref{RVD}, \Chi \ and  \NDLi. Choose any  $x \in M$, $r>R_\infty''\ep(x)^{\up_1}$ and $\psi(R_\infty''r^{\up/\up_1})<t \le \phi(ar)$. Let $\vt:=\vt_2(t,C_{12}r)$ and $n_0:=\lceil C_{12}r/(2\vt)\rceil$. 
	By \eqref{RVD}, there exists $y \in B(x, 2\ell r) \setminus B(x,2r)$. Since \Chi \ holds and 
	\begin{equation}\label{e:Chi}
		n_0R_\infty(\ep(x) \vee \ep(y))^\up \le \frac{C_{12}R_\infty r}{2\psi^{-1}(t)}(\ep(x) + 2 \ell r)^\up \le \frac{r^{1-\up/\up_1}}{2^{1+\up} \ell^\up} ( r^{\up/\up_1}+ (2\ell)^\up r^\up ) <2r,
	\end{equation}
	there is a sequence $(z_i)_{0 \le i \le n_0} \subset M$  such that $z_0=x$, $z_{n_0}=y$ and $d(z_{i-1},z_i) \le 2 \ell A r/n_0$ for all $1 \le i \le n_0$. Observe that \eqref{e:EPchain} is still valid. Hence, we have that for all $0 \le i \le n_0$,
	\begin{align*}
		2^{-1}\eta \phi^{-1}(t/n_0) &> 2 n_0^{-1}\ell A r \ge 4 C_{12}^{-1} \ell A \vt\ge 4 C_{12}^{-1} \ell A \psi^{-1}(t) > 4 \ell A R_\infty r^{\up/\up_1}\\
		&> R_\infty \ep(x)^\up + (2\ell A)^\up R_\infty r^\up > R_\infty(\ep(x)+2\ell Ar)^\up \ge R_\infty \ep(z_i)^\up.
	\end{align*}
	Then using \NDLi \ and \VRDi, we  get \eqref{e:EPdown} by repeating the calculation \eqref{e:EPLOW}. 
	
	When  \Tailig \ and  $\U^{R_\infty}(\psi, \beta_3, C_U')$ also hold, using Lemma \ref{l:exp1}(ii) and Proposition \ref{p:EPL0}, we get \eqref{e:EPdown2} as in (i). \qed
	
	\subsection{ Estimates on oscillation of parabolic functions and zero-one law}\label{s:tp}
	In this subsection, we establish  estimates on oscillation of  bounded parabolic functions for large distances under Assumption \ref{assinf2} (Proposition \ref{p:phr}). Then, as an application, we obtain a zero-one law for shift-invariant events in Proposition \ref{p:law01}.
	
	\medskip

	{\it  Throughout  this subsection, we always  assume that  $\L^{R_\infty}(\phi, \beta_1, C_L)$ and  \NDLi \  hold. }
	
	\medskip
	\noindent With the constants $\beta_1>0$ and $C_L \in (0,1]$ in $\L^{R_\infty}(\phi, \beta_1, C_L)$, we set  
	\begin{equation}\label{e:defvk}
		\vk:=(2/C_L)^{\beta_1} \in [1, \infty).
	\end{equation}
	Then by  $\L^{R_\infty}(\phi, \beta_1, C_L)$, we have
	\begin{equation}\label{e:doublevk}
		\phi(\vk r) \ge 2\phi(r) \quad \text{for all} \;\; r \in (R_\infty, \infty).
	\end{equation}
	
	With the constant $\vk$ in \eqref{e:defvk}, and  $\eta \in (0,1)$ in \NDLi,  we define open cylinders $\sQ$, $\sQ^-$, $\sQ^+$ and $\sI$ as follows: for $x \in M$, $r>0$ and $t \ge  \phi(2  \vk \eta r)$,
	\begin{align}\label{e:defQ}
		\sQ(t,x,r) &:= (t - \phi(2\vk \eta r), \, t) \times B(x, 2\vk r), \quad \sQ^-(t,x,r):= (t-\phi(2\vk \eta r), \, t- \phi(\vk \eta r)) \times B(x,\eta^2  r)  \nn\\[1pt]
		\sI(t,x,r)&:=(t - \phi(2 \eta r), \, t ) \times B(x, \vk r),\qquad\, \sQ^+(t,x,r):= (t-\phi( \eta r), \, t) \times B(x,\eta^2  r).
	\end{align}
	Note that $\sQ^-(t,x,r) \cup \sQ^+(t,x,r) \cup \sI(t,x,r) \subset \sQ(t,x,r)$ and $\sQ^+(t,x,r) \cap \sQ^-(t,x,r)=\emptyset$.

	Denote by $Z = (Z_s)_{s \ge 0} = (V_s,X_s)_{s \ge 0}$ a time-space stochastic process corresponding to $X$ with $V_s:= V_0 - s$.  For $D \in \sB ([0, \infty) \times M)$, we write $\sigma^Z_D:=\inf\{t>0: Z_t \in D\}$ and $\tau^Z_D:=\inf\{t>0:Z_t \in [0, \infty) \times M_\partial \setminus D\}$. 
	
	Let $dt \otimes \mu$ be the product measure of the Lebesgue measure on $[0, \infty)$ and $\mu$ on $M$.

	\begin{lem}\label{l:3.7}
		There exists a constant $c_0>0$ such that for all $x \in M$, $r>\eta^{-1}R_\infty \ep(x)^\up$, $t \ge \phi(2\vk \eta r)$ and any compact set $D \subset \sQ^-(t,x,r)$,
		$$ \inf_{w \in \sQ^+(t,x,r)}	 \P^w\big(\sigma^Z_D < \tau^Z_{\sQ(t,x,r)}\big) \ge  \frac{c_0  (dt \otimes \mu)(D)}{\big(\phi(2 \vk \eta r)- \phi( \vk \eta r) \big)V(x, 2 \vk r)}, $$
		where the constant $\vk$ and cylinders $\sQ^-, \sQ^+$ and $\sQ$ are defined as  \eqref{e:defvk} and \eqref{e:defQ}.
	\end{lem}
	\pf 
	For $u>0$, we set  $D_u:= \{ y \in M : (u,y) \in D\}$. Write $\tau_r:= \tau_{\sQ(t,x,r)}^Z$.
	We first observe that, by the fact that $D \subset \sQ^-(t,x,r)$,  for  any $w = (s,z) \in \sQ^+(t,x,r)$,
	\begin{align}
		&\big(\phi(2 \vk \eta r)  - \phi( \vk \eta r) \big)  \P^w (\sigma^Z_D < \tau_r)  \ge \big(\phi(2 \vk \eta r)- \phi( \vk \eta r) \big)  \P^w\bigg( \int_0^{\tau_r} \1_D (s-u, X_u) du > 0 \bigg)  \nn\\
		& \ge \int_0^{\phi(2 \vk \eta r)- \phi(\vk \eta r) } \P^{z} \bigg( \int_0^{\tau_r} \1_D (s-u, X_u) du > a \bigg) da=  \E^{z} \bigg[ \int_0^{\tau_r}  \1_D(s-u,X_u)du  \bigg]\nn\\
		&=\int_{\phi( \vk \eta r) - (t-s)   }^{\phi(2 \vk \eta r) - (t-s)}  \P^z\big(X_u^{B(x, 2\vk r)} \in D_{s-u} \big) du  = \int_{\phi( \vk \eta r)}^{\phi(2 \vk \eta r)}   \int_{D_{t-u}} p^{B(x,2\vk r)} \big(u - (t-s) ,z,y\big) \mu(dy) du.\label{e:l3.14}
	\end{align}
	By the definition of $\sQ^+$,  $t-s \in (0, \phi( \eta r))$. Thus,  by \eqref{e:doublevk}, we see that for all $u \in [\phi(\vk \eta r), \phi(2\vk \eta r)]$,  $ 2 \vk \eta r \ge  \phi^{-1}(u-t+s) \ge \phi^{-1}( \phi(\vk \eta r) -\phi( \eta r)) \ge \eta r$. 
	Since \NDLi \ holds and $z \in B(x,\eta^2r)$, it follows that for all $u \in [\phi(\vk \eta r), \phi(2 \vk \eta r)]$ and $y \in B(x, \eta^2 r )$,
	\begin{align*}
		&p^{B(x,2 \vk r)}(u - (t-s),z,y) \ge p^{B(x,\eta^{-1}\phi^{-1}(u-t+s))}(u -t+s,z,y) \ge  \frac{c_1}{V(x,2 \vk  r)}. 
	\end{align*}
	Thus, by combining this with \eqref{e:l3.14}, we get
	\begin{align*}
		\big(\phi(2 \vk \eta r)  - \phi( \vk \eta r) \big)  \P^w (\sigma^Z_D < \tau_r)  \ge  \frac{c_1}{V(x, 2 \vk  r)} \int_{\phi(\vk \eta r)}^{\phi(2\vk \eta r)} \int_{D_{t-u}} \mu(dy)du = \frac{c_1 (dt \otimes \mu)(D)}{V(x,2 \vk r)}.
	\end{align*}
	This yields the desired inequality. \qed

	We say that a Borel measurable function $q(t,x)$ on $[0,\infty) \times M$ is \textit{parabolic} on $D = (a,b) \times B(x,r)$ for the  process $X$, if for every open set $U  \subset D$ with $\overline U \subset D$,  it holds $q(t,x) = \E^{(t,x)} \big[q(Z_{\tau_U^Z})\big]$ for every $(t,x) \in U \cap ( [0,\infty) \times M )$.
	
	In the following proposition, we let $r_\infty$ be the constant defined as \eqref{e:defrinf}.
	
	\begin{prop}\label{p:phr}		
		Suppose that \asss \ holds. Let $\up_1 \in (\up,1)$. Then, there exist constants $c_0, b>0$ such that for all $x \in M$, $r >(8\vk^2 \eta^{-2}r_{\infty})^{\up_1/(\up_1-\up)} \ep(x)^{\up_1}$, $t \ge \phi(2\vk \eta r)$, and any function $q$ which is non-negative in $[t-\phi(2 \vk \eta r),t] \times M$ and parabolic in $\sQ(t,x,r)$, it holds for all $(s_1,y_1), (s_2,y_2) \in \sI(t,x,r)$ that 
		\begin{equation}\label{e:phrinf}
			\big|q(s_1,y_1) - q(s_2,y_2)\big|\le c_0 \Vert q \Vert_{L^\infty([t-\phi(2 \vk \eta r),\,t] \times M; \, dt \otimes \mu)} \Big( \frac{\phi^{-1}(|s_1-s_2|) + d(y_1,y_2) + r^{\up/\up_1}}{r} \Big)^{b},
		\end{equation}
		where the constant $\vk$ and cylinders $\sQ$ and $\sI$ are defined as  \eqref{e:defvk} and \eqref{e:defQ}, respectively.
	\end{prop}
	\pf Before giving the main argument for \eqref{e:phrinf}, we first set up some constants and functions.
	
	According to Lemma \ref{l:3.7}, since \VRDi \ holds, there exists  $c_1 \in (0,1)$ such that for all $z \in M$, $l>\eta^{-1}R_\infty \ep(z)^{\up}$, $u\ge \phi(2 \vk \eta l)$ and any compact set $D \subset \sQ^-(u,z,l)$ satisfying $(dt \otimes \mu)(D) \ge 3^{-1}(dt \otimes \mu)(\sQ^-(u,z,l))$, 
	\begin{equation}\label{e:phr2-1}
		\inf_{w \in \sQ^+(u,z,l)} \P^{w}\big(\sigma_D < \tau^Z_{\sQ(u,z,l)}\big) \ge c_1. 
	\end{equation}
	Besides, by the L\'evy system \eqref{e:system} and  Proposition \ref{p:E}(ii), since $\Tail^{R_\infty}(\phi, \le,\up)$ and $\U^{R_\infty}(\phi, \beta_2, C_U)$ hold, we see that for all $z_1,z \in M$ and $r_\infty( \ep(z) \vee \ep(z_1))^\up <l \le \wt l/2$,
	\begin{align}
		\P^{z_1}\big(X_{\tau_{B(z,l)}} \in M_\partial \setminus B(z, \wt l)\big) &\le \E^{z_1} \Big[ \int_0^{\tau_{B(z,l)}}  J\big(X_s,M_\partial \setminus B(X_s,\wt l/2) \big)ds \Big] \nn \\ &\le  \frac{c_2 \E^{z_1}[\tau_{B(z_1,\, l+d(z,z_1))}]}{\phi(\wt l/2)} \le \frac{c_3\phi(l + d(z,z_1))}{\phi(\wt l)}, \label{e:phr2-2}
	\end{align} 
	for some  $c_3>1$. The first inequality in \eqref{e:phr2-2} is valid since for all $y \in B(z,l)$, it holds  that $R_\infty\ep(y)^\up\le R_\infty \ep(z)^\up+R_\infty l^\up < l/2 +l/2  \le \wt l/2$ because  $l \ge r_\infty> r_\infty^{1-\up}>2R_\infty$. Using the constants $c_1, c_3$ in \eqref{e:phr2-1} and \eqref{e:phr2-2}, we define
	\begin{equation*}
		\xi =  \Big(1-\frac{c_1}{4}\Big)^{1/2}, \qquad  N_1:= \bigg\lceil  \frac{ \log(2^{\beta_1}/C_L) + \log (  (c_1+4c_3) / (c_1\xi))}{\beta_1 \log (2\vk/\eta^2)} \bigg\rceil.
	\end{equation*}
	Choose any $x \in M$, $r>(8\vk^2 \eta^{-2}r_\infty)^{\up_1/(\up_1-\up)} \ep(x)^{\up_1}$ and $t \ge \phi(2 \vk  \eta r)$.  Then we set $l_0:=r$ and 
	\begin{equation}\label{e:defr}
		l_{i}:=2^{-1}\vk^{-1}\eta^2  l_{i-1} \quad \text{ for } \; i \ge 1.
	\end{equation}
	Note that $l_1 > 4 \vk r_\infty r^{\up/\up_1}$. Let $N_2:=\max\{ n\ge 1  :  l_n > 4 \vk r_\infty r^{\up/\up_1}\}$. Since $\L^{R_\infty}(\phi, \beta_1, C_L)$ holds, 
	\begin{equation}\label{e:ratio}
		\phi(4 \vk l_{i+N_1}) \le 2^{\beta_1} C_L^{-1}\Big( \frac{\eta^2}{2\vk}\Big)^{N_1 \beta_1} \phi(2\vk l_{i}) \le   \frac{c_1\xi}{c_1+4c_3}  \phi(2\vk l_{i}), \quad \;\;1 \le i \le N_2-N_1.
	\end{equation}
	For a given $(u, z) \in \sI(t,x,r)$, we set for $i \ge 1$,
	$$Q_i = Q_i(u,z) := \sQ(u,z, l_i), \quad Q^+_i:= \sQ^+(u,z,l_i) \quad \mbox{ and } \quad Q^-_i:= \sQ^-(u,z,l_i).$$ 
	Since  $\phi$ is increasing, by \eqref{e:defr} and \eqref{e:doublevk}, we have $\phi(2 \vk \eta l_1) < \phi(2\eta r) \le \phi(2\vk \eta r) - \phi(2\eta r)$ and $2\vk l_1 < \vk r$.
	Thus,  $\overline{Q_1(u,z)} \subset \sQ(t,x,r)$  and $Q_{i+1} \subset Q_i^+$ for all $i \ge 1$.
	
	\smallskip
	
	Now, we are ready to begin the main step. 
	Without loss of generality, by considering a constant multiple of $q$, we assume  $\Vert q \Vert_{L^\infty([t-\phi(2 \vk \eta r),t] \times M; \, dt \otimes \mu)}  = 1$.
	We claim that for all $(u,z) \in \sI(t,x,r)$,
	\begin{equation}\label{e:phr2-3}
		\sup_{Q_{1+ j N_1}(u,z)}\, q - \inf_{Q_{1 + j N_1}(u,z)}\, q \le \xi^{j}, \quad  \;\; 0 \le j <N_2/N_1.
	\end{equation}
	In the followings, we prove \eqref{e:phr2-3} by induction. 
	Let $a_j = \inf_{Q_{1+jN_1}(u,z)} q$ and $b_j = \sup_{Q_{1+jN_1}(u,z)} q$. First, we have $b_0 - a_0 \le \Vert q \Vert_{L^\infty([t-\phi(2 \vk \eta r),t] \times M; \, dt \otimes \mu)}  = 1$. 
	
	For the next step, we suppose that $b_j - a_j \le \xi^{j}$ for all $0 \le j \le k < \lceil N_2/ N_1 \rceil-1$. 
	Set 
	$$D' := \big\{ z \in Q_{1+k N_1}^- : q(z) \le 2^{-1}(a_k+b_k) \big\}.$$ We may  assume  $(dt \otimes \mu)(D')\ge 2^{-1}(dt \otimes \mu)( Q_{1+kN_1}^- )$, or else using $1-q$ instead of $q$. 
	Let $D$ be a compact subset of $D'$ such that
	$(dt \otimes \mu)(D)\ge 3^{-1}(dt \otimes \mu)( Q_{1+kN_1}^- )$. Note that we have $l_{1+kN_1} \ge l_{N_2} >4 \vk r_\infty r^{\up/\up_1} > 4 \vk \eta^{-1}R_\infty  (r^{1/\up_1} + r)^\up \ge  \eta^{-1}R_\infty  (\ep(x) + \vk r)^\up \ge  \eta^{-1}R_\infty \ep(z)^\up$ since $r_\infty \ge 2\eta^{-2}R_\infty$. Write $\tau_k := \tau^Z_{Q_{1+kN_1}}$. Then by \eqref{e:phr2-1}, since $Q_{1+(k+1)N_1} \subset Q_{1+kN_1}^+$, we have
	\begin{equation}\label{e:phr2-4}
		\inf_{w \in Q_{1+(k+1)N_1}} \P^w (\sigma_D < \tau_k) \ge c_1.
	\end{equation}

	For a given arbitrary $\eps>0$, find $w_1, w_2 \in Q_{1+(k+1)N_1}$ such that $q(w_1) \ge b_{k+1} -\eps/2$ and $q(w_2) \le a_{k+1} + \eps/2$. Since $q$ is parabolic in $\sQ(t,x,r)$ and $\overline{Q_{1+kN_1} \setminus D} \subset \sQ(t,x,r)$,  we  obtain
	\begin{align*}
		b_{k+1} - a_{k+1} &-\eps 
		\le q(w_1) - q(w_2)  = \E^{w_1} \big[ q(Z_{\sigma_D \land \tau_k} )  - q(w_2) \big] \nn\\[6pt]
		&= \E^{w_1} \big[ q(Z_{\sigma_D}) - q(w_2) ; \, \sigma_D < \tau_k  \big]  + \E^{w_1} \big[ q(Z_{\tau_k} ) - q(w_2)  ; \, \sigma_D \ge \tau_k, Z_{\tau_k} \in Q_{1+(k-1)N_1}  \big] \nn\\
		&\quad + \sum_{j=2}^{k+1} \, \E^{w_1} \big[  q(Z_{\tau_k}) - q(w_2) ; \, \sigma_D  \ge  \tau_k, Z_{\tau_k} \in Q_{1+(k-j)N_1} \setminus Q_{1+(k+1 -j)N_1} \big]\nn\\
		&=: K_1 + K_2 +K_3,
	\end{align*}
	where $Q_{1-N_1}:=(t-\phi(2\vk \eta r),t) \times M_\partial$. 
	
	First,  by \eqref{e:phr2-4} and the induction hypothesis, since $D \subset D'$ and $\xi \in (0,1)$, we have
	\begin{align}\label{e:K12}
		K_1 + K_2 &\le ( 2^{-1}(a_k+b_k) - a_k ) \P^{w_1} (\sigma_D <  \tau_{k}) + (b_{k-1} - a_{k-1}) \P^{w_1} (\sigma_D \ge  \tau_{k}) \nn\\[2pt]
		&= 2^{-1}(b_k-a_k ) \P^{w_1} (\sigma_D <  \tau_{k}) + (b_{k-1} - a_{k-1}) \big(1-\P^{w_1} (\sigma_D <  \tau_{k})\big) \nn\\[2pt]
		& \le \xi^{k-1} \big( 2^{-1}\xi \P^{w_1} (\sigma_D < \tau_k)  +  1- \P^{w_1} (\sigma_D < \tau_k)\big) \le \xi^{k-1} (1-2^{-1}c_1).
	\end{align} 
	On the other hand, observe that for each $j \ge 2$, the  process $Z$ can  enter the set $Q_{1+(k-j)N_1} \setminus Q_{1+ (k+1 -j)N_1}$ at time $\tau_k$ only  through a jump. Write $w_1=(u_1, z_1)$ and observe that
	\begin{equation}\label{e:phr_1}
		r_\infty (\ep(z) \vee \ep(z_1))^\up \le r_\infty (\ep(x) +2 \vk r)^\up < r_\infty r^{\up/\up_1} + 2\vk r_\infty r^\up < 4\vk r_\infty r^{\up/\up_1}<l_{N_2} \le l_{1+kN_1}.
	\end{equation}
	By  \eqref{e:phr2-2}, \eqref{e:ratio} and \eqref{e:phr_1}, since $\Vert q \Vert_{L^\infty([t-\phi(2\vk \eta r),t] \times M; \, dt \otimes \mu)} =1$ and $d(z,z_1) \le 2\vk l_{1+kN_1}$, we get
	\begin{align*}
		K_3 &\le  \sum_{j=2}^{k-1} (b_{k-j} - a_{k-j} ) \P^{z_1} \big(X_{\tau_{B(z, 2 \vk l_{1+kN_1})}} \in M_\partial \setminus B(z, 2 \vk l_{1+ (k+1-j)N_1})\big) \\
		& \quad + \P^{z_1} \big(X_{\tau_{B(z, 2 \vk l_{1+kN_1})}} \in M_\partial \setminus B(z, 2 \vk l_{1+N_1})\big)  \\
		& \le c_3\sum_{j=2}^{k} \xi^{k-j} \frac{\phi(4 \vk l_{1+kN_1})}{\phi(2 \vk l_{1+(k+1-j)N_1})}\le c_3\xi^{k-1}\sum_{j=2}^{k}  \Big(\frac{c_1 }{c_1+4c_3}\Big)^{j-1} \le c_3\xi^{k-1}\sum_{j=1}^{\infty}   \Big(\frac{c_1 }{c_1+4c_3}\Big)^{j} =\frac{c_1}{4}\xi^{k-1}.
	\end{align*}
	Since we can choose $\eps$ arbitrarily small, 
 by combining the above with \eqref{e:K12},  we arrive at
	$$b_{k+1} - a_{k+1} \le (K_1+K_2)+K_3 \le \xi^{k-1}\big( 1- 2^{-1}c_1 + 4^{-1}c_1 \big) =\xi^{k+1}.$$
Therefore,  the claim \eqref{e:phr2-3} holds by induction.

	Let $(s_1,y_1), (s_2,y_2) \in \sI(t,x,r)$ be two different pairs such that $s_1 \ge s_2$. Fix  $s_0 \in (s_1, t)$ such that $\phi^{-1}(s_0-s_2) \le 2\phi^{-1}(s_1-s_2) +  d(y_1, y_2)$. 
	 We set 
	\begin{equation*}
		j_0:=\max \big\{j \in \Z \, : \, -1 \le j <N_2/N_1, \; (s_2, y_2) \in Q_{1+jN_1}(s_0,y_1) \big\}.
	\end{equation*}
	Then,  by \eqref{e:phr2-3}, since $\Vert q \Vert_{L^\infty([t-\phi(2 \vk \eta r),t] \times M; \, dt \otimes \mu)}  = 1$, we obtain 
	\begin{equation}\label{e:qregular}
		|q(s_1,y_1)-q(s_2,y_2)| \le 2\xi^{j_0}
		=2\xi^{-1} (2\vk\eta^{-2})^{-(j_0+1)\log (\xi^{-1})/ \log (2\vk\eta^{-2})}.
	\end{equation}
	If $j_0 \ge N_2/N_1-1$, then we see that $ l_{1+(j_0+1)N_1} \le  l_{1+N_2} \le 4 \vk r_\infty r^{\up/\up_1}$.
	Otherwise, if $j_0<N_2/N_1-1$, then since $(s_2, y_2) \notin Q_{1+(j_0+1)N_1}(s_0,y_1)$, it must holds either $\phi^{-1}(s_0-s_2) \ge 2 \vk \eta l_{1+(j_0+1)N_1}$ or $\eta d(y_1,y_2) \ge 2\vk \eta  l_{1+(j_0+1)N_1}$.
	Thus, by  \eqref{e:defr},  whether $j_0 \ge N_2/N_1-1$ or not, we obtain
	\begin{align*}
		\phi^{-1}(s_1-s_2) +  d(y_1, y_2) + 4\vk^2 \eta r_\infty r^{\up/\up_1}  &\ge 2^{-1}\big( \phi^{-1}(s_0-s_2) + \eta d(y_1, y_2) \big)  + 4\vk^2 \eta r_\infty r^{\up/\up_1} \\
		& \ge \vk \eta  l_{1+(j_0+1)N_1} = 2^{-1}\eta^3 r (2\vk\eta^{-2})^{-(j_0+1)N_1}.
	\end{align*}
	Therefore, by  \eqref{e:qregular}, we deduce that
	\begin{equation*}
		\big|q(s_1,y_1)-q(s_2,y_2)\big| \le 2\xi^{-1} 
		\Big(\frac{\phi^{-1}(s_1-s_2) + d(y_1, y_2) + 4\vk^2 \eta r_\infty r^{\up/\up_1}}{2^{-1}\eta^3 r} \Big)^{\log (\xi^{-1}) / (N_1 \log (2\vk\eta^{-2}))}
	\end{equation*}
The proof is complete. \qed

	An event $G$ is called \textit{shift-invariant}  if $G$ is a tail event (i.e. $\cap_{t>0}^\infty \sigma(X_s:s>t)$-measurable), and $\P^y(G)= \P^y(G \circ \theta_t)$ for all $y \in M$ and $t>0$.

	\begin{prop}\label{p:law01}
		Suppose that \asss \ holds. Then, for every shift-invariant event $G$, it holds either $\P^z(G)=0$ for all $z \in M$ or else $\P^z(G) = 1$ for all $z \in M$.
	\end{prop}
	\pf  
	Fix $\up_1 \in (\up,1)$, $x_0 \in M$ and choose any $\eps \in (0,1)$. By  Propositions \ref{p:E}(ii)  and \ref{p:EP}(ii), there exist constants $T_0, C_6>1$ such that  for all $t_0>T_0$,
	\begin{equation}\label{e:law01_inf-1}
		\P^{x_0} \big( \sup_{s \le t_0} d(x_0,X_s) > \phi^{-1}(\eps^{-1}C_6 t_0)\big) =\P^{x_0}\big(\tau_{B(x_0,\phi^{-1}(\eps^{-1} C_6t_0))} \le t_0\big) \le \eps.
	\end{equation}
Note that  the map $(t,x) \mapsto P_tf(x):=\E^xf(X_t)$ is parabolic in $\sQ(t_1 + 1, x_0, \phi^{-1}(t_1))$ for all $t_1>0$. Hence, by Proposition \ref{p:phr}, for each $t_0 >T_0$, by taking $t_1>t_0$ large enough, one can get that for all non-negative $f \in L^\infty(M;\mu)$ and $x \in M$ with $d(x,x_0) \le \phi^{-1}(\eps^{-1}C_6 t_0)$, it holds
	\begin{equation}\label{e:law01_inf-2}
		|P_{t_1}f(x)- P_{t_1}f(x_0)| \le c_1  \sup_{t>0} \Vert P_t f \Vert_{L^\infty(M; \mu)} \Big( \frac{\phi^{-1}(\eps^{-1}C_6t_0) +  \phi^{-1}(t_1)^{\up/\up_1} }{ \phi^{-1}(t_1) } \Big)^{b} \le \eps \lVert f \Vert_{L^\infty(M; \mu)}.
	\end{equation}
	Indeed, we see that $(t_1, x), (t_1, x_0) \in \sI(t_1 + 1, x_0, \phi^{-1}(t_1))$ for all $t_1>0$ large enough.	Observe that \eqref{e:law01_inf-1} and \eqref{e:law01_inf-2} can play the same roles as \cite[(8.3) and (8.4)]{BB99} (or  \cite[(A.6) and (A.7)]{KKW17}). Hence, by following the proof of \cite[Theorem 8.4]{BB99} (see also the proof of \cite[Theorem 2.10]{KKW17}), we can deduce that $\P^{z}(G) \in \{0,1\}$ for all $z \in M$.

	Now, let $\sO= \{  z \in M :  \P^z (G) = 1  \}$  (which is Borel measurable) and suppose that both $\sO$ and $ \sO^c$ are  nonempty. Since at least one of $\mu(\sO)$ and $\mu(\sO^c)$ is  positive, by considering the event $G^c$ instead of $G$, we assume that $\mu(\sO) >0$ without loss of generality.
	Choose any $x_0 \notin \sO$.  Since $\mu(\sO)>0$, there exists $R > R_\infty \ep(x_0)^\up$ such that $\mu(B(x_0, \eta^2 R) \cap \sO)>0$.
	Then since $G$ is shift-invariant, we get from \NDLi \ and the Markov property that 
	\begin{align*}
		0 &= \P^{x_0}(G) = \P^{x_0}(G \circ \theta_{\phi(\eta R)}) \ge \P^{x_0} \big(\P^{X_{\phi(\eta R)}}(G) ; X_{\phi(\eta R)} \in B(x_0,\eta^2 R) \cap \sO  \big) \\
		&= \P^{x_0} \big(X_{\phi(\eta R)} \in B(x_0,\eta^2 R) \cap \sO\big) = \int_{ B(x_0,\eta^2 R) \cap \sO} p(\phi(\eta R), x_0, y) \mu(dy) \ge \frac{c_2\mu(B(x_0, \eta^2 R) \cap \sO)}{V(x_0, R)},
	\end{align*}
	which is a contradiction. Thus, we get that either $\sO = \emptyset$ or $\sO^c = \emptyset$. This completes the proof.
	\qed

	\section{Proofs of  Theorems \ref{t:limsup0-1}, \ref{t:limsupinf-1}, \ref{t:supprecise1} and \ref{t:supprecise2}}\label{s:limsup}

	We give a version of conditional Borel-Cantelli lemma which will be used in the proofs of our limsup LILs.
	
	\begin{prop}\label{p:law001}
		Let $(\tO, \tP, \GG, (\GG_t)_{t \ge 0})$ be a filtered probability space.

		\noindent (i) Suppose that there exist a decreasing sequence $(t_n)_{n \ge 1}$ such that $\lim_{n \to \infty} t_n = 0$ and a sequence of events $(A_n)_{n \ge 1}$ satisfying the following conditions:

		\setlength{\leftskip}{3mm}
		
		\noindent (Z1) $A_n \in \GG_{t_n}$ for all $n \ge 1$.
		
		\noindent (Z2) There exist a sequence of non-negative numbers $(a_n)_{n \ge 1}$, a sequence of events $(G_n)_{n \ge 1}$ and a constant $p_1 \in (0,1)$ such that $\sum_{n=1}^\infty a_n = \infty$ and for all $n \ge 1$,
		\begin{align*}
			\tP(G_n^c) \le p_1 \quad \text{and} \quad \tP(A_n \, | \,  \GG_{t_{n+1}}) \ge a_n \1_{G_n}, \quad \tP\text{-a.s.}
		\end{align*}
		
		\setlength{\leftskip}{0mm}
		\noindent Then, $\tP(\limsup A_n) \ge 1-p_1$. In particular, if $\lim_{n \to \infty} \tP(G_n^c)=0$, then $\tP(\limsup A_n) = 1$.

		\vspace{1mm}
		
		\noindent (ii) Suppose that there exist an increasing sequence $(s_n)_{n \ge 1}$ such that $\lim_{n \to \infty} s_n = \infty$ and	a sequence of events $(B_n)_{n \ge 1}$ satisfying the following conditions:
		
		\setlength{\leftskip}{3mm}
		
		\noindent (I1) $B_n \in \GG_{s_n}$ for all $n \ge 1$.
		
		\noindent (I2) There exist a sequence of non-negative numbers $(b_n)_{n \ge 1}$, a sequence of events $(H_n)_{n \ge 1}$ and a constant $p_2 \in (0,1)$ such that $\sum_{n=1}^\infty b_n = \infty$ and for all $n \ge 1$,
		\begin{align*}
			\tP(H_n^c) \le p_2 \quad \text{and} \quad \tP(B_{n+1} \, | \, \GG_{s_{n}}) \ge b_{n+1} \1_{H_{n+1}}, \quad \tP\text{-a.s.}
		\end{align*}
		
		\setlength{\leftskip}{0mm}
		\noindent	Then,  $\tP(\limsup B_n) \ge 1-p_2$. In particular, if $\lim_{n \to \infty} \tP(H_n^c)=0$, then $\tP(\limsup B_n) = 1$.	
	\end{prop}
	\pf Since the proofs are similar, we only give the proof for (i). 
	Choose any $n$ and $\eps \in (0, 1-p_1)$. It suffices to show that
	$\tP (\cup_{k \ge 0} A_{n+k}) \ge 1-p_1-\eps$.
	To prove this, we fix any $N$ such that $\sum_{k=0}^{N}a_{n+k} > \eps^{-1}$ and observe that
	\begin{equation}\label{e:setdecomp}
		\tP (\cup_{k \ge 0} A_{n+k}) \ge \sum_{k=0}^{N} \tP\big(A_{n+k} \cap (\cup_{i=1}^{N-k} A_{n+k+i})^c\big).
	\end{equation}
	Let $L_{n,k}:= \cup_{i=1}^{N-k} A_{n+k+i} = \cup_{i=k+1}^{N}A_{n+i}$. By (Z1) and (Z2), it holds that for all $0 \le k \le N$,
	\begin{align}\label{e:setcal}
		\tP(A_{n+k} \cap L_{n,k}^c) &= \tE \big[\,\tE [\1_{A_{n+k}} \1_{L_{n,k}^c} \, |\, \GG_{t_{n+k+1}} ] \, \big]= \tE \big[\, \1_{L_{n,k}^c}\tE [\1_{A_{n+k}}  \, |  \, \GG_{t_{n+k+1}} ] \, \big] \nn\\
		& \ge a_{n+k} \tE[ \1_{L_{n,k}^c}1_{G_{n+k}}] \ge a_{n+k} ( \tP(L_{n,k}^c)  - p_1).
	\end{align}
	In the last inequality, we used the fact that $\tP(A \cap B) \ge \tP(A) - \tP(B^c)$.
	
	Suppose that $\tP(L_{n,0}^c) > p_1 + \eps$. Then $\tP(L_{n,k}^c) \ge \tP(L_{n,0}^c) >p_1 + \eps$ for all $0 \le k \le N$. Thus, by \eqref{e:setdecomp} and \eqref{e:setcal},  $\tP (\cup_{k \ge 0} A_{n+k}) \ge \sum_{k=0}^{N} a_{n+k}( \tP(L_{n,k}^c)  - p_1) \ge \eps\sum_{k=0}^{N} a_{n+k} >1$, 
	which is a contradiction. Therefore, we conclude that $\tP(L_{n,0}^c) \le p_1+\eps$ and hence $\tP(\cup_{k \ge 0} A_{n+k}) \ge \tP(L_{n,0}) \ge 1-p_1-\eps$. 
	\qed
	
	We know that a mixed stable-like process and a Brownian-like jump process have totally different limsup LILs. Indeed, the former enjoys a limsup LIL of type \eqref{e:ieq2}, while the latter enjoys a one of type \eqref{e:ieq1}. (See \cite{KKW17, BKKL19a}.)
	Under our setting, the process $X$ may behave like a mixed stable-like process in some time range, while it may behave like  a  Brownian-like jump process in  some other time range. See \cite[Section 4.2]{CK202} for an explicit example of a subordinator which behaves like this.  Hence,
	one must overcome significant technical difficulties to obtain limsup LILs for  our $X$.  The following proofs of Theorems \ref{t:limsup0-1}(ii) and \ref{t:limsupinf-1}(ii), together with Propositions \ref{p:EP} and \ref{p:EPL}, are the most delicate part of this paper.
	
	\vspace{2mm}

	\noindent \textbf{Proof of Theorem \ref{t:limsup0-1}.} Choose any $x \in U$ and let $\eps_x:=2^{-1}(\updelta_U(x) \land 1) \in (0, 1)$. 
	
	\smallskip

 (i) By \eqref{e:phipsi} and \eqref{limsup0-case1}, we have $\phi(r) \asymp \psi(r)$ for $r \in (0,1]$. Then by \eqref{e:EP_0} and \eqref{e:EPdown0}, there exists  $c_1 \ge 1$ such that for all $z \in U$,  $0<r<3^{-1}\big(R_0 \land (C_1 \updelta_U(z))\big)$ and $0<t \le \phi(r)$, 
	\begin{equation}\label{e:Main1}
		\frac{c_1^{-1}t}{\phi(r)} \le \P^z(\tau_{B(z,r)} \le t) \le \frac{c_1t}{\phi(r)}.
	\end{equation}
	
 (i-a) First, assume that $\int_0 \phi(\Psi(t))^{-1} dt  = \infty$. Choose any $K\ge1$ and  define
	\begin{align*}
		A_n:=\Big\{ \sup_{s \in [2^{-(n+1)},2^{-n}]}d(X_0, X_s) \ge K\big(\Psi(2^{-n})\vee\phi^{-1}(2^{-n})\big)\Big\} \quad \text{and} \quad  G_n:=\big\{\tau_{B(X_0, \eps_x)} \ge 2^{-n}\big\}.
	\end{align*}
	Observe that by the Markov property and the triangle inequality,  we have
	\begin{align*}
		\P^x(A_n \, |  \, \FF_{2^{-(n+1)}}) &\ge \inf_{y \in B(x,\eps_x)}\P^y\Big( \sup_{s \in [0,2^{-(n+1)}]}d(X_0, X_s) \ge 2K\big(\Psi(2^{-n})\vee\phi^{-1}(2^{-n})\big)\Big)  \1_{G_n} =: a_n \1_{G_n}.
	\end{align*}
	According to \eqref{e:Main1}, since $\U_{R_0}(\phi, \beta_2, C_U)$ holds,  we see that for all $n$ large enough,
	\begin{equation}\label{e:supjumplow}
		a_n \ge \frac{c_1^{-1}2^{-(n+1)}}{\phi\big(2K(\Psi(2^{-n})\vee\phi^{-1}(2^{-n}))\big)} = \frac{(2c_1)^{-1}2^{-n}}{\phi(2K\Psi(2^{-n})) \vee \phi(2K\phi^{-1}(2^{-n}))}  \ge c_2K^{-\beta_2} \Big(1 \land \frac{2^{-n}}{\phi(\Psi(2^{-n}))}\Big).
	\end{equation}
	Since $\int_0 \phi(\Psi(t))^{-1} dt  = \infty$, it follows  $\sum_{n=1}^{\infty}a_n=\infty$. Moreover, we also get from  \eqref{e:Main1} that  
	\begin{equation*}
		\lim_{n \to \infty} \P^x(G_n^c) = \lim_{n \to \infty} \P^x(\tau_{B(x, \eps_x)} < 2^{-n}) \le c_1\lim_{n\to \infty}2^{-n} = 0.
	\end{equation*} Thus, by Proposition \ref{p:law001}(i), we obtain that
	\begin{equation*}
		\limsup_{t \to 0}d(x,X_t)/\Psi(t) \ge \lim_{n \to \infty}\sup_{s \in [2^{-(n+1)}, 2^{-n}]}d(x,X_s)/\Psi(2^{-n}) \ge K, \quad \P^x\text{-a.s.}
	\end{equation*}
	Since the above holds with any $K \ge 1$, we conclude that  $\limsup_{t \to 0}d(x,X_t)/\Psi(t) = \infty$, $\P^x$-a.s.
	
	 (i-b) Next, assume that $\int_0\phi(\Psi(t))^{-1}dt<\infty$. Choose any $\lambda \in (0, 1)$ and define
	$$
	B_n:=\big\{ \sup_{s \in [0,2^{-n+1}]}d(X_0, X_s) \ge \lambda \Psi(2^{-n})\big\}.$$
	By \eqref{e:Main1}, since $\U_{R_0}(\phi, \beta_2, C_U)$ holds,  for all $n$ large enough, we have $\P^x(B_n) \le  c_12^{-n}/\phi(\lambda \Psi(2^{-n})) \le c_3 2^{-n}/\phi(\Psi(2^{-n}))$. Since $\int_0 \phi(\Psi(t))^{-1}dt<\infty$, it follows $\sum_{n=1}^\infty \P^x(B_n)<\infty$. Hence, by the Borel-Cantelli lemma, we get $\limsup_{t \to 0}\sup_{0<s \le t}d(x,X_s)/\Psi(t) \le \lim_{n \to 0} \sup_{0<s \le 2^{-n+1}}d(x,X_s)/\Psi(2^{-n})  \le  \lambda$, $\P^x$-a.s. Since we can choose $\lambda$ arbitrarily small, this finishes the proof for (i).
	
	\medskip
	
	(ii)  Fix any $\delta>0$ and choose a decreasing sequence $(s_n)_{n \ge 4} \subset (0, 1)$ such that 
	\begin{equation}\label{e:defsn}
		n\phi(s_n) (\log n)^{2+\delta}  \le \psi(s_n \log n) \quad \text{ and } \quad 4\phi(s_{n+1}) \le \phi(s_n) \quad \;\; \text{for all} \; n\ge 4.
	\end{equation}
	Such sequence $(s_n)_{n \ge 4}$  exists because  $\limsup_{r \to 0}\psi(r)/\phi(r)=\infty$. 
	Then we define  $\Psi\in \sM_+$ as 
	\begin{equation}\label{e:defPsi}
		\Psi(t):= \sum_{n=4}^\infty (s_n \log n)  \1_{(t_{n+1},t_n]}(t) + (s_4 \log 4) \cdot  \1_{(t_4, \infty)}(t) \quad \text{where} \;\; t_n:=\phi(s_n)\log n.
	\end{equation}
	We claim that there exist constants $0<q_1 \le q_2<\infty$ such that for all $x \in U$,
	\begin{equation}\label{e:sup0}
		q_1 \le \limsup_{t \to0} \frac{d(x, X_t)}{\Psi(t)} \le \limsup_{t \to0} \frac{\sup_{0<s \le t}d(x, X_s)}{\Psi(t)} \le q_2,~\qquad\,
		\P^x\mbox{-a.s.}
	\end{equation} 
	The second inequality above is obvious. Recall that by Proposition \ref{p:E}(i), under the present setting, \as \ hold with $R_0=r_0$ where the constant $r_0$ is  defined as \eqref{e:defr0}.
	
	First, we prove the third inequality in \eqref{e:sup0}.  Set $q_2:=2C_8^{-1}+C_9^{-1}+1$ where $C_8, C_9$ are constants in Proposition \ref{p:EP}(i). According to Proposition \ref{p:EP}(i),  for all $n$ large enough,
	\begin{align}\label{e:limsupup}
		&\P^x\big(d(x,X_t)>q_2\Psi(t) \;\; \text{for some} \;\; t \in (t_{n+1}, t_{n}] \big) \nn\\
		&\le \P^x \big(\tau_{B(x, q_2 s_n \log n)} \le t_n \big) \le C_7 \bigg(\frac{t_n}{\psi(q_2s_n \log n)} +  \exp\Big(-\frac{q_2C_8  s_n \log n}{\vt_1(t_n,q_2C_9  s_n \log n)}\Big)\bigg).
	\end{align}
	Since $\psi$ is increasing, by \eqref{e:defsn}, we have that $\sum_{n=4}^\infty t_n/\psi(q_2 s_n \log n) \le \sum_{n=4}^\infty n^{-1}(\log n)^{-1-\delta}<\infty$. Besides, since $\psi$ and $\phi$ are increasing, and $\U_{R_0}(\psi, \beta_3, C_U')$ holds, by \eqref{e:defsn}, we also have that 
	\begin{equation}\label{e:checkrange}
		\psi^{-1}(t_n) \le \psi^{-1}(n^{-1}\psi(s_n\log n))  \le  s_n <\phi^{-1}(t_n/2) \quad \text{for all $n$ large enough.}
	\end{equation}
	Hence, for all  $n$ large enough, since $q_2C_9 > 1$ and $t_n s_n/\phi(s_n)= s_n \log n \in [\Up_1(t_n), \Up_2(t_n)]$, we have $\vt_1(t_n, q_2C_9s_n \log n) \le \vt_1(t_n, s_n \log n) \le s_n$. It follows that
	\begin{equation*}
		\sum_{n=4}^\infty \exp\Big(-\frac{q_2C_8s_n \log n}{\vt_1(t_n, q_2C_9s_n \log n)}\Big) \le c_1 + c_2 \sum_{n=4}^\infty n^{-q_2C_8} \le  c_1 + c_2 \sum_{n=4}^\infty n^{-2}<\infty.
	\end{equation*}
	Therefore, by the Borel-Cantelli lemma,  we can see that the upper bound in \eqref{e:sup0} holds.

	Next, we prove the first inequality in \eqref{e:sup0}.
	Choose $C_{11},C_{12}$ according to  Proposition \ref{p:EPL}(i) with $a=C_1'$ and set $q_1:=4^{-1}(C_{11}+C_{12})^{-1}$ where $C_1'$ is the constant in Proposition \ref{p:E}. Define
	\begin{equation*}
		E_n:=\Big\{ \sup_{s \in [t_{n+1},t_n]}d(X_0, X_s) \ge q_1 \Psi(t_n) \Big\} \quad \text{and}\quad F_n:=\big\{\tau_{B(X_0, \eps_x)} \ge t_n\big\}.
	\end{equation*}
By the Markov property and the triangle inequality, since $t_n > 2t_{n+1}$ for all $n \ge 4$ due to the second inequality in  \eqref{e:defsn}, we have that
	\begin{equation*}
		\P^x(E_n \, | \, \FF_{t_{n+1}}) \ge \inf_{y \in B(x,\eps_x)} \P^y\big( \tau_{B(y, 2q_1\Psi(t_n))} \le 2^{-1}t_n \big) \cdot \1_{F_n} =:b_n \1_{F_n} \quad \text{for all} \;\; n \ge 4.
	\end{equation*}
Observe that for all $n$ large enough, $s_n \in [\psi^{-1}(2^{-1}t_n), \phi^{-1}(2^{-1}t_n)]$ by  \eqref{e:checkrange}. Besides, for all $n \ge 4$, since $2^{-1}t_ns_n/\phi(s_n) = 2^{-1}s_n \log n \ge 2q_1C_{12} s_n \log n$,  we have 
	$$\vt_2(2^{-1}t_n, 2q_1 C_{12} s_n \log n) \ge \vt_2(2^{-1}t_n, 2^{-1}t_ns_n/\phi(s_n)) \ge s_n.$$ Thus,  using   Proposition \ref{p:E}(i) and \ref{p:EPL}(i), we get that for all $n$ large enough, whether $2^{-1}t_n \le \phi(2q_1C_1 \Psi(t_n))$ or not,
	\begin{align}\label{e:checklow}
		b_n &\ge (1-e^{-a_2}) \wedge \big(C_{10}  \exp (- 2q_1C_{11}\Psi(t_n)/s_n )\big) \ge c_3 n^{-2q_1 C_{11}} \ge c_3 n^{-1/2}.
	\end{align}
	for large $n$ satisfying $2q_1\Psi(t_n) < 2^{-1}c'\big(C_1^{-1}R_0 \land \big( (C_V \land C_1 \land C_2) \updelta_U(x) \big) \big)$.
	Thus, $\sum_{n=4}^\infty b_n =\infty$.  Note that  $\lim_{n \to \infty}\P^x(F_n^c)=0$ by Proposition \ref{p:EP}(i).  Therefore, by Proposition \ref{p:law001}(i) and the monotone property of $\Psi$,  we  conclude that the first inequality in \eqref{e:sup0} holds.

	In the end, we deduce \eqref{e:limsuplaw1} from \eqref{e:sup0} and the Blumenthal's zero-one law. \qed

	\noindent \textbf{Proof of Theorem \ref{t:limsupinf-1}.} 
	Choose any $x,y \in M$ and put $\up_1:=\up^{2/3} \in (\up,1)$. 
	
	\smallskip
	
	 (i)  Since  $\phi(r) \asymp \psi(r)$ for $r \ge 1$ by \eqref{e:phipsi} and \eqref{limsupinf-case1}, according to \eqref{e:EP_0} and \eqref{e:EPdown0}, there exist $c_1,c_2 \ge 1$ such that for all $z \in M$ and  $r>c_1\ep(z)^{\up_1}$, 
	\begin{equation}\label{e:supEP}
		\P^z(\tau_{B(z,r)} \le t) \ge \frac{c_2^{-1}t}{\phi(r)} \quad \text{for} \;\; 0<t \le \phi(r)
	\end{equation}
and
	\begin{equation}\label{e:supEP2}
 \P^z(\tau_{B(z,r)} \le t) \le \frac{c_2t}{\phi(r)}\quad \text{for} \;\; t \ge \phi(2r^{\up/\up_1}).
\end{equation}
	
	 (i-a) First, assume that $\int^\infty \phi(\Psi(t))^{-1}dt=\infty$. Choose any $K \ge 1$ and define 
	\begin{align*}
		A_n':=\Big\{ \sup_{s \in [2^{n-1},2^{n}]}d(X_0, X_s) \ge K\big(\Psi(2^n)\vee\phi^{-1}(2^n)\big) \Big\} \;\; \text{and} \;\; G_n':=\big\{\tau_{B(X_0,\phi^{-1}(2^{n})^{\up_1/\up} )} \ge C_U2^{\beta_2+n}\big\}.
	\end{align*}
	Then by the Markov property and the triangle inequality, we get
	\begin{align*}
		\P^y(A_{n}' \, |  \, \FF_{2^{n-1}}) &\ge \inf_{z \in B(y,\phi^{-1}(2^{n})^{\up_1/\up} )}\P^z\Big( \sup_{s \in [0,2^{n-1}]}d(X_0, X_s) \ge 2K\big(\Psi(2^{n})\vee\phi^{-1}(2^{n})\big)\Big)  \1_{ G_n'} =: a_n' \1_{G_n'}.
	\end{align*}
	Observe that for all $n$ large enough and $z \in B(y, \phi^{-1}(2^{n})^{\up_1/\up} )$, since $\up_1^2/\up = \up^{1/3}<1$, we have 
	$$c_1 \ep(z)^{\up_1} \le c_1( \ep(y)^{\up_1} +  \phi^{-1}(2^{n})^{\up_1^2/\up} ) < \phi^{-1}(2^n) < 2K\big(\Psi(2^{n})\vee\phi^{-1}(2^{n})\big).$$
	Hence, by \eqref{e:supEP} and $\U^{R_\infty}(\phi, \beta_2, C_U)$, we get that for all $n$ large enough (cf. \eqref{e:supjumplow}),
	$$ a_n' \ge \frac{c_2^{-1}2^{n-1}}{\phi\big(2K(\Psi(2^n) \vee \phi^{-1}(2^n))\big)} \ge c_3K^{-\beta_2} \Big(1 \land \frac{2^{n}}{\phi(\Psi(2^{n}))}\Big)$$
	so that $\sum_{n=1}^\infty a_n'=\infty$ since $\int^\infty \phi(\Psi(t))^{-1}dt = \infty$. We also see from \eqref{e:supEP},  $\U^{R_\infty}(\phi, \beta_2, C_U)$ and  $\L^{R_\infty}(\phi, \beta_1, C_L)$ that
	$$\lim_{n \to \infty} \P^y( G_n'^{\,c}) \le \lim_{n \to \infty} \frac{c_2C_U2^{\beta_2} \phi\big(\phi^{-1}(2^n)\big)}{\phi\big(\phi^{-1}(2^n)^{\up_1/\up}\big)} \le  c_2C_U2^{\beta_2} C_L^{-1} \lim_{n \to \infty} \phi^{-1}(2^n)^{-\beta_1(\up_1-\up)/\up}=0.$$  
	Therefore, according to Proposition \ref{p:law001}(ii), we obtain that
	\begin{equation*}
		\limsup_{t \to \infty} \frac{d(x,X_t)}{\Psi(t)} \ge \limsup_{n \to \infty} \frac{ \sup_{s \in [2^{n-1}, 2^{n}]}d(y, X_s)}{\Psi(2^n)} - \limsup_{t \to \infty}\frac{d(x,y)}{\Psi(t)} \ge K- \limsup_{t \to \infty}\frac{d(x,y)}{\Psi(t)}, \;\; \P^y\text{-a.s.}
	\end{equation*}
	Since we can choose $K$ arbitrarily large, this yields that $\limsup_{t \to \infty} d(x,X_t)/\Psi(t)=\infty$, $\P^y$-a.s.
	
 (i-b) Now, assume that $\int^\infty \phi(\Psi(t))^{-1}dt<\infty$. Define 
	$\Psi_\land(t):= \Psi(t) \land \phi^{-1}(t)^{\up_1/\up}$. Then we have that, by $\L^{R_\infty}(\phi, \beta_1, C_L)$ and $\U^{R_\infty}(\phi, \beta_2, C_U)$,
	\begin{equation}\label{e:intPsi}
		\int^\infty \frac{dt}{\phi(\Psi_\land(t))} \le \int^\infty \frac{dt}{\phi(\Psi(t))} + \int^\infty \frac{t^{-1}\phi(\phi^{-1}(t))}{\phi(\phi^{-1}(t)^{\up_1/\up})}dt \le c_4\bigg(1+ \int^\infty \frac{dt}{t^{1+\beta_1\beta_2^{-1}(\up_1-\up)/\up}} \bigg)  <\infty.
	\end{equation}
	For $\lambda \in (0,1/2)$, let $B_n':=\{ \sup_{s \in [0, 2^{n+1}]} d(X_0, X_s) \ge \lambda \Psi_\land(2^n) \}$. Note that  $\phi(2\lambda \Psi_\land(2^n)^{\up/\up_1}) \le (\phi \circ \phi^{-1})(2^n) <2^{n+1}$. Hence, by \eqref{e:supEP2},  $\P^y(B_n') \le c_22^{n+1}/\phi(\lambda \Psi_\land(2^n))$ for all $n$ large enough so that $\sum_{n=1}^\infty \P^y(B_n')<\infty$ by $\U^{R_\infty}(\phi, \beta_2, C_U)$ and \eqref{e:intPsi}. 
	Therefore, by the Borel-Cantelli lemma, since $\lim_{t \to \infty}\Psi_\land(t)=\infty$ in view  of \eqref{e:intPsi}, we conclude that  
	$$\limsup_{t \to \infty}\frac{\sup_{0<s \le t}d(x,X_s)}{\Psi(t)} \le \limsup_{n \to \infty}\frac{d(x,y) + \sup_{0<s \le 2^{n+1}}d(y,X_s)}{\Psi_\land(2^n)}\le \lambda, \quad \P^y\text{-a.s.}$$
	Since the above inequalities hold with any $\lambda \in (0,1/2)$,  we finish the proof for (i).
	\medskip
	
	(ii) Define $\psi_\land(r):=\psi(r) \land \phi(1+r^{\sqrt{\up_1/\up}}) \in [\phi(r), \psi(r)]$. Then  $\limsup_{r \to \infty} \psi_\land(r)/\phi(r)=\infty$ since $\limsup_{r \to \infty}\psi(r)/\phi(r) = \lim_{r \to \infty} \phi(r^{\sqrt{\up_1/\up}})/\phi(r)= \infty$ due to $\L^{R_\infty}(\phi, \beta_1, C_L)$. Fix any $\delta>0$ and  find an increasing sequence $(s_n)_{n \ge 4}$ such that $s_n \ge e^n$, $n\phi(s_n)(\log n)^{2+\delta} \le \psi_\land(s_n \log n)$  and $\phi(s_{n+1}) \ge 4\phi(s_n)$ for all $n$.  We  put 
	\begin{equation*}
		\Psi(t):= \sum_{n=5}^\infty (s_n \log n)  \1_{(t_{n-1},t_n]}(t) \quad \text{where} \;\; t_n:=\phi(s_n)\log n.
	\end{equation*}

	Let $C_8, C_9$ be the constants in Proposition \ref{p:EP}(ii) and choose $C_{11},C_{12}$ according to Proposition \ref{p:EPL}(ii) with $a=C_1'$ where $C_1'$ is the constant in Proposition \ref{p:E}. Then  we  set   $q_2:=2C_8^{-1}+C_9^{-1}+1$ and $q_1:=4^{-1}(C_{11}+C_{12})^{-1}$, and we fix $x, y \in M$.
	
	By Proposition \ref{p:EP}(ii),  for all $n$ large enough, using  $B(x, 2q_2 s_n\log n)\supset B(y, q_2s_n \log n)$, we get
	\begin{align}\label{e:limsupinf_1}
		&\P^y\big(d(x,X_t)>2q_2\Psi(t) \;\; \text{for some} \;\; t \in (t_{n-1}, t_{n}] \big)  \le \P^y\big(\tau_{B(x, 2q_2 s_n \log n )} \le t_n \big)  \nn\\
		&\le \P^y \big(\tau_{B(y, q_2 s_n \log n )} \le t_n \big) \le C_7 \bigg(\frac{t_n}{\psi_\land (q_2s_n \log n)} +  \exp\Big(-\frac{q_2C_8  s_n \log n}{\vt_1(t_n,q_2C_9  s_n \log n)}\Big)\bigg).
	\end{align}
	Indeed, for all $n$ large enough, since $\up<\up_1$, $s_n \ge e^n$ and $\phi$ is increasing,  we see  that
	\begin{equation}\label{e:limsupinf_2}
		\psi_\land\big(2 (q_2 s_n \log n)^{\up/\up_1}\big) \le \phi\big(2^{1+\sqrt{\up_1/\up}}(q_2s_n\log n)^{\sqrt{\up/\up_1}}\big) \le \phi(s_n) < t_n.
	\end{equation}
	Using \eqref{e:limsupinf_1}, by following the proof of Theorem \ref{t:limsup0-1}(ii), we obtain 
	\begin{equation}\label{e:limsupinf_3}
		\limsup_{t \to \infty} \frac{\sup_{0<s\le t} d(x,X_s)}{\Psi(t)} \le 2q_2, \quad \P^y\text{-a.s.}
	\end{equation}

	On the other hand, define
	\begin{equation*}
		E_n':=\Big\{ \sup_{s \in [t_{n-1},t_n]}d(X_0, X_s) \ge q_1 \Psi(t_n) \Big\} \quad \text{and}\quad  F_n':=\big\{\tau_{B(X_0, s_n^{\,\up_1/\up})} \ge C_U2^{\beta_2}t_{n}\big\}.
	\end{equation*}
	Then by the Markov property and the triangle inequality, since $t_n > 2t_{n-1}$ for all $n$, we have
	\begin{equation*}
		\P^y(E_n' \, | \, \FF_{t_{n-1}}) \ge \inf_{z \in B(y,s_n^{\,\up_1/\up})} \P^z\big( \tau_{B(y, 2q_1\Psi(t_n))} \le 2^{-1}t_n \big) \cdot \1_{F_n'} =:b_n' \1_{F_n'} \quad \text{for all} \;\; n \ge 5.
	\end{equation*}
	For any constant $R>0$, by $\L^{R_\infty}(\phi, \beta_1, C_L)$, $\U^{R_\infty}(\phi, \beta_2, C_U)$ and similar calculations to \eqref{e:limsupinf_2}, we have that, for all $n$ large enough and $z \in B(y, s_n^{\,\up_1/\up})$, since $\up_1^2/\up<1$,
	\begin{equation*}
		R\ep(z)^{\up_1}<R\ep(y)^{\up_1} + R s_n^{\up_1^2/\up} < 2q_1 s_n \log n =  2q_1\Psi(t_n) \quad \text{and} \quad \psi_\land\big( R (2q_1 \Psi(t_n))^{\up/\up_1}  \big)<2^{-1}t_n.
	\end{equation*}
	Hence, using \eqref{e:EPdown}, by following the proof of Theorem \ref{t:limsup0-1}(ii), we get $\sum_{n=5}^\infty \P^y(E_n')=\infty$. Moreover, we also get from \eqref{e:EP_0} and $\L^{R_\infty}(\phi, \beta_1, C_L)$ that 
	$$\lim_{n \to \infty} \P^y(F_n') \le c_1 \lim_{n \to \infty} \phi(s_n)\log s_n/ \phi(s_n^{\up_1/\up})=0.$$
	It follows from Proposition \ref{p:law001}(ii) that $\P^y(\limsup E_n')=1$ so that
	\begin{equation}\label{e:limsupinf_4}
		\limsup_{t \to \infty} \frac{d(x, X_t)}{\Psi(t)} \ge \limsup_{t \to \infty} \frac{d(y,X_t)-d(x,y)}{\Psi(t)} = \limsup_{t \to \infty} \frac{d(y, X_t)}{\Psi(t)}\ge q_1, \quad \P^y\text{-a.s.}
	\end{equation}
	
	Since $\lim_{t \to \infty} \Psi(t)=\infty$,  one can verify that  events 
	$\{\limsup_{t \to \infty} \sup_{0 < s \le t}d(x,X_s)/\Psi(t) \ge \lambda  \}$ and $\{\limsup_{t \to \infty} d(x,X_t)/\Psi(t) \ge \lambda  \}$ are shift-invariant for any $\lambda>0$.  Eventually, by \eqref{e:limsupinf_3}, \eqref{e:limsupinf_4} and Proposition \ref{p:law01}, we conclude the result. \qed

	\noindent \textbf{Proofs of Theorems \ref{t:supprecise1} and \ref{t:supprecise2}.} Since the proofs are similar, we only prove Theorem \ref{t:supprecise1}. 
	
	Set $p:=2(4/C_L)^{1/\beta_1}$ and  $q_2:=2C_8^{-1}+(C_9C_L)^{-1}+1$.  Since $\L_{R_0}(\phi, \beta_1, C_L)$ holds, we see that  $\phi(\lambda) \ge 4 \phi(\lambda/p)$ for all $\lambda \in (0, R_0)$. Define $s_n:=p^{-n}$, $t_n:=\phi(s_n) \log n$ and 
	$$\Psi_0(t):= \sum_{n=4}^\infty (s_n \log n) \cdot 1_{(t_{n+1}, t_n]}(t).$$ Then since $\L_{R_0}(\phi, \beta_1, C_L)$ and $\U_{R_0}(\phi, \beta_2, C_U)$ hold, one can see that for all  small enough $t$,
	\begin{equation}\label{e:stepcompare}
		2^{-1}C_U^{-1}\phi^{-1}(t/\log|\log t|) \log |\log t| \le  \Psi_0(t) \le 2C_L^{-1}p\phi^{-1}(t/\log|\log t|) \log |\log t|.
	\end{equation}
	
	 (i) We claim that  there exist constants $0<q_1 \le q_2<\infty$ such that for all $x \in U$, \eqref{e:sup0} is valid with $\Psi_0$ instead of $\Psi$. Indeed, we can see that \eqref{e:limsupup} still works (with redefined $s_n$ and $t_n$), and since $\phi$ and $\psi$ are increasing and $\U_{R_0}(\phi, \beta_2, C_U)$ holds, it still holds that
	\begin{equation}\label{e:jumppart}
		\sum_{n=4}^\infty \frac{t_n}{\psi(q_2 s_n \log n)} \le c_1\sum_{n=4}^\infty \frac{\phi(p^{-n}) \log |\log p^{-n}|}{\psi(p^{-n} \log |\log p^n|)} \le c_2 \int_0^{p^{-4}} \frac{\phi(r)\log |\log r|}{r \psi(r \log |\log r|)}<\infty.
	\end{equation}
	
	Note that $s_n \le \phi^{-1}( t_n)$ for all $n \ge 4$. Thus, for all $n$ large enough,  if $s_n \ge \psi^{-1}( t_n)$, then   $\vt_1( t_n, q_2C_9 s_n \log n) \le \vt_1( t_n, s_n \log n) = \vt_1(\phi(s_n)\log n, s_n \log n)  \le   s_n$.  Otherwise, if $ s_n <\psi^{-1}(t_n)$, then since $\L_{R_0}(\phi, \beta_1, C_L)$ holds with $\beta_1 \ge 1$, we get
	$$
	q_2C_9 s_n \log n =q_2 C_9  t_n\frac{s_n}{\phi(s_n)} \ge q_2C_9 C_L t_n \frac{ \psi^{-1}( t_n)}{\phi(\psi^{-1}( t_n))}>t_n \frac{ \psi^{-1}( t_n)}{\phi(\psi^{-1}( t_n))}
	$$
	and hence  $\vt_1( t_n, q_2C_9 s_n \log n) \le \vt_1( t_n,  t_n \psi^{-1}( t_n)/\phi(\psi^{-1}( t_n)))  =  \psi^{-1}( t_n)$. Thus, for all $n$ large enough, whether $ s_n \ge \psi^{-1}( t_n)$ or not, we have that  
	$$\vt_1( t_n, q_2C_9 s_n \log n) \le  s_n \vee \psi^{-1}( t_n).$$ 
	Therefore, since $\U_{R_0}(\psi, \beta_3, C_U')$ holds, by \eqref{e:jumppart}, we get that
	\begin{align}\label{e:Psi0}
		&\sum_{n=4}^\infty \exp\Big(-\frac{q_2C_8 s_n \log n}{\vt_1(t_n, q_2C_9 s_n \log n)}\Big) \le c_3 \sum_{n=4}^\infty \Big(n^{-q_2C_8} + \exp\Big(-\frac{q_2C_8 s_n \log n}{\psi^{-1}( t_n)}\Big) \Big) \nonumber\\
		&\le  c_4 + c_5\sum_{n=4}^\infty \Big(\frac{\psi^{-1}( t_n)}{ s_n \log n} \Big)^{\beta_3} \le c_4 + c_6 \sum_{n=4}^\infty \frac{t_n}{\psi( s_n \log n)}<\infty.
	\end{align}
	Hence, by the Borel-Cantelli lemma, we can deduce that the third inequality in \eqref{e:sup0} still holds. 
	
	On the other hand, by using the fact that $\vt_2(2^{-1}t_n, r) \ge \psi^{-1}(2^{-1}t_n)$ for all $r$, we can follow the proof for the first inequality in \eqref{e:sup0} line by line, whether $s_n \ge \psi^{-1}(2^{-1}t_n)$ or not. Here, we mention that the condition \Tailog \ is unnecessary in that proof.
	
	Finally, in view of  \eqref{e:sup0} (with $\Psi_0$) and \eqref{e:stepcompare}, by using the Blumenthal's zero-one law again,  we conclude the desired result.
	
	\smallskip
	
 (ii) Since $\U_{R_0}(\phi, \beta_2, C_U)$ and  $\U_{R_0}(\psi, \beta_3, C_U')$ hold, for every $K \ge 1$, we have
	\begin{equation*}
		\sum_{n=4}^\infty \frac{t_n}{\psi(Ks_n \log n)} \ge c_1 \sum_{n=4}^\infty \frac{\phi( s_n) \log |\log s_n|}{\psi(s_n \log |\log  s_n|)} \ge c_2 \int_0^{1/16}\frac{\phi(r)\log|\log r|}{r\psi(r\log |\log r|)} dr = \infty.
	\end{equation*}
	Then by a similar proof to the one for Theorem \ref{t:limsup0-1}(i) (the case when $\int_0 dt/\phi(\Psi(t)) =\infty$), using \eqref{e:EPdown0},  one can deduce that for all $x\in U$, it holds $\limsup_{t \to 0} d(x,X_t)/ \Psi_0(t) =\infty$,  $\P^x$-a.s. In view of \eqref{e:stepcompare}, this shows that for all $x \in U$,  the limsup in \eqref{e:1.14} is infinte, $\P^x$-a.s.

	Now, we prove \eqref{e:suppre1}.  Let $\Psi \in \sM_+$ be a function satisfying \eqref{e:limsuplaw1}. (According to Theorem \ref{t:limsup0-1}(ii), at least one such $\Psi$ exists.) Then we define $f_0(t):=\Psi(t)/(\phi^{-1}(t/\log |\log t|)\log|\log t|)$. From the definition, one can see that $\limsup_{t \to 0}f_0(t)=\infty$  since the limsup in \eqref{e:1.14} is infinite,  $\P^x$-a.s for all $x \in U$. Thus, it remains to show that $\liminf_{t \to 0}f_0(t)<\infty$.
	
	Suppose  $\liminf_{t \to 0}f_0(t)=\infty$. Then by \eqref{e:stepcompare}, we have $\liminf_{t \to 0}\Psi(t)/\Psi_0(t)=\infty$. Fix any $x \in U$ and let $E_n(\eps):=\P^x(d(x,X_t)>\eps \Psi(t) \; \text{for some} \; t \in (t_{n+1}, t_n])$  for $\eps>0$. According to Proposition \ref{p:EP}(i), it holds  that  for all $n$ large enough,
	\begin{equation*}
		\P^x(E_n(\eps))\le c_3 \Big(\frac{t_{n}}{\psi( \eps \Psi(t_{n+1}))} +  \exp\Big(-\frac{\eps C_8  \Psi(t_{n+1})}{\vt_1(t_n,\eps C_9  \Psi(t_{n+1}))}\Big)\Big)=:c_3(I_{n,1}+I_{n,2}).
	\end{equation*}
	
	Observe that $\int_0 dt/\psi(\Psi(t))<\infty$. Indeed, if this integral is infinite, then by repeating  the proof for Theorem \ref{t:limsup0-1}(i) again, one can see that for all $x \in U$, it holds  $\limsup_{t \to 0} d(x,X_t)/ \Psi(t) =\infty$,  $\P^x$-a.s., which is contradictory to \eqref{e:limsuplaw1}.  Since $\U_{R_0}(\psi, \beta_3, C_U')$ holds, it follows that  $\int_0 dt/\psi(\eps \Psi(t))<\infty$ and hence $\sum_{n=1}^\infty I_{n,1}<\infty$. 
	
	On the other hand, note that since $\L_{R_0}(\phi, \beta_1, C_L)$ holds, in view of \eqref{e:stepcompare}, one can check that  there exists $c_4>1$ such that $c_4\Psi_0(t_{n+1}) \ge \Psi_0(t_n)$ for all $n \ge 4$. Since  $\vt_1(t, \cdot)$ is non-increasing and $\liminf_{t \to 0} \Psi(t)/\Psi_0(t)=\infty$,  there exists $N(\eps)$ depending on $\eps$ such that 
	\begin{equation*}
		I_{n,2}\le \exp\Big(-\frac{c_4 q_2C_8\Psi_0(t_{n+1})}{\vt_1(t_n, c_4 q_2C_9\Psi_0(t_{n+1}))}\Big)  \le \exp\Big(-\frac{q_2C_8\Psi_0(t_{n})}{\vt_1(t_n, q_2C_9\Psi_0(t_{n}))}\Big) \quad \text{for all} \;\; n \ge N(\eps).
	\end{equation*}
	Therefore, by \eqref{e:Psi0}, we also obtain $\sum_{n=1}^\infty I_{n,2}<\infty$.

	In the end, 
	by the Borel-Cantelli lemma, we deduce that  $\limsup_{t \to 0} d(x, X_t)/\Psi(t)\le \eps$, $\P^x$-a.s. for every $\eps>0$. This contradicts \eqref{e:limsuplaw1}. Hence, we conclude that $\liminf_{t \to 0}f(t)<\infty$. \qed

	\appendix
	\section{Proofs of Theorems \ref{t:SDsup0} and \ref{t:SDsup2}}\label{s:A}

	Recall that the function $F$ is  increasing continuous on $(0,\infty)$ satisfying \eqref{e:Fscale} and the function $F_1$ is defined in \eqref{e:F1}.
	
The next lemma follows from \cite[(3.42) and Lemma 3.19]{GT12} and the inequality $e^{-x} \le \beta^{\beta} x^{-\beta}$, $x, \beta>0$.
	\begin{lem}\label{l:expF}
		(i) For any constants $k_1, k_2>0$, there exists a constant  $c=c(k) \ge 1$ such that
		\begin{equation*}
			\exp\big(-k_1F_1(r,s)\big) \le c\big(s/F(r)\big)^{k_2} \quad \text{for all} \;\; r,s>0.
		\end{equation*}
		
		\noindent (ii) $F_1(ar, s) \ge aF_1(r,s)$ for all $a \ge 1$ and  $r,s>0$.
	\end{lem}
	
 By a standard argument, using Lemma \ref{l:expF}(i), one can see from \eqref{e:diffusion} that there exist constants $\eps_1,\eta_1 \in (0,1)$, $c_1>0$ such that   for all $x \in M$ and $r \in (0, \eps_1R_1)$, 
	\begin{equation}\label{e:ndlq}
		q^{B(x,r)}(t,y,z) \ge  \frac{c_1}{V(x, F^{-1}(t))} \quad \text{for all} \;\; t \in (0, F(\eta_1  r)], \; y,z \in B(x, \eta_1 F^{-1}(t)).
	\end{equation}
	
	Recall the definitions of $\Phi, \Theta$ and $\Pi$ from \eqref{e:subaux1}. 
	We note that by $\mathrm{VRD}_{\bar R}(M)$,  there exists a constant $\ell_1>1$ such that
	\begin{equation}\label{e:dset2}
		V(x, r) \ge 2V(x, r/\ell_1) \quad \text{for all} \;\; x \in M, \; r \in (0, \bar R).
	\end{equation}	
	
	\begin{lem}\label{l:subSP}
	There exist constants $\eps \in (0,1)$, $a_1,a_2>0$ and $k_1\ge1$ such that for all $x \in M$ and $0<r<\eps R_1$,
		\begin{equation}\label{e:subSP}
		e^{-a_1n} \le	\P^x(\tau_{B(x,r)} \ge n \phi(k_1r) ) \le  e^{-a_2n} \quad \text{for all} \;\; n \ge 1.
	\end{equation}
Therefore, $\mathrm{E}_{\eps R_1}(\phi, M)$ holds.
	\end{lem}
	\pf Let $\eps_1, \eta_1$ be the constants in  \eqref{e:ndlq}. By \eqref{e:ndlq}, \eqref{e:dset2} and  $\L(F, \gamma_1, c_L)$,  similar calculations  to  \eqref{e:E_0-2} show  that  there exist constants $\delta_1 \in (0, \eps_1]$,  $c_1 \in (0,1)$ and  $k_1\ge 1$ such that for every $x \in M$ and $r \in (0, \delta_1 R_1)$,
	$$\P^y\big(\sup_{0<s \le 2^{-1}F(k_1  r)}d(x,Z_s)> r\big) \ge c_1 \quad \text{for all} \;\; y \in B(x,r).$$ 
	Besides, by  \cite[Proposition 2.4]{Mi16}, it holds that  $\P(S_t \ge 2^{-1}/\phi_1^{-1}(t^{-1})) \ge 1-e^{-1/2}$ for all $t>0$. Thus,  for every $x \in M$ and  $r \in (0, \delta_1R_1)$,
	\begin{align}\label{e:subSP1}
		\inf_{y \in B(x,r)}\P^y\big( \tau_{B(x,r)}\le \Phi(k_1  r)\big) &\ge \inf_{y \in B(x,r)} \P^y\big( \sup_{0<s \le \Phi(k_1  r)}d(x, Z_{S_s}) > r \big)  \nn\\
		&\ge c_1 \P\big(S_{\Phi(k_1  r)} \ge 2^{-1}F(k_1  r)\big) \ge c_1(1-e^{-1/2}).
	\end{align}
Then using Markov property,  we deduce that the upper bound in \eqref{e:subSP} holds as in \eqref{e:E_0-4}.

	On the other hand, by \eqref{e:ndlq} and \VRDM, we see that for all $x \in M$ and $r \in (0, \delta_1R_1)$, 
	\begin{equation*}
		\P^y\big(Z_{F(\eta_1 r)}^{B(x,r)} \in B(x,\eta_1^2 r)\big) \ge \frac{c_2 V(x, \eta_1^2 r)}{V(x, \eta_1 r)} \ge c_3 \quad \text{for all}  \,\; y \in B(x, \eta_1^2r).
	\end{equation*}
	Moreover, by  \cite[Lemmas 2.11 and 2.4(ii)]{CK20}, there exists a universal constant $p_0\in (0,1)$ such that $\P(S_t \le 1/\phi_1^{-1}(t^{-1})) \ge p_0$ for all $t>0$.
	Since $S$ is independent of $Z$ and $F(r) = 1/\phi_1^{-1}(\Phi(r)^{-1})$,  using the Markov property and $\U(F, \gamma_2, c_U)$, we get that for all $x \in M$, $r \in (0, \delta_1R_1)$ and $n \ge 1$,
	\begin{align*}
		&\P^x\big( \tau_{B(x,r)} \ge n\Phi(k_1 r) \big) \ge \P^x\big( Z^{B(x,r)}_{S_{n\Phi(k_1 r)}} \in B(x,r) \big) \ge  \P\big(S_{n\Phi(k_1r)} \le nF(k_1r)\big) \, \P^x\big( Z^{B(x,r)}_{nF(k_1r)} \in B(x,r)\big)\\[2pt]
		&\ge \P\big(S_{\Phi(k_1r)} \le F(k_1r)\big)^n  \Big[\inf_{y \in B(x, \eta_1^2 r)}\P^y\big(Z^{B(x,r)}_{F(\eta_1r)}  \in B(x, \eta_1^2 r) \big)\Big]^{\lceil n F(k_1r)/F(\eta_1 r)\rceil}  \ge (p_0c_3^{c_U (k_1/\eta_1)^{\gamma_2}})^n.
	\end{align*}
In the third inequality above, we used the fact that $S_{n\Phi(k_1r)}-S_{(n-1)\Phi(k_1r)}$ and $S_{\Phi(k_1r)}$ have the same law for all $n \ge 1$. This finishes the proof for \eqref{e:subSP}. Then  $\mathrm{E}_{\eps R_1}(\phi, M)$ follows from \eqref{e:SP-E-1} and \eqref{e:SP-E-2} with $k_1$ instead of $C_1'$.  \qed

	\begin{lem}\label{l:subTail}
		\noindent (i) $\mathrm{TJ}_{R_1}(\Theta,  \le,M)$ holds. In particular, by \eqref{e:aux2}, $\mathrm{TJ}_{R_1}(\Phi, \le,M )$ also holds.
		
		\noindent(ii) If $\U_{F(R_2)}(w^{-1}, \alpha_2, c_0)$ holds, then there exists $\eps \in (0,1)$ such that  $\mathrm{TJ}_{\eps (R_1 \land R_2)}(\Pi, M)$ holds.
		
		\noindent(iii) If  $\bar R=R_1=\infty$ and $\U^{F(R_3)}(w^{-1}, \alpha_4, c_0)$  holds, then   $\mathrm{TJ}^{R_3}(\Pi)$ holds.
	\end{lem}
	\pf  Since 
	 \VRDM, $\U(F, \gamma_2, c_U)$ and $\L(F, \gamma_1, c_L)$  hold with $\gamma_1>1$, by \eqref{e:diffusion} and Lemma \ref{l:expF},   we see that for all $x \in M$, $r \in (0, R_1)$ and $s \in (0, F(r))$ (cf. \cite[Lemma 3.9]{Ba98}),
	\begin{equation*}
		\begin{split}
			&\int_{M_\partial \setminus B(x,r)} q(s,x,y)\mu(dy) = \sum_{k=0}^\infty \int_{ B(x, 2^{k+1}r) \setminus B(x,2^kr)} q(s,x,y) \mu(dy) \\
			&\le \frac{c_1\exp\big(-2^{-1}c_2 F_1(r,s) \big)}{V(x, F^{-1}(s))} \sum_{k=0}^\infty \mu\big( B(x, 2^{k+1}r) \setminus B(x,2^kr)\big)\exp\big(-2^{-1}c_2 F_1( 2^kr,s) \big)\\
			& \le \frac{c_3V(x, r)\exp\big(-2^{-1}c_2 F_1(r,s) \big)}{V(x, F^{-1}(s))} \sum_{k=0}^\infty 2^{kd_2}\Big(\frac{s}{2^kF(r)}\Big)^{2d_2} \\
			&\le c_4\Big(\frac{F(r)}{s}\Big)^{d_2/\gamma_1}\Big(\frac{s}{F(r)}\Big)^{2d_2}\exp\big(-2^{-1}c_2 F_1(r,s) \big)  \le c_4\exp\big(-2^{-1}c_2 F_1(r,s) \big).
		\end{split}
	\end{equation*}
	By \cite[Lemma 3.9(ii)]{BKKL19b}, there exists  $c_5>0$ such that  $F_1(r, s) \le c_5$ for all $r \in (0, R_1)$ and $s \ge F(r)$. Hence, by taking $c_4$ larger than $e^{c_2c_5/2}$, we get that $\int_{B(x,r)^c} q(s,x,y)\mu(dy) \le c_4\exp(-2^{-1}c_2 F_1(r,s) )$  for all $x \in M$, $r \in (0, R_1)$ and $s>0$. In view of \eqref{e:PJ} and Tonelli's theorem, it follows that for all $x \in M$ and $r \in (0, R_1)$, 
	\begin{align*}
		J\big(x, M_\partial \setminus B(x,r)\big) &= \int_0^\infty \int_{M_\partial \setminus B(x,r)}q(s,x,y) \mu(dy) \nu(ds) \le c_4\int_0^\infty \exp \big(- 2^{-1}c_2F_1(r,s) \big)\nu(ds)\\
		&\le c_4\int_{(0, F(r)]} \exp \big(- 2^{-1}c_2F_1(r,s) \big)\nu(ds)+ c_4 w(F(r))=:J_1 + c_4w(F(r)).
	\end{align*}
	
 (i) By \eqref{e:aux2} and Lemma \ref{l:expF}(i), using the inequality $1-e^{-\lambda}-\lambda e^{-\lambda} \ge (2e)^{-1}(\lambda^2 \land 1)$ for $\lambda \ge 0$, we get that
	\begin{align*}
		J_1 + c_4 w(F(r))&\le c_5 \int_{(0, F(r)]} \Big(\frac{s}{F(r)}\Big)^{2}\nu(ds) + 
		c_4w(F(r)) \le (c_4+c_5) \int_0^\infty \Big( \Big(\frac{s}{F(r)}\Big)^{2} \land 1\Big) \nu(ds) \\
		&\le 2e(c_4+c_5) \int_0^\infty \big(1-e^{-s/F(r)}-sF(r)^{-1}e^{-s/F(r)}\big)\nu(ds)=\frac{2e(c_4+c_5)}{\Theta(r)}.
	\end{align*}

	(ii) By Lemma \ref{l:expF}(i), $\U_{F(R_2)}(w^{-1}, \alpha_2, c_0)$ and the integration by parts,  for all  $r\in (0, R_2)$,
	\begin{align}\label{e:TailJ1}
		J_1 &\le  c_6\int_{(0, F(r)]} \Big(\frac{s}{F(r)}\Big)^{\alpha_2+1} \nu(ds) \le  c_6(\alpha_2+1)F(r)^{-\alpha_2-1}\int_{(0, F(r)]}s^{\alpha_2}w(s)ds \nn\\
		& \le c_7 F(r)^{-\alpha_2-1}F(r)^{\alpha_2} w(F(r)) \int_{(0,F(r)]}ds = c_7w(F(r)).
	\end{align}
Next, let  $c_8>0$ be such that $w\big(F(c_8\ell_1 (R_1 \land R_2))\big) \ge 2 w\big(F(R_1)\big)$ where $\ell_1$ is the constant in \eqref{e:dset2}. By \eqref{e:PJ}, \eqref{e:diffusion}, \eqref{e:dset2},  $\U(F, \gamma_2, c_U)$ and $\U_{F(R_2)}(w^{-1}, \alpha_2, c_0)$, 
	 we get that  for all $x \in M$ and $r \in (0, c_8 (R_1 \land R_2))$,
	\begin{align}\label{e:J1int2}
		J\big(x, M_\partial \setminus B(x,r)\big) &\ge c_9\int_{F(\ell_1 r)}^{F(R_1)} \int_{B(x, F^{-1}(s) ) \setminus B(x,r)} \frac{\mu(dy)}{V(x, F^{-1}(s))} \nu(ds) \nn\\
		&\ge  \frac{c_9}{2}\int_{F(\ell_1 r)}^{F(R_1)} \nu(ds) \ge \frac{c_9}{4} w\big(F(\ell_1 r)\big) \ge  c_{10}w\big(F(r)\big).
	\end{align}
	
	 (iii) By similar calculations to \eqref{e:TailJ1},  since $\U^{F(R_3)}(w^{-1}, \alpha_4, c_0)$ holds, we see that for all $r>R_3$,
	\begin{equation*}
		\begin{split}
			J_1 &\le  c_{11}(\alpha_4+1)F(r)^{-\alpha_4-1}\Big(\int_{(0, F(R_3)]}s^{\alpha_4}w(s)ds  + \int_{(F(R_3), F(r)]}s^{\alpha_4}w(s)ds  \Big)\\[2pt]
			& \le  c_{12}F(r)^{-\alpha_4-1}\big(1 + F(r)^{\alpha_4+1}w(F(r)) \big) \le c_{13}w(F(r)).
		\end{split}
	\end{equation*}
	Moreover, by similar calculations to \eqref{e:J1int2}, we also have that for all $r>R_3$,
	\begin{equation*}
		J\big(x, M_\partial \setminus B(x,r)\big) \ge c_{14}\int_{F(\ell_1 r)}^{\infty} \frac{\mu\big(B(x, F^{-1}(s)) \setminus B(x,r)\big)}{V(x, F^{-1}(s))} \nu(ds) \ge c_{15}w\big(F(\ell_1 r)\big) \ge c_{16} w\big(F(r)\big). 
	\end{equation*} \qed
	
	\begin{lem}\label{l:subNDL} {\rm (cf. \cite[Lemma 4.5]{BKKL19b})}	
		\noindent(i) If $\L^1(\phi_1, \alpha_1, c_0)$  holds, then  $\mathrm{NDL}_{\eps_1(R_1 \wedge 1)}(\Phi, M)$  holds where $\eps_1 \in (0,1)$ is the constant in \eqref{e:ndlq}.
		
		\noindent(ii) If $\bar R=R_1=\infty$ and $\L_1(\phi_1, \alpha_3, c_0)$  holds, then  $\mathrm{NDL}^{R}(\Phi)$  holds for some $R\ge1$.
	\end{lem}
	\pf (i)  Since  $\L^1(\phi_1, \alpha_1, c_0)$ holds and $\phi_1^{-1}$ is increasing, we see that, for any $a>0$, there exists $c_1>0$ such that $\U^{a}(\phi_1^{-1}, 1/\alpha_1, c_1)$ holds (see \cite[Remark 2.1]{BKKL19a}). Thus, by \cite[Proposition 2.4]{Mi16},  there exist constants $c_2>1$, $c_3>0$ such that for all $t \in (0, \Phi(R_1 \wedge 1))$,
	\begin{equation}\label{e:meansub}
		\P\big( S_t \in [2^{-1}\phi_1^{-1}(t^{-1})^{-1}, c_2\phi_1^{-1}(t^{-1})^{-1} ]\big) \ge c_3.
	\end{equation}
	Let $\eta_2:=(2c_2/c_L)^{1/\gamma_1} \eta_1$.
	Here, $c_L, \gamma_1$ and $\eta_1$ are the constants in  $\L(F, \gamma_1, c_L)$ and \eqref{e:ndlq}.
	By \eqref{e:ndlq} and \eqref{e:meansub}, since \VRDM \ holds,  we get that for all $r \in (0, \eps_1(R_1 \land 1))$,  $x \in M$ and $y,z \in B(x,\eta_2^2r)$,
	\begin{align*}
		p^{B(x,r)}(\Phi(\eta_2 r), y,z) &\ge \int_0^\infty q^{B(x,r)}(s,y,z) \P(S_{\Phi(\eta_2r)} \in ds) \\
		&\ge c_3\inf\big\{ q^{B(x,r)}(s,y,z) : 2^{-1}\phi_1^{-1}(\Phi(\eta_2r)^{-1})^{-1} \le  s \le c_2\phi_1^{-1}(\Phi(\eta_2r)^{-1})^{-1}\big\}\\
		& \ge \frac{c_4}{V\big(x, F^{-1}(c_2 F(\eta_2 r))\big)}  \ge \frac{c_5}{V(x, r)}.
	\end{align*}
	Indeed, since $\phi_1^{-1}(\Phi(\eta_2r)^{-1})^{-1}=F(\eta_2 r)$ and $\L(F, \beta_1, c_L)$ holds,  we have  $c_2\phi_1^{-1}(\Phi(\eta_2r)^{-1})^{-1}  \le c_2c_L^{-1}(\eta_2/\eta_1)^{\beta_1}F(\eta_1 r) \le F(\eta_1 r)$ and $\eta_1 F^{-1}(2^{-1}F(\eta_2r)) \ge \eta_1 \eta_2 (c_L/2)^{1/\beta_1}r > \eta_2^2r$. 
	
(ii) Similarly, since $\L_1(\phi_1, \alpha_3, c_0)$ holds, for all $a>0$,  $\U_{a}(\phi_1^{-1}, 1/\alpha_3, c_6)$ holds with some $c_6>0$. Thus, by applying \cite[Proposition 2.4]{Mi16} again,  we get that  \eqref{e:meansub} holds for all $t \in (\Phi(1), \infty)$. Define $\eta_2$ as in (i). Then by the same argument as that for (i), we conclude that  $\mathrm{NDL}^{1/\eta_2}(\Phi)$ holds.\qed
	
	Now, before giving proofs of Theorems  \ref{t:SDsup0} and \ref{t:SDsup2}, we mention  that although $\Pi$ is not an increasing continuous function in general, under $\U_{F(R_2)}(w^{-1}, \alpha_2, c_0)$ (resp. $\U^{F(R_3)}(w^{-1}, \alpha_4, c_0)$), there exists an increasing continuous function $\Pi_1$ such that $\Pi(r) \asymp \Pi_1(r)$ for all $r \in (0, R_2/2)$ (resp. $r \in (R_3, \infty)$). Since our tail condition is invariant even if we change the scale function $\psi$ to some other function comparable with $\psi$, this observation allows us to assume that $\Pi$ is  increasing and continuous without loss of generality.
	Indeed, for example,  define an  function $\Pi_1$ on $(0,\infty)$ as 
	\begin{equation}\label{e:continu}
		\Pi_1(r):= \begin{cases}
			\big(\int_{1}^\infty e^{1-u}w(  F(ru)) du + e^{-r}w(F(2^{-1}R_2))\big)^{-1},  &\mbox{ if $\U_{F(R_2)}(w^{-1}, \alpha_2, c_0)$ holds},\\[5pt]
			\big(\int_{1}^\infty e^{1-u} w(F(ru)) du\big)^{-1}, \; &\mbox{ if $\U^{F(R_3)}(w^{-1}, \alpha_4, c_0)$ holds}.
		\end{cases}
	\end{equation}
	By the monotone property of $w$, we have  $\int_{1}^\infty e^{1-u}w(F(ru)) du \le w(F(r))$ and $\int_{1}^\infty e^{1-u}w(F(ru)) du  \ge w(F(2r)) \int_1^2 e^{1-u}du $ for all $r>0$.
	Note that under $\U_{F(R_2)}(w^{-1}, \alpha_2, c_0)$ (resp. $\U^{F(R_3)}(w^{-1}, \alpha_4, c_0)$), $w(F(r)) \asymp w(F(2r))$ for all $r \in (0, R_2/2)$ (resp. $r \in (R_3, \infty)$).   Thus, $\Pi_1(r) \asymp \Pi(r)$ for  $r \in (0, R_2/2)$ (resp. $r \in (R_3, \infty)$) under $\U_{F(R_2)}(w^{-1}, \alpha_2, c_0)$ (resp. $\U^{F(R_3)}(w^{-1}, \alpha_4, c_0)$).

	\bigskip

	\noindent \textbf{Proof of Theorem \ref{t:SDsup0}.} (i) Suppose that $\limsup_{r \to \infty}\phi_1(r)/H_0(r)<\infty$. Then   $\phi_1(\lambda) \asymp H_0(\lambda)$ for all $\lambda \ge 1$. Since  $\lim_{r \to \infty}r/H_0(r) \ge \lim_{r \to \infty}r/\phi_0(r)=\infty$, we have that $b=0$ and hence $\phi_1(\lambda)=\phi_0(\lambda)$. Note that $(\lambda^{-1}\phi_0(\lambda))'=-\lambda^{-2}H_0(\lambda)$. Thus, it holds that for $\lambda \ge 1$,
	\begin{equation}\label{e:rv}
		\lambda^{-1}H_0(\lambda) \asymp \lambda^{-1}\phi_1(\lambda)= \lambda^{-1}\phi_0(\lambda) = \int_\lambda^\infty u^{-2}H_0(u)du.
	\end{equation} 
	According to \cite[Corollary 2.6.4]{BGT89}, \eqref{e:rv} implies $\U^1(H_0, 1-\eps, c_1)$  with some constants $\eps, c_1>0$. Hence, by \cite[Lemma 2.6]{Mi16}, there exists $c_2>0$ such that  $w(r) \asymp H_0(r^{-1}) \asymp \phi_0(r^{-1})$ for $r \in (0, c_2)$. In particular,  $\limsup_{r \to \infty}\phi_1(r)/w(r^{-1})<\infty$, and  $\U_{c_2}(w^{-1},1, c_3)$ holds. Finally, by Lemmas \ref{l:subSP} and \ref{l:subTail}(ii), we conclude the result from Theorem \ref{t:limsup0-1}(i).

(ii) By \eqref{e:aux} and Lemmas \ref{l:subSP}, \ref{l:subTail}(i,ii)  and \ref{l:subNDL}(i),  we obtain the result from Theorems \ref{t:limsup0-1}(ii) and \ref{t:supprecise1} (see the paragraph below Lemma \ref{l:subNDL}). \qed

	\noindent \textbf{Proof of Theorem \ref{t:SDsup2}.} (i) Suppose that $\limsup_{r \to 0} \phi_1(r)/H_0(r)<\infty$. Since $\lim_{r \to 0} \phi_0(r)/r \in (0, \infty]$ and $\phi_0(\lambda) \ge H_0(\lambda)$ for all $\lambda$,  \eqref{e:rv} holds for $\lambda \in (0,1]$. Let $f(s):=sH_0(1/s)$. By the change of the variables, we have $\int_0^s u^{-1}f(u)du = s \phi_0(1/s) \asymp f(s)$ for  $s \ge 1$. By \cite[Corollary 2.6.2]{BGT89}, it follows that $\L^1(f, \eps, c_1)$ holds with some constants $\eps, c_1>0$ so that $\U_1(H_0,  1-\eps, c_1^{-1})$ holds. Thus, by \cite[Lemma 2.1(iii) and (iv)]{CK20},  $w(r) \asymp H_0(r^{-1}) \asymp \phi_1(r^{-1})$ for all $r \ge 1$. Therefore, $\limsup_{r \to 0} \phi_1(r)/w(r^{-1})<\infty$ and $\U^{1}(w^{-1},1, c_2)$ holds. Now we get the result from Theorem \ref{t:limsupinf-1}(i) and Lemmas \ref{l:subSP} and \ref{l:subTail}(iii).
	
	(ii) By Lemmas \ref{l:subSP}, \ref{l:subTail}(i,iii) and \ref{l:subNDL}(ii), the result follows from Theorems \ref{t:limsupinf-1}(ii)  and \ref{t:supprecise2} (see the paragraph below Lemma \ref{l:subNDL}). \qed
	
	\section{Proof of Theorem \ref{t:Hunt0}}\label{s:B}
	Throughout this subsection, we assume that Assumption L holds, and that, without loss of generality, $\beta_1 \le \gamma_2$. Here, $\beta_1$ is the constant in (L1) and $\gamma_2$ is the one in \eqref{e:Fscale}. When $R_4 < \infty$, we extend $\psi_0$ to all $r \in (0,\infty)$ by setting $\psi_0(r) =\psi_0(R_4-)R_4^{-2\gamma_2} r^{2\gamma_2}$ for $r \ge R_4$. Note that	\begin{equation}\label{e:scalepsi0}
		\text{$\L(\psi_0, \beta_1, c_1)$ \, and \,  $\U(\psi_0, \beta_2 + 2\gamma_2, c_2)$ \,  hold.}
	\end{equation}

	\begin{lem}\label{l:JX}  $\Tail_{(\bar R \, \land R_4)/\ell_1}(\psi_0,\, \sU)$  holds where $\ell_1$ is the constant in \eqref{e:dset2}.
	\end{lem}
	\pf  Choose any $x \in \sU$ and $0<\ell_1 r<\bar R \land  R_4 \land (\kappa_0 \updelta_\sU(x))$. Then, by (L2), (L3), \eqref{e:scalepsi0} and \cite[Lemma 2.1]{CKW16a}, since $\mathrm{VRD}_{\bar R}(M)$ holds, we see that
	\begin{equation*}
		J_X\big(x, M_\partial \setminus B(x,r) \big) \le  \int_{\{r \le d(x,y) <  R_4 \land (\kappa_0 \updelta_\sU(x))\}} \frac{c_1 \mu(dy)}{V(x, d(x,y))\psi_0(d(x,y))} + \frac{c_2}{\psi_0(r)} \le \frac{c_3}{\psi_0(r)}
	\end{equation*}
	and by (L2), 	\eqref{e:dset2}, $\mathrm{VRD}_{\bar R}(M)$ and \eqref{e:scalepsi0},
	\begin{equation*}
		J_X\big(x, M_\partial \setminus B(x,r) \big) \ge  \int_{\{r \le d(x,y) < \ell_1 r\}} \frac{c_4 \mu(dy)}{V(x, d(x,y))\psi_0(d(x,y))} \ge \frac{c_4 V(x,r)}{V(x, \ell_1 r)\psi_0(\ell_1 r)} \ge \frac{c_5}{\psi_0(r)}.
	\end{equation*}
	\qed
	
	In the followings,  we establish local versions of stability results obtained in \cite{GH14, 
		 CKW16a, CKW19}. Then, by using those and comparing the process $X$ with a suitable subordinate process, we show that $\mathrm{NDL}_{R}(\Phi_1, \sU)$  holds with some $R \in (0,\infty]$ (Proposition \ref{p:NDLX}).

	Let $(\EE,\FF)$ be a regular Dirichlet form on $L^2(M;\mu)$ without killing term and having the jumping kernel $J(dx,dy)$ (Cf,  \cite[(1.1)]{CKW19}), and $\bar X = \{\bar  X_t, t \ge 0; \bar {\P}^x, x \in M \setminus \bar \sN\}$ be a Hunt process associated with $(\EE, \FF)$ where $\bar \sN$ is the properly exceptional set for $\bar X$. Write the strongly part of $(\EE, \FF)$ as $\EE^{(c)}$, and let $\Gamma$ and $\Gamma^{(c)}$ be the energy measures for $(\EE,\FF)$ and  $(\EE^{(c)}, \FF)$,  respectively.  Set $\FF_b := \FF \cap L^\infty(M;\mu)$. For an open set $D \subset M$, denote by $\FF_D$ the $\EE_1$-closure of $\FF\cap C_c(D)$ in $\FF$. 
	
	Here,  we introduce local (and interior) versions of a \textit{cut-off Sobolev inequality}, the \textit{Faber-Krahn inequality} and the \textit{(weak) Poincar\'e inequality} for $(\EE,\FF)$. See  \cite{CKW19} for global versions.

	For open sets $D \subset U \subset M$ with $\overline D \subset U$, we say a measurable function $\varphi$ on $M$ is a \textit{cut-off function} for $D \subset U$, if $\varphi = 1$ on $D$, $\varphi=0$ on $M \setminus U$ and $0 \le \varphi \le 1$ on $M$.

	\begin{definition}\label{d:CS}
		{\rm	For an open set $U \subset M$ and $R_0 \in (0,\infty]$, we say that a \textit{local cut-off Sobolev inequality}  $\CS_{R_0}(\phi,U)$ (for $(\EE,\FF)$) holds (with $C$)  if there exist constants $C \in (0,1)$,  $a \in (0,1]$ and $c>0$ such that for almost all $x \in U$, all $0<r \le R <R_0 \land (C\updelta_U(x))$ and any $f \in \FF$, there exists a cut-off function $\varphi \in \FF_b$ for $B(x_0,R) \subset B(x_0,R+r)$ satisfying the following inequality:
			\begin{align*}
				\int_{B(x_0,R+(1+a)r)} f^2 d\Gamma(\varphi,\varphi) &\le  \frac{c}{\phi(r)} \int_{B(x_0,R+(1+a)r)} f^2 d\mu +  c\int_{B(x_0,R+r)} \varphi^2 d\Gamma^{(c)}(f,f)\\
				&\quad + c\int_{B(x_0,R+r) \times B(x_0,R+(1+a)r)} \varphi^2(x) (f(x) - f(y))^2 J(dx,dy). 
			\end{align*}	
		}
	\end{definition}
	
	\begin{definition}\label{d:FK}
		{\rm	For an open set $U \subset M$ and $R_0 \in (0,\infty]$, we say that a \textit{local Faber-Krahn inequality} $\FK_{R_0}(\phi,U)$  (for $(\EE,\FF)$) holds (with $C$ and $\nu$) if there exist constants  $C, \nu \in (0,1)$ and  $c>0$  such that for all $x \in U$, $0<r<R_0 \land (C\updelta_U(x))$ and any open set $D \subset B(x,r)$,
			\begin{equation}\label{e:FK}
				\lambda_1(D) \ge \frac{c}{\phi(r)} \left( \frac{V(x,r)}{\mu(D)}\right)^\nu,
			\end{equation}
			where $\lambda_1(D):= \inf \{ \EE(f,f) : f \in \FF_D \mbox{ with } \Vert f \Vert_2 = 1  \}$.
		}
	\end{definition}
	
	\begin{definition}\label{d:PI}
		{\rm
			For an open set $U \subset M$ and $R_0 \in (0,\infty]$, we say that a \textit{(weak) Poincar\'e inequality} $\PI_{R_0}(\phi,U)$  (for $(\EE,\FF)$)  holds (with $C$) if there exist constants $C \in (0,1)$  and $c>0$ such that for all $x \in U$, $0<r<R_0 \land (C\updelta_U(x))$ and  any $f \in \FF_b$,
			\begin{equation*}
				\int_{B(x,r)} (f- \bar{f}_{B(x,r)})^2 d\mu \le c\phi(r) \Big(\int_{B(x,r/C)} d\Gamma^{(c)}(f,f) + \int_{B(x,r/C) \times B(x,r/C)} (f(y) - f(x))^2 J(dx,dy) \Big), 
			\end{equation*}
			where $\bar{f}_{B(x,r)} := \frac{1}{V(x,r)} \int_{B(x,r)} f d\mu$ is the average value of $f$ on $B(x,r)$.
		}
	\end{definition}
	
	\smallskip

	We observe that $\FK_{R_0}(\phi,U)$ is equivalent to an existence and upper bound \eqref{e:Dup} of Dirichlet heat kernel for small balls contained in $U$ (and a local Nash inequality). When  $U=M$ and $R_0=\infty$, the following result was established in  \cite[Section 5]{GH14} (when $\phi(r)=r^\beta$ for some $\beta>0$) and \cite[Proposition  7.3]{CKW16a}.  Even though the proof is similar, we give the full proof for the reader's convenience.
	
	For an open set $D \subset M$, we write  $\bar \tau_D$, 
	and $\bar p^D(t,x,y)$ for the first exit time of $\bar X$ from $D$ and heat kernel of killed process $\bar X^D$, respectively.
	
	\begin{lem}\label{l:7.3}
		Let $U \subset M$ be an open set, and $R \in (0,\infty]$ and $\nu,\kappa \in (0,1)$ be constants. Then, the followings are equivalent.
		
		\smallskip
		
		\setlength{\leftskip}{0.4cm}
		
		\noindent {\rm	(1)} $\FK_{R}(\phi,U)$ holds with  $\kappa$ and $\nu$. 
		
		\noindent {\rm	(2)} There exists a constant $K_1>0$ such that for all $x \in U$ and  $0<r<R \land (\kappa\updelta_U(x))$,
		\begin{equation}\label{e:Nashtype} \frac{V(x,r)^\nu}{\phi(r)} \Vert u \Vert_2^{2+2\nu}  \Vert u \Vert_1^{-2\nu}  \le K_1 \EE(u,u), \quad  \;\;  \forall \, u \in \FF_{B(x,r)}.
		\end{equation}

		\noindent {\rm	(3)} There exists a constant $K_2>0$ such that for all $x \in U$ and $0<r<R \land (\kappa\updelta_U(x))$, the (Dirichlet) heat kernel $\bar p^{B(x,r)}(t, \cdot, \cdot)$ exists and satisfies that
		\begin{equation}\label{e:Dup}
			\esssup_{y,z \in B(x,r)} \bar p^{B(x,r)}(t,y,z) \le \frac{K_2}{V(x,r)} \Big( \frac{\phi(r)}{t} \Big)^{1/\nu}, \quad \;\;  \forall \, t>0.
		\end{equation}
	\end{lem}
	\pf Choose any $x \in U$ and $0<r<R \land (\kappa\updelta_U(x))$, and denote by $B = B(x,r)$. 
	
	\smallskip
	
	(1) $\Rightarrow$ (2). We first assume that $0 \le u \in \FF \cap C_c(B)$. Set $B^s := \{y \in B : u(y) >s \}$ for $s>0$. Then $B^s \subset B$ is open because $u$ is continuous. By the Markovian property of $\EE$, since $\FK_{R}(\phi,U)$ holds,  we have that, for any $t>s>0$,
	\begin{align}\label{e:7.3.0}
		\EE(u,u) &\ge \EE\big((u - t)_+, (u-t)_+\big)  \ge \lambda_1(B^s) \int_{B^s} ((u-t)_+)^2 d\mu \nn\\
		&\ge \frac{c}{\phi(r)} \Big( \frac{V(x,r)}{\mu(B^s)}\Big)^\nu\int_{B^s} ((u-t)_+)^2 d\mu.
	\end{align}
	Observe that since $u \ge 0$, we have $((u-t)_+)^2 \ge u^2-2tu$. Thus, we get that $\int_{B^s} ((u-t)_+)^2 d\mu  = \int_{M} ((u-t)_+)^2 d\mu  \ge  \int_M (u^2 - 2tu) d\mu =  \Vert u \Vert_2^2 - 2t\Vert u \Vert_1$. Moreover, by the definition of $B^s$, it holds that $\mu(B^s) \le s^{-1} \int_{B^s} u d\mu \le s^{-1} \Vert u \Vert_1$. Therefore, we get from \eqref{e:7.3.0} that for any $t>s>0$,
	\begin{equation*}
		\EE(u,u) \ge \frac{c  V(x,r)^\nu }{\phi(r)}\Big( \frac{s}{\Vert u \Vert_1} \Big)^\nu ( \Vert u \Vert_2^2 - 2t\Vert u \Vert_1).
	\end{equation*}
	By letting $t \downarrow s$ and taking $s = \Vert u \Vert_2^2/(4\Vert u \Vert_1)$ in the above inequality, we obtain that 
	\begin{align*}
		\EE(u,u) &\ge \frac{cV(x,r)^\nu  }{\phi(r)} \Big( \frac{\Vert u \Vert_2}{4\Vert u \Vert_1} \Big)^{2\nu}  \frac{\Vert u \Vert_2^2}{2}   =  \frac{c 2^{-4\nu-1}V(x,r)^\nu}{\phi(r) } \Vert u \Vert_2^{2+2\nu} \Vert u \Vert_1^{-2\nu}.
	\end{align*}
	Thus, \eqref{e:Nashtype} holds for any $0 \le u \in \FF \cap C_c(B)$. Then since $\EE(|u|,|u|) \le \EE(u,u)$ by the Markovian property, 
	we can see that \eqref{e:Nashtype} holds for any signed $u \in \FF \cap C_c(B)$. 
	
	Now, we consider a general $u \in \FF_B$. By the definition of $\FF_B$, there exists a sequence $(u_n)_{n \ge 1} \subset \FF \cap C_c(B)$ such that $\EE(u_n - u, u_n - u) + \Vert u_n - u \Vert_2^2 \rightarrow 0$. Since $\mu(B) < \infty$, by Cauchy inequality, it also  holds $\Vert u_n - u \Vert_1 \le \sqrt{\mu(B)} \Vert u_n - u \Vert_2 \to 0$. Thus, since  \eqref{e:Nashtype} holds for each $u_n \in \FF \cap C_c(B)$, by passing to the limit as $n \to \infty$, we conclude that \eqref{e:Nashtype} holds for any $u \in \FF_B$.
	
	\smallskip
	
	(2) $\Rightarrow$ (3). Let $0 \le f \in L^2(B;\mu)$ with $\Vert f \Vert_1 = 1$.
	Set $u_t(\cdot) =\bar P_t^B f(\cdot):=\E^{\cdot}f(\bar X^B_t)$ and denote $m(t) = \Vert u_t \Vert_2^2$ for $t>0$. Then since $(\EE,\FF)$ is regular, we have $u_t \in \FF_B$ for every $t>0$. Hence, by (2), because $\bar P_t^B$ is a $L^1$-contraction,  it holds that
	$$- \frac{K_1}{2} \frac{dm}{dt} =K_1\EE (u_t, u_t) \ge  \frac{V(x,r)^\nu}{\phi(r)} m(t)^{1+\nu} \lVert u \rVert_1^{-2\nu}\ge  \frac{V(x,r)^\nu}{\phi(r)} m(t)^{1+\nu} =:a(r) m(t)^{1+\nu}.$$
	It follows $\frac{d}{dt}(m(t)^{-\nu}) \ge 2\nu^{-1}K_1^{-1}a(r)$ so that $m(t) \le 2^{-1/\nu}\nu^{1/\nu}K_1^{1/\nu}a(r)^{-1/\nu}t^{-1/\nu}$. Hence, we obtain $\lVert \bar P_t^B \rVert_{L^1 \to L^2} \le c_1 a(r)^{-1/(2\nu)}t^{-1/(2\nu)}$. Thus, by \cite[Lemma 3.7]{GH14}, we conclude that (3) holds.
	
	\smallskip
	
	(3) $\Rightarrow$ (1) Let $D \subset B$ be an open set. By (3), we see that for $\mu$-almost every $y \in D$, 
	\begin{equation}\label{e:FK3}
		\P^y(\bar \tau_D>t)
		= \int_D \bar  p^D(t,y,z) \mu(dz) \le \int_D \bar p^B(t,y,z) \mu(dz) \le \frac{K_2\mu(D)}{V(x,r)} \Big( \frac{\phi(r)}{t}\Big)^{1/\nu}. 
	\end{equation}
	Set $T: = \big( \mu(D)/V(x,r)\big)^\nu \phi(r)$. It follows from \eqref{e:FK3} that for   $\mu$-almost every $y \in D$, 
	\begin{equation}\label{e:E<}
		\bar \E^y \bar\tau_D \le T+ \int_T^\infty \bar \P^y(\bar \tau_D>t) dt  \le T + K_2 T^{1/\nu}\int_T^\infty t^{-1/\nu}  dt = \big(1 + \nu(1-\nu)^{-1} K_2\big)T.
	\end{equation}
	Thus, by \cite[Lemma 6.2]{GH14}, we obtain that $\lambda_1(D) \ge (\esssup_{y \in D} \bar\E^y \bar \tau_D)^{-1} \ge c_2T^{-1}$. \qed

	\begin{cor}\label{c:E<}
		Suppose that one of the equivalent conditions of Lemma \ref{l:7.3} is satisfied.
		Then there exists a constant $c>0$ such that 
		\begin{equation*}
			\bar\E^{x}[\bar \tau_{B(x,r)}] \le c \phi(r) \quad \text{for all} \;\; x \in U \setminus \bar \sN \;\; \text{and} \;\; 0<r<R \land (\kappa\updelta_U(x)).
		\end{equation*}
	\end{cor}
	\pf Let $x \in U \setminus \bar \sN$ and $0<r<R \land (\kappa\updelta_U(x))$. By taking $D=B(x,r)$ in \eqref{e:E<}, we see that $\esssup_{y \in B(x,r)} \bar\E^y \bar \tau_{B(x,r)} \le c_1\phi(r)$ for  $c_1>0$ independent of $x$ and $r$. Hence, by the strong Markov property, since the heat kernel $\bar p^{B(x,r)}(t, \cdot, \cdot)$ exists due to Lemma \ref{l:7.3}, we obtain
	\begin{align*}
		\bar\E^x[\bar \tau_{B(x,r)}] &\le \phi(r) + \bar\E^x[ \bar\tau_{B(x,r)} - \phi(r); \bar\tau_{B(x,r)} \le \phi(r)] +  \bar\E^x[\bar \tau_{B(x,r)} - \phi(r); \bar\tau_{B(x,r)} >\phi(r)] \\
		&\le  \phi(r) + \int_{B(x,r)} \bar\E^z[\bar \tau_{B(x,r)}] \bar p^{B(x,r)}(\phi(r), x, z) \mu(dz) \le (1+c_1)\phi(r).
	\end{align*}
	\qed
	
	A global version of the following lemma can be founded in \cite[Proposition 7.6]{CKW16a}.

	\begin{lem}\label{l:7.4}
		\noindent Let $U \subset M$ be an open set, $R_0 \in (0,\infty]$ a constant, and $\phi$ an increasing function satisfying $\U(\phi,\gamma,c)$ for some $\gamma,c>0$. If $\PI_{R_0}(\phi,U)$  holds, then  $\FK_{\eps(R_0 \land \bar R)}(\phi,U)$  holds for some $\eps \in (0,1)$. 
	\end{lem}
	\pf  Without loss of generality, we may  assume that $R_0 \le \bar R$. By Lemma \ref{l:7.3}, it suffices to prove that $\PI_{R_0}(\phi,U)$ implies  (2) therein with $R=\eps R_0$ for some $\eps \in (0, 1)$.  Let $C_V \in (0,1)$ be the constant that $\mathrm{VRD}_{\bar R}(M)$ holds with. 
	
	By following the proof of \cite[Theorem 2.1]{Sa92} or \cite[Proposition 2.3]{CG98}, since  $\PI_{R_0}(\phi,U)$ holds, one can see that the following Sobolev-type  inequality holds: There is constants $\nu \in (0,1)$, $\kappa_1 \in (0, C_V)$ and $c_1>0$ such that for all $x \in U$ and $0<r<R_0 \land (\kappa_1 \updelta_U(x))$,
	\begin{equation}\label{e:Nash} \Vert u \Vert_2^{2+2\nu} \le \frac{c_1}{V(x,r)^\nu} \Vert u \Vert_1^{2\nu} \big(  \phi(r) \EE(u,u)+  \Vert u \Vert_2^2  \big),\quad  \;\;  \forall \, u \in \FF_{B(x,r)}. 
	\end{equation}
	Indeed, even though they only proved \eqref{e:Nash} when $\phi(r)=r^2$, with simple modifications, one can easily follow their proof with general $\phi$ satisfying $\U(\phi,\gamma,c)$. 
	
	Now, we adopt a method in the proof of \cite[Proposition 2.3]{CG98}. Choose a constant  $C>1$ such that
	\begin{equation*}
		2c_1 V(x,r)^\nu \le V(x,Cr)^\nu \quad \text{for all} \;\; x \in U, \; 0<r<R_0 \land (\kappa_1 \updelta_U(x)).
	\end{equation*}
	Since ${\rm VRD}_{\bar R}(M)$ holds, $R_0 \le \bar R$ and $\kappa_1< C_V$, such constant $C$ exists. By \eqref{e:Nash} and  Cauchy inequality, we get that for all  $x \in U$, $0<r<C^{-1}\big(R_0 \land (\kappa_1 \updelta_U(x))\big)$ and $u \in \FF_{B(x,r)} \subset \FF_{B(x,Cr)}$,
	\begin{align*}
		\Vert u \Vert_2^{2+2\nu} &\le \frac{c_1}{V(x,Cr)^\nu} \Vert u \Vert_1^{2\nu} \big(\phi(Cr) \EE(u,u) +  \Vert u \Vert_2^2 \big) \\
		&\le \frac{\phi(Cr) }{2V(x,r)^\nu} \Vert u \Vert_1^{2\nu} \EE(u,u) +   \frac{\Vert u \Vert_1^{2\nu}}{2V(x,r)^\nu}  \Vert u \Vert_2^2 \le \frac{\phi(Cr) }{2V(x,r)^\nu} \Vert u \Vert_1^{2\nu} \EE(u,u) +   \frac{1}{2}  \Vert u \Vert_2^{2+2\nu}.
	\end{align*}
	Thus, since $\U(\phi,\gamma,c)$ holds, we conclude that  (2) in Lemma \ref{l:7.3} holds with $R=C^{-1}R_0$.
	\qed

	Now, let us define a L\'evy measure $\nu_0$ and a Bernstein function $\phi_1$ as (cf. \cite[(4.6)]{BKKL19b})
	\begin{equation}\label{e:nuphi0}
		\nu_0(ds):=\frac{ds}{ s\psi_0(F^{-1}(s))}, \;\; \quad  \phi_1(\lambda)  :=\ll \lambda + \int_0^\infty (1-e^{-\lambda s})\nu_0(ds).
	\end{equation}
	Since (L1) and $\U(F, \gamma_2, c_U)$ hold, one can see $\nu_0((0,\infty)) = \infty$ and $\int_0^\infty (1 \land s)\nu_0(ds)<\infty$. Hence, the above    $\phi_1$ is well-defined.

	Let $S_t$ be a subordinator with the  Laplace exponent $\phi_1$, and set $Y_t:=Z_{S_t}$. Define $J_Y(x,y)$ as the right hand side of the second equality in \eqref{e:PJ} with $\nu=\nu_0$.
	Then  $Y_t$ is a symmetric Hunt process associated with the following regular Dirichlet form $(\EE^Y, \FF^Y)$ on $L^2(M;\mu)$:
	\begin{align}\label{e:EY}
		\EE^Y(f,f)&= \ll \EE^{Z}(f,f) + \int_{M \times M \setminus \text{diag}} (f(x)-f(y))^2J_Y(x,y)\mu(dx)\mu(dy), \;\;\; f \in \FF^Y \hspace{-1.5mm}, \; \\
		\FF^Y&=\overline{\big\{f \in C_c(M): \EE^Y(f,f)<\infty\big\}}^{\EE^Y_1}.
	\end{align}
	Moreover, according to \cite{AR05}, we have   $\FF^Z \subset \FF^Y$, and $\FF^Z=\FF^Y$ if $\ll>0$.

	Recall the definition of $\Phi_1$ from \eqref{d:HuntPhi}.

	\begin{lem}\label{l:Phi-comp}
		(i) It holds that $\nu_0((r,\infty)) \asymp \psi_0(F^{-1}(r))^{-1}$ on $(0,\infty)$.
		
		\noindent (ii) $\L(\phi_1, \alpha_1, c)$ holds for some $\alpha_1, c>0$.
		
		\noindent (iii)  It holds that $\Phi_1(r) \asymp \phi_1(F(r)^{-1})^{-1}$ on $(0,\infty)$.
	\end{lem}
	\pf (i) Observe that  $\L(\psi_0 \circ F^{-1}, \gamma_2^{-1} \beta_1,c_1)$ and  $\U(\psi_0 \circ F^{-1}, \gamma_1^{-1} (\beta_2 + 2 \gamma_2),c_2)$ hold for some constants $c_2 \ge 1 \ge c_1>0$. Thus, by \cite[Lemma 2.3(1)]{CK202}, we get the result.
	
	 (ii) Let $\phi_0(\lambda):=\phi_1(\lambda)-\ll \lambda$. Since $\L(\psi_0 \circ F^{-1}, \gamma_2^{-1} \beta_1,c_1)$ holds, by  a similar argument to the one for \cite[Lemma 2.3(3)]{CK202},  we see that  $\L(\phi_0, \alpha_1, c_3)$ holds for some $\alpha_1 \in (0,1)$. Thus, if $\ll=0$, then the result follows. Now, assume that $\ll>0$.  Then we see that  $\phi_1(\lambda) \asymp \phi_0(\lambda)$ on $(0,1]$, and $\phi_1(\lambda)  \asymp \ll \lambda$ on $[1,\infty)$ (see Example \ref{e:mixedtype}). Since $\Lambda \lambda$ is linear, we  deduce the  result from $\L(\phi_0, \alpha_1, c_3)$.

	 (iii) By  \cite[Lemma 2.1(i)]{CK20}, the above (i) and the change of the variables, we obtain
	\begin{equation*}
		\phi_1\big(F(r)^{-1}\big)^{-1} \asymp  F(r) \Big(\int_0^{F(r)}\nu_0((s,\infty))ds + \ll\Big)^{-1} \asymp \Phi_1(r) \quad \text{for}\;\, r>0.
	\end{equation*}
	\qed
	
	\noindent As a consequence of the above  Lemma \ref{l:Phi-comp}(ii, iii), since $\U(\phi_1, 1, 1)$ holds, 
	\begin{equation}\label{e:scalePhi1}
		\text{$\L(\Phi_1, \alpha_1\gamma_1, c_{L}')\;$  and  $\;\U(\Phi_1,  \gamma_2, c_U')\;$ hold with some constants $c_L',c_U'>0$.}
	\end{equation}
	
	\begin{lem}\label{l:JY}  (i)	There exists  $a_1 \ge 1$ such that for all $x,y \in M$ with $d(x,y) < F^{-1}(2^{-1}(F(R_1))$, 
		\begin{equation}\label{e:JY}
			\frac{a_1^{-1}}{V(x,d(x,y))\psi_0(d(x,y))} \le J_Y(x,y) \le \frac{a_1}{V(x,d(x,y))\psi_0(d(x,y))}.
		\end{equation}
		
		\noindent (ii) There exists $a_2 \ge 1$ such that for all $x,y \in \sU$  with $d(x,y) < F^{-1}(2^{-1}F(R_1)) \land R_4 \land (\kappa_0 \updelta_\sU(x))$,
		\begin{equation}\label{e:JY1}
			a_2^{-1}J_X(x,y) \le J_Y(x,y) \le a_2J_X(x,y).
		\end{equation}
	\end{lem}
	\pf (i) Fix any $x,y \in M$ with $d(x,y) < F^{-1}(2^{-1}(F(R_1))$ and let $l:=d(x,y)$.
	
	First, we have that, by  \eqref{e:Fscale}, \eqref{e:diffusion}, \eqref{e:scalepsi0} and \eqref{e:nuphi0}, since ${\rm VRD}_{\bar R}(M)$ holds, 
	\begin{equation*}
		J_Y(x,y) \ge \int_{F(l)}^{2F(l)} q(s,x,y)\nu_0(ds)
		\ge \frac{c_1 }{V(x,l) F(l) \psi_0(l)} \int_{F(l)}^{2F(l)} ds = \frac{c_1}{V(x,l)\psi_0(l)}.
	\end{equation*}

	Next, to prove the upper bound in \eqref{e:JY}, we  claim that there exists a constant $c_2>0$ independent of $x$ such that
	\begin{equation}\label{e:uhkd-diffusion}
		\sup_{z \in M} q(t,x,z) \le \frac{c_2}{V(x,F^{-1}(t) \land R_1)} \quad \text{for all} \;\;  t>0.
	\end{equation}
	Indeed, in view of  \eqref{e:diffusion}, the above \eqref{e:uhkd-diffusion} holds for $t \in (0, F(R_1))$. Hence, it suffices to consider the case $R_1<\infty$. In such case, by the semigroup property and ${\rm VRD}_{\bar R}(M)$, since $Z$ is $\mu$-symmetric, we see that for all $t \ge F(2^{-1}R_1)$ and $z \in M$,
	\begin{align*}
		q(t,x,z) \le \frac{c_3}{V(x,2^{-1}R_1)} \int_M q(t-F(2^{-1}R_1),z,w)\mu(dw)  \le \frac{c_4}{V(x,R_1)}.
	\end{align*}
	Thus,  \eqref{e:uhkd-diffusion} holds.
	
	Now, by \eqref{e:Fscale}, \eqref{e:diffusion}, \eqref{e:scalepsi0}, \eqref{e:nuphi0},  \eqref{e:uhkd-diffusion},  ${\rm VRD}_{\bar R}(M)$   and Lemmas \ref{l:expF}(i) and \ref{l:Phi-comp}(i), we obtain
	\begin{align*}
		&J_Y(x,y) = \Big(\int_0^{F(l)} + \int_{F(l)}^{F(R_1)} + \int_{F(R_1)}^{\infty} \Big) q(s,x,y)\nu_0(ds) \\
		&\le \int_0^{F(l)} \frac{c_5\exp\big(-c_6 F_1(l,s) \big)}{s  V(x,F^{-1}(s))\psi_0(F^{-1}(s))} ds  + \int_{F(l)}^{F(R_1)} \frac{c_7 ds }{s V(x,F^{-1}(s)) \psi_0(F^{-1}(s))}  + \frac{c_8}{V(x,R_1)\psi_0(R_1)}
		\\
		&\le \int_0^{F(l)} \frac{c_9}{F(l) V(x,F^{-1}(s)) \psi_0(l)} \Big( \frac{s}{F(l)} \Big)^{d_2 / \gamma_1} ds +  \frac{c_{10}}{V(x,l) \psi_0(l)}  + \frac{c_8}{V(x,R_1)\psi_0(R_1)} \le \frac{c_{11}}{V(x,l)\psi_0(l)}.
	\end{align*}
	We used a  convention that $V(x, \infty) = \psi_0(\infty)=\infty$ in the above. 
	
 (ii) The result follows from \eqref{e:JY} and (L2). \qed

	\begin{prop}\label{p:EY}
		There exists a constant $\eps_0 \in (0,1)$ such that  $\CS_{\eps_0 R_1}(\Phi_1,M)$, $\FK_{\eps_0 R_1}(\Phi_1,M)$ and $\PI_{\eps_0 R_1}(\Phi_1,M)$ for $(\EE^Y,\FF^Y)$ hold.
	\end{prop}
	\pf By Lemmas \ref{l:subTail}(i), \ref{l:subNDL}(i,ii) and  \ref{l:Phi-comp}(ii,iii), we see that $\Tail_{R_1}(\Phi_1,\le,M)$ and   $\NDL_{c_1 R_1}(\Phi_1,M)$ for $Y$ hold for some $c_1 \in (0,1)$. Hence,  by Proposition \ref{p:E}(i) and \eqref{e:scalePhi1},  $Y$ satisfies  Assumption \ref{ass0} with $U=M$. For an open set $D \subset M$, denote by $\tau^Y_D$ the first exit time of $Y$ from $D$, and $(P^{Y,D}_t)_{t \ge 0}$ the semigroup of the subprocess $Y^{D}$ in $L^2(M;\mu)$.
	By Proposition \ref{p:EP}(i) and \eqref{e:scalePhi1}, it follows  that there exist   $c_2 \in (0,1]$ and $c_3>0$ such that for all  $z \in M$ and $r \in (0, c_2 R_1)$, 
	\begin{equation}\label{e:EPY1}
		\underset{y \in B(z, \frac{1}{4}r)}{\essinf}\,\, P_t^{Y,B(z,r)} \1(y) \ge \underset{y \in B(z,\frac{1}{4}r)}{\essinf} \,\, \P^y ( \tau^Y_{B(y,\frac{3}{4}r)} > t  )  \ge 1/2, \quad \;\; \forall \, 0< t  \le c_3 \Phi_1(r).
	\end{equation}
	By \eqref{e:EPY1}, since $\Tail_{R_1}(\Phi_1,\le,M)$ for $Y$ holds, one  can follow the proofs of \cite[Lemma 2.8]{GHH18} and \cite[Proposition 2.5]{CKW19} in turn to deduce that  $\CS_{c_4R_1}(\Phi_1,M)$ for $(\EE^Y,\FF^Y)$ holds for some $c_4 \in (0, c_2)$.
	
	Next, since $\NDL_{c_1R_1}(\Phi_1,M)$ for $Y$ holds, by following the proof of \cite[Proposition 3.5(i)]{CKW16b}, we can deduce that  $\PI_{c_5R_1}(\Phi_1,M)$ for $(\EE^Y,\FF^Y)$ holds for some $c_5 \in (0,1)$.  Indeed, even though only pure-jump type Dirichlet forms are considered in \cite{CKW16b}, after  redefining the form $\bar \EE$ given in the fifth line of the proof of \cite[Proposition 3.5(i)]{CKW16b} by 
	$$ \bar\EE(u,v) = \ll \int_{B(x_0,r)} d\Gamma^Z(u,v) + \int_{B(x_0,r) \times B(x_0,r)} (u(x) - u(y))(v(x) - v(y)) J_Y(x,y)\mu(dx)\mu(dy),$$ one can follow the rest of the proof.
	
	Lastly, since $R_1 \le \bar R$, by Lemma \ref{l:7.4},  $\FK_{c_6 R_1}(\Phi_1,M)$ for $(\EE^Y,\FF^Y)$ holds with some $c_6 \in (0,c_5)$.
	In the end, we finish the proof by taking $\eps_0= c_4 \land c_6$.
	\qed
	
	\begin{lem}\label{l:Domain}
		It holds that	$\FF^X_\sU = \FF^Y_\sU$.
	\end{lem}
	\pf By (L3), (L4) and \eqref{e:JY1}, since $\Tail_{R_1}(\Phi_1,  \le,M)$ for $Y$ holds due to Lemmas \ref{l:subTail}(i) and \ref{l:Phi-comp}(iii), one sees that for every $f \in C_c(\, \sU)$, $\EE^X(f,f)<\infty$ if and only if $\EE^Y(f,f)<\infty$. Therefore, we have $\FF^X \cap C_c(\, \sU) = \FF^Y \cap C_c(\, \sU)$.

	Let $f \in \FF^X_\sU$. Then there exists a sequence $(f_n)_{n \ge 1} \subset \FF^X\cap C_c(\,\sU)$ such that $\EE^X_1(f-f_n,f-f_n) \to 0$. Fix a compact set $S \subset \sU$ such that $\text{supp}(f),\text{supp}(f_n) \subset S$ and let $\delta:=\inf_{y \in S, z \in M \setminus\sU}d(y,z)>0$. Then by \eqref{e:EY}, (L4) and  \eqref{e:JY1}, since $\Tail_{R_1}(\Phi_1, \le,M)$ for $Y$ holds,  we have that 
	\begin{align*}
		&\EE^Y(f-f_n,f-f_n) \\
		&\le c_1\EE^X(f-f_n,f-f_n) +\int_{\sU} \int_{M \setminus B(x, c_2 \delta)} \big( (f(x)-f_n(x))-(f(y)-f_n(y))\big)^2 J_Y(x,y) \mu(dx)\mu(dy)\\
		& \quad +  \int_{M \setminus \sU} \int_{S} \big( (f(x)-f_n(x))-(f(y)-f_n(y))\big)^2 J_Y(x,y) \mu(dx)\mu(dy) \\
		& \le c_1\EE^X(f-f_n, f-f_n) + \frac{c_3}{\Phi_1((c_2 \land 1)\delta)} \lVert f-f_n \rVert_{L^2(M;\mu)} \to 0 \quad \text{as} \;\; n \to \infty.
	\end{align*}
	Hence, $f \in \FF^Y_\sU$. Conversely, by the same argument, one can see that $\FF^Y_\sU\subset \FF^X_\sU$. \qed

	\begin{prop}\label{p:EX}
		There exists a constant $R_0 \in (0,\infty]$ such that  $\CS_{R_0}(\Phi_1,\sU)$, $\FK_{R_0}(\Phi_1,\sU)$ and $\PI_{R_0}(\Phi_1,\sU)$ for $(\EE^X,\FF^X_\sU)$ hold. In particular, if $\sU = M$ and $R_1 = R_4 = \infty$, then $\CS_{\infty}(\Phi_1,M)$, $\FK_{\infty}(\Phi_1,M)$ and $\PI_{\infty}(\Phi_1,M)$ for $(\EE^X,\FF^X)$ hold.
	\end{prop} 
	\pf  By Proposition \ref{p:EY}, ${\rm CS}_{\eps_0R_1}(\Phi_1,M)$, ${\rm FK}_{\eps_0R_1}(\Phi_1,M)$ and $\PI_{\eps_0 R_1}(\Phi_1,M)$ for $(\EE^Y,\FF^Y)$ hold with a  constant  $\kappa \in (0,1)$ for some $\eps_0 \in (0,1)$. Let $\kappa_0\in (0,1)$ be the constant in Assumption L.  Let $C\ge2$ be a constant  chosen to be later. We  set $\wt \kappa:= \kappa_0 \kappa/(3C)$ and  $R_5:=\eps_0 \, \wt \kappa \, (F^{-1}(2^{-1}F( R_1)) \land R_4)$.  Choose any $x \in \sU$ and $r \in (0, R_5 \land (\wt\kappa\updelta_\sU(x)))$.

	For all $y, z \in B(x, Cr/\kappa)$, we have $d(y,z) <2Cr/\kappa< 2CR_5/\kappa< F^{-1}(2^{-1}F(R_1)) \land R_4$ and 
	\begin{equation*}
		\frac{d(y,z)}{\kappa_0}< \frac{2Cr}{\kappa_0\kappa}< \frac{2 C\wt \kappa}{\kappa_0\kappa} \updelta_\sU(x) = \frac{2}{3}\updelta_\sU(x) \le \updelta_\sU(x)- \frac{C\wt \kappa}{\kappa}\updelta_\sU(x)< \updelta_\sU(x)- \frac{Cr}{\kappa} <  \updelta_\sU(y).
	\end{equation*}
	Hence, by (L4) and \eqref{e:JY1}, there exists $c_1>1$ independent of $x$ and $r$ such that for any $f \in \FF^X_\sU$,
	\begin{align}\label{e:comXY}
		&\ll \int_{B(x, Cr/ \kappa)} d\Gamma^Z(f,f) + \int_{B(x, Cr/ \kappa) \times B(x, Cr/ \kappa)} (f(y) - f(z))^2 J_Y(y,z)\mu(dy)\mu(dz) \nn\\
		&\le c_1 \bigg(\ll \int_{B(x,r/\wt \kappa )} d\Gamma^{X,(c)}(f,f) + \int_{B(x, r/\wt \kappa ) \times B(x, r/\wt \kappa )} (f(y) - f(z))^2 J_X(y,z)\mu(dy)\mu(dz) \bigg). 
	\end{align}
	Thus,   we deduce  $\PI_{R_5}(\Phi_1,\sU)$ for $(\EE^X,\FF^X_\sU)$ from $\PI_{\eps_0 R_1}(\Phi_1,M)$ for $(\EE^Y,\FF^Y)$ and Lemma \ref{l:Domain}.
	
	Next, by \eqref{e:EY}, \eqref{e:comXY}, (L4),  Lemmas \ref{l:subTail}(i), \ref{l:Phi-comp}(iii) and \ref{l:Domain},  \eqref{e:scalePhi1} and $\mathrm{FK}_{\eps_0 R_1}(\Phi_1, M)$ for $(\EE^Y,\FF^Y)$, we have  that,  for any open set $D \subset B(x,r)$ and  $f \in \FF^X_\sU$ with $\Vert f\Vert_2 = 1$,
	\begin{align}\label{e:XFK}
		\EE^X(f,f)&\ge c_1^{-1} \EE^Y(f,f) - 4\Vert f \Vert_2^2 \sup_{y \in B(x,r)} \int_{B(x,Cr/\kappa)^c} J_Y(y,z)\mu(dz) \nn\\
		& \ge c_1^{-1} \EE^Y(f,f) - \frac{c_2}{\Phi_1((C/\kappa -1)r)} \ge \frac{c_1^{-1}c_3}{\Phi_1(r)} \left( \frac{V(x,r)}{\mu(D)}\right)^\nu - \frac{c_2c_L^{-1}}{(C/\kappa-1)^{\alpha_1\gamma_1}\Phi_1(r)}\nn\\[3pt]
		& \ge \frac{c_1^{-1}c_3 - c_2 c_L'^{-1}(C/\kappa-1)^{-\alpha_1\gamma_1}}{\Phi_1(r)} \left( \frac{V(x,r)}{\mu(D)}\right)^\nu.
	\end{align}
	The last inequality above is valid since $V(x,r)/\mu(D) \ge 1$. Here, we point out that the constants $c_1,c_2$ and $c_3$ are independent of $x$ and $r$. Now, we choose $C=1 + (2c_1c_2/(c_3c_L'))^{1/(\alpha_1\gamma_1)}$. Then  we get from \eqref{e:XFK} that  $\FK_{R_5}(\Phi_1,\sU)$ for $(\EE^X,\FF^X_\sU)$ holds.

	Now, we show that  $\CS_{R_5}(\Phi_1,\sU)$ for $(\EE^X,\FF^X_\sU)$ holds.  Observe that for all $x_0 \in \sU$, $f \in \FF^X_\sU$, $0<u \le s <R_5 \land (\wt \kappa\updelta_\sU(x_0))$ and $a \in (0,1]$,   and any cut-off function $\varphi$ for $B(x_0,s) \subset B(x_0,s+u)$, by (L3), (L4), \eqref{e:JY1} and Lemma \ref{l:JX}, it holds that
	\begin{align}\label{e:CSX}
		&\int\limits_{B(x_0,s+(1+a)u)} \hskip -0.1in f^2 d\Gamma^X(\varphi,\varphi) \nn\\
		& = \int\limits_{B(x_0,s+(1+a)u)} \hskip -0.1in f^2 d\Gamma^{X,(c)}(\varphi,\varphi)+ \int\limits_{B(x_0,s+(1+a)u) \times M} \hskip -0.1in  f(y)^2 (\varphi(y)- \varphi(z))^2 J_X(y,z)\mu(dy)\mu(dz) \nn\\
		&\le c_4 \int\limits_{B(x_0,s+(1+a)u)} \hskip -0.1in f^2 d\Gamma^Y(\varphi,\varphi) +  \int\limits_{B(x_0,s+(1+a)u)} \hskip -0.1in f(y)^2 J_X\Big(y,M \setminus B\big(y, (\eps_0^{-1}\wt \kappa^{-1} R_5) \land (\kappa_0\updelta_\sU(y))\big) \Big) \mu(dy) 
		\nn\\ 
		&\le c_4 \int\limits_{B(x_0,s+(1+a)u)} \hskip -0.1in f^2 d\Gamma^Y(\varphi,\varphi) +  \frac{c_5}{\psi_0(u)} \int\limits_{B(x_0,s+(1+a)u)} \hskip -0.1in f^2 d\mu.
	\end{align}
	Therefore, since $\CS_{\eps_0R_1}(\Phi_1, M)$ for $(\EE^Y,\FF^Y)$ hold and $\Phi_1(u) \le \psi_0(u)$ for all $u \in (0, R_4)$, we can deduce from Lemma \ref{l:Domain}, \eqref{e:CSX}, (L4) and \eqref{e:JY1} that $\CS_{R_5}(\Phi_1, \sU)$ for $(\EE^X, \FF^X_\sU)$ holds.

	Lastly, since $R_5=\infty$ if $R_1=R_4=\infty$, and $\FF^X_\sU = \FF^X$ if $\sU = M$, the latter assertion  holds. \qed

	\begin{prop}\label{p:NDLX}
		There exists a constant $R \in (0, \infty]$ such that $\mathrm{NDL}_{R}(\Phi_1, \sU\cap M_0)$ for $X$ holds. Moreover, if $\sU=M$ and $R_1=R_4=\infty$, then $\mathrm{NDL}_\infty(\Phi_1,M_0)$ for $X$ holds.
	\end{prop}
	\pf  According to Lemma \ref{l:JX} and Proposition \ref{p:EX},  since $\Phi_1(r) \le \psi_0(r)$ for $r \in (0, R_4)$, we see that  $\Tail_{R_0}(\Phi_1,\le,\sU)$, $\FK_{R_0}(\Phi_1,\sU)$ and $\CS_{R_0}(\Phi_1,\sU)$  hold with $\kappa$ for some $R_0 \in (0, \infty]$ and $\kappa\in(0,1)$. Hence, by  \eqref{e:scalePhi1}, one can  follow the proofs of \cite[Lemmas 4.15 and 4.17]{CKW16a} to see that  there are constants $R_6 \in (0, R_0]$ and $c_1,c_2\in (0,1)$ such that for any $x \in \sU \cap M_0$ and $r \in (0,R_6 \land (\kappa \updelta_\sU(x)))$,
	\begin{equation}\label{e:2.8.1}
		\P^x(\tau_{B(x,r)} \ge \Phi_1(2c_1r)) \ge c_2.
	\end{equation}

	Fix any $x_0 \in \sU \cap M_0$, $0<r<R_6 \land (\kappa\updelta_\sU(x_0))$ and let $B = B(x_0,r)$. 
	By Lemma \ref{l:7.3}, the heat kernel $p^B(t,x,y)$ of $X^B=(X^\sU)^B$ exists. Moreover, since $X^B$ is symmetric, by \eqref{e:2.8.1} and Cauchy inequality, we have that for all $x \in B(x_0, r/2) \setminus \sN$ and  $t \le \Phi_1(c_1r)$,
	\begin{align}\label{e:2.8.2}
		p^B(t,x,x) &= \int_{B} p^B(t/2,x,y)^2 \mu(dy) \ge \frac{1}{\mu(B)} \left(\int_B p^B(t/2,x,y)\mu(dy) \right)^2 \nn \\ &
		\ge \frac{1}{V(x_0,r)} \P^{x}\big(\tau_{B(x,r/2)} \ge \Phi_1(c_1r)\big)^2\ge \frac{c_2^2}{V(x_0,r)}.  
	\end{align}

	We claim that after assuming that $R_6$ and $\kappa$ are sufficiently small, there are constants $c_3,\theta>0$ and $\eta_0 \in (0,1)$ independent of $x_0$ and $r$  such that for all $t \le \Phi_1(c_1r)$ and $x,y \in B(x_0, \eta_0\Phi_1^{-1}(t)) \setminus \sN$,
	\begin{equation}\label{e:EHR}
		\big|\, p^{B(x_0, c_1^{-1}\Phi_1^{-1}(t))}(t,x,x) - p^{B(x_0, c_1^{-1}\Phi_1^{-1}(t))}(t,x,y)\big| \le \frac{c_3}{V(x_0,\Phi_1^{-1}(t))}\left( \frac{d(x,y)}{\Phi_1^{-1}(t)} \right)^\theta.
	\end{equation}
	
	To prove \eqref{e:EHR}, we first obtain a local version of elliptic Hölder regularity (EHR) of harmonic functions. First, since $\CS_{R_0}(\Phi_1,\sU)$ holds, we can follow the proofs of  \cite[Proposition 2.9]{CKW19} and \cite[Proposition 4.12]{CKW16b} to get local versions of those two proposition, by applying $\Tail_{R_0}(\Phi_1,\le,\sU)$ in the third inequality in the display of the proof for \cite[Proposition 2.9]{CKW19},  and the last inequality in \cite[p.3787]{CKW16b}. By using those two results, one can see that a local version of  \cite[Proposition 5.1]{CKW19} holds by following its proof line by line, after redefining the notation $\mathrm{Tail}_\phi(u;x_0,r)$ therein by $\phi(r) \int_{M_\partial \setminus B(x_0,r)} |u(z)|J_X(x_0, dz)$. 
	Besides, since $\FK_{R_0}(\Phi_1, \sU)$ and $\CS_{R_0}(\Phi_1,\sU)$ hold, by following the proofs of  \cite[Lemmas 4.6, 4.8 and 4.10]{CKW16a} and using  $\Tail_{R_0}(\Phi_1, \le,\sU)$ whenever  $\mathrm{J}_{\phi,\le}$ therein used,  one can obtain their local versions hold (so-called Caccioppoli inequality, comparison inequality over balls and $L^2$-mean value inequality for subharmonic functions). In the end, by using local versions of \cite[Proposition 5.1]{CKW19} and \cite[Lemmas 4.8, 4.10]{CKW16a}, since $\Tail_{R_0}(\Phi_1,\le,\sU)$, $\CS_{R_0}(\Phi_1, \sU)$ and $\PI_{R_0}(\Phi_1, \sU)$ hold, one can follow the proofs of \cite[Corollary 4.13 and Proposition 4.14]{CKW16b} in turn to obtain a local EHR. Now, by using the local EHR, Corollary \ref{c:E<},  \eqref{e:Dup} and \eqref{e:scalePhi1}, since $\mathrm{VRD}_{\bar R}(M)$ holds, we can follow the proof of \cite[Lemma 4.8]{CKW16b} and deduce that  \eqref{e:EHR} holds.
	
	Eventually, by choosing $\wt \eta \in (0, c_1\eta_0/2)$ sufficiently small, since $\mathrm{VRD}_{\bar R}(M)$ holds, we get from  \eqref{e:2.8.2} and \eqref{e:EHR} that for all $y,z \in B(x_0, \wt \eta^2 r) \setminus \sN$,
	\begin{align*}
		p^{B}(\Phi_1(\wt\eta r), y,z) &\ge p^{B(x_0, c_1^{-1} \wt \eta r)}(\Phi_1(\wt\eta r), y,z)\ge  p^{B(x_0, c_1^{-1} \wt \eta r)}(\Phi_1(\wt\eta r), y,y) -   \frac{c_3}{V(x_0,\wt \eta r)}\left( \frac{2\wt \eta^2 r}{\wt \eta r} \right)^\theta  \\
		& \ge \frac{1}{V(x_0, c_1^{-1}\wt \eta r)}\left(c_2^2 - \frac{2^\theta \wt \eta^\theta c_3V(x_0,c_1^{-1}\wt \eta r)}{V(x_0, \wt \eta r)}\right) \ge \frac{2^{-1}c_2^2}{V(x_0, r)}.
	\end{align*}
	This shows that $\mathrm{NDL}_{R_6}(\Phi_1, \sU \cap M_0)$ for $X$ holds. Then in view of the latter assertion in Proposition \ref{p:EX}, by the same proof, we can see that the second claim in the proposition also holds.
	\qed

	Finally, the proof for Theorem \ref{t:Hunt0} is straightforward.

	\medskip

	\noindent \textbf{Proof of Theorem \ref{t:Hunt0}.} 
	By \eqref{e:scalepsi0}, \eqref{e:scalePhi1}, Lemma \ref{l:JX}, Propositions \ref{p:NDLX} and \ref{p:E}, the theorem follows from 
	Theorems  \ref{t:limsup0-1}--\ref{t:supprecise2}  in Section \ref{s:intro}. \qed


\bigskip
	
	\small
	\begin{thebibliography}{10}
		
		\bibitem{AR05}
		S. Albeverio, B. Rüdiger.
		\newblock Subordination of symmetric quasi-regular Dirichlet forms.
		\newblock{\em 	Random Oper. Stochastic Equations} 13 (2005), no. 1, 17–38.
		
				\bibitem{BKKL19a}
		J.~Bae, J.~Kang, P.~Kim,  J.~Lee.
		\newblock Heat kernel estimates for symmetric jump processes with mixed
		polynomial growths.
		\newblock {\em Ann. Probab.} 47 (2019), no. 5, 2830–2868.
		
		\bibitem{BKKL19b}
		J.~Bae, J.~Kang, P.~Kim,  J.~Lee.
		\newblock Heat kernel estimates and their stabilities for symmetric jump
		processes with general mixed polynomial growths on metric measure spaces.
		\newblock available at arXiv:1904.10189v3.
		

		
		\bibitem{Ba98}
		M.~T. Barlow.
		\newblock Diffusions on fractals.
		\newblock In {\em Lectures on probability theory and statistics
			({S}aint-{F}lour, 1995)}, volume 1690 of {\em Lecture Notes in Math.}, pages
		1--121. Springer, Berlin, 1998.
		
		
		\bibitem{Ba04}
		M.~T. Barlow.
		\newblock Random walks on supercritical percolation clusters.
		\newblock {\em Ann. Probab.} 32 (2004), no. 4, 3024–3084.
		
		
		
		
		\bibitem{BB99}
		M.~T. Barlow, R.~F. Bass.
		\newblock Brownian motion and harmonic analysis on {S}ierpinski carpets.
		\newblock {\em Canad. J. Math.}  51 (1999), no. 4, 673–744.
		
		\bibitem{BB99graph} M.~T. Barlow, R.~F. Bass.
		\newblock 	Random walks on graphical Sierpinski carpets. 
		\newblock  Random walks and discrete potential theory (Cortona, 1997), 26–55,
		Sympos. Math., XXXIX, Cambridge Univ. Press, Cambridge, 1999.
		
		\bibitem{BGK09}
		M.~T. Barlow, A.~Grigor'yan,  T.~Kumagai.
		\newblock Heat kernel upper bounds for jump processes and the first exit time.
		\newblock {\em J. Reine Angew. Math.} 626 (2009), 135–157.
		
		\bibitem{BP88}
		M.~T. Barlow, E.~A. Perkins.
		\newblock Brownian motion on the {S}ierpi\'{n}ski gasket.
		\newblock {\em Probab. Theory Related Fields} 79 (1988), no. 4, 543–623.
		
		\bibitem{BK}
		R.~F. Bass, T.~Kumagai.
		\newblock Laws of the iterated logarithm for some symmetric diffusion	processes.
		\newblock {\em Osaka J. Math.} 37 (2000), no. 3, 625–650.
		
		\bibitem{BJ73}
		A.~Benveniste, J.~Jacod.
		\newblock Syst\`emes de {L}\'{e}vy des processus de {M}arkov.
		\newblock {\em Invent. Math.} 21 (1973), 183–198.
		
		
		\bibitem{BGT89}
		N.~H. Bingham, C.~M. Goldie, J.~L. Teugels.
		\newblock  Regular variation.
		\newblock Cambridge University Press, Cambridge, 1987.
		
		
		
		\bibitem{Bo84}
		N. Bouleau. 
		\newblock Quelques résultats probabilistes sur la subordination au sens de Bochner.
		\newblock {\em Seminar on potential theory, Paris, No. 7},  54–81, Lecture Notes in Math., 1061, Springer, Berlin, 1984.
		
		
		\bibitem{CKS87}
		E.~A. Carlen, S.~Kusuoka, D.~W. Stroock.
		\newblock Upper bounds for symmetric {M}arkov transition functions.
		\newblock {\em Ann. Inst. H. Poincaré Probab. Statist.} 23 (1987), no. 2, 245–287.
		
		
		
		\bibitem{CKW18}
		X.~Chen, T.~Kumagai, J.~Wang.
		\newblock Random conductance models with stable-like jumps: quenched invariance principle.
		\newblock {\em Ann. Appl. Probab.}  31 (2021), no. 3, 1180–1231.
		
		
		
		\bibitem{CKW20-1}
		X.~Chen, T.~Kumagai, J.~Wang.
		\newblock Random conductance models with stable-like jumps: heat kernel estimates and Harnack inequalities. 
		\newblock {\em J. Funct. Anal.}  279 (2020), no. 7, 108656, 51 pp.
		
		
		
		\bibitem{CKW16a}
		Z.-Q.~Chen, T.~Kumagai, J.~Wang.
		\newblock Stability of heat kernel estimates for symmetric non-local
		{D}irichlet forms.
		\newblock {\em Mem. Amer. Math. Soc.}  271 (2021), no. 1330, v+89 pp.
		
		\bibitem{CKW19}
		Z.-Q.~Chen, T.~Kumagai, J.~Wang.
		\newblock Heat kernel estimates and parabolic {H}arnack inequalities for
		symmetric {D}irichlet forms.
		\newblock {\em Adv. Math.}  374 (2020), 107269, 71 pp.
		
		\bibitem{CKW20}
		Z.-Q. Chen, T.~Kumagai, J.~Wang.
		\newblock Heat kernel estimates for general symmetric pure jump Dirichlet
		forms.
		\newblock  
		{\em Ann. Sc. Norm. Super. Pisa Cl. Sci. (5)}  23 (2022), 1091–1140

		
		\bibitem{CKW16b}
		Z.-Q. Chen, T.~Kumagai, J.~Wang.
		\newblock Stability of parabolic Harnack inequalities for symmetric non-local
		Dirichlet forms.
		\newblock {\em J. Eur. Math. Soc.}  22 (2020), no. 11, 3747–3803.
		
		
	
		
		
		
		\bibitem{CK20}
		S.~Cho, P.~Kim.
		\newblock Estimates on the tail probabilities of subordinators and applications
		to general time fractional equations.
		\newblock {\em Stochastic Process. Appl.}  130 (2020), no. 7, 4392–4443.
		
		\bibitem{CK202}
		S.~Cho, P.~Kim.
		\newblock Estimates on transition densities of subordinators with jumping
		density decaying in mixed polynomial orders.
		\newblock {\em Stochastic Process. Appl.}  139 (2021), 229–279.
		
		\bibitem{CKL}
		S.~Cho, P.~Kim, J. Lee.
		\newblock 		General Law of iterated logarithm for Markov processes: Liminf laws.
		To appear in  {\em Trans. Amer. Math. Soc.} available at arXiv:2206.08850v2.
		
		
		
	
		\bibitem{CG98}
		T.~Coulhon, A.~Grigoryan.
		\newblock Random walks on graphs with regular volume growth.
		\newblock {\em Geom. Funct. Anal.}  8 (1998), no. 4, 656–701.
		
		
		\bibitem{Dum08}
		H.~Duminil-Copin.
		\newblock Law of the Iterated Logarithm for the random walk on the infinite percolation cluster.
		\newblock {Master thesis}, 
		\newblock {available at arXiv:0809.4380.}
		
		
		
	\bibitem{EL05}		U.~Einmahl, D.~Li.	\newblock Some results on two-sided {LIL} behavior.	\newblock {\em Ann. Probab.}  33 (2005), no. 4, 1601–1624.
		
		
		
		\bibitem{Fe45}	W.~Feller.
		\newblock The fundamental limit theorems in probability.
		\newblock {\em Bull. Amer. Math. Soc.}  51 (1945), 800–832.
		
		\bibitem{Fe68}
		W.~Feller.
		\newblock An extension of the law of the iterated logarithm to variables
		without variance.
		\newblock {\em J. Math. Mech.}  18 1968/1969 343–355.
		
		\bibitem{Fr67}
		B.~Fristedt.
		\newblock Sample function behavior of increasing processes with stationary,
		independent increments.
		\newblock {\em Pacific J. Math.}  21 (1967), 21–33.
		
		\bibitem{Fr71}
		B.~Fristedt.
		\newblock Upper functions for symmetric processes with stationary, independent
		increments.
		\newblock {\em Indiana Univ. Math. J.}  21 (1971/72), 177–185.
		
		
	
		\bibitem{FOT11}
		M.~Fukushima, Y.~Oshima,  M.~Takeda.
		\newblock {\em Dirichlet forms and symmetric {M}arkov processes}, volume~19 of
		{\em De Gruyter Studies in Mathematics}.
		\newblock Walter de Gruyter \& Co., Berlin, Second revised and extended edition, 2011.
		
		\bibitem{Gn43}
		B.~Gnedenko.
		\newblock Sur la croissance des processus stochastiques homog\`enes \`a
		accroissements ind\'{e}pendants.
		\newblock {\em Bull. Acad. Sci. URSS. Sér. Math.}  7, (1943). 89–110.
		
		\bibitem{Gri85}
		P.~S. Griffin.
		\newblock Laws of the iterated logarithm for symmetric stable processes.
		\newblock {\em Z. Wahrsch. Verw. Gebiete}  68 (1985), no. 3, 271–285.
		
		\bibitem{GHH18}
		A.~Grigor'yan, E.~Hu,  J.~Hu.
		\newblock Two-sided estimates of heat kernels of jump type {D}irichlet forms.
		\newblock {\em Adv. Math.}  330 (2018), 433–515.
		
		\bibitem{GH14}
		A.~Grigor'yan, J.~Hu.
		\newblock Upper bounds of heat kernels on doubling spaces.
		\newblock {\em Mosc. Math. J.}  14 (2014), no. 3, 505–563, 641–642.
		
	
		
		\bibitem{GT02}
		A. Grigor'yan, A. Telcs.
		\newblock 	Harnack inequalities and sub-Gaussian estimates for random walks.
		\newblock{\em 		Math. Ann.} 324 (2002), no. 3, 521–556.
		
		
		\bibitem{GT12}
		A.~Grigor'yan, A.~Telcs.
		\newblock Two-sided estimates of heat kernels on metric measure spaces.
		\newblock {\em Ann. Probab.}  40 (2012), no. 3, 1212–1284.
		
		
		
		\bibitem{HK99}
		B.~M. Hambly,  T.~Kumagai.
		\newblock Transition density estimates for diffusion processes on post
		critically finite self-similar fractals.
		\newblock {\em Proc. London Math. Soc.}  (3) 78 (1999), no. 2, 431–458.
		
		\bibitem{HW41}
		P.~Hartman, A.~Wintner.
		\newblock On the law of the iterated logarithm.
		\newblock {\em Amer. J. Math.}  63 (1941), 169–176.
		
		\bibitem{He69}
		C.~C. Heyde.
		\newblock A note concerning behaviour of iterated logarithm type.
		\newblock {\em Proc. Amer. Math. Soc.}  23 (1969), 85–90.
		
		
		\bibitem{JP87}
		N.~C. Jain, W.~E. Pruitt.
		\newblock Lower tail probability estimates for subordinators and nondecreasing
		random walks.
		\newblock {\em Ann. Probab.}  15 (1987), no. 1, 75–101.
		
		\bibitem{Ke72}
		H.~Kesten.
		\newblock Sums of independent random variables without moment conditions.
		\newblock {\em Ann. Math. Statist.}  43 (1972), 701–732.
		
		
		
		
		
		\bibitem{Kh33}
		A.~Khintchine.
		\newblock {\em Asymptotische Gesetze der Wahrscheinlichkeitsrechnung}.
		\newblock Springer, Berlin, 1933.
		
		\bibitem{Kh38}
		A.~Khintchine.
		\newblock {Zwei S\"atze \"uber stochastische Prozesse mit stabilen
			Verteilungen}.
		\newblock {\em Mat. Sb.} 3 (1938), 577-584.
		
		\bibitem{KKW17}
		P.~Kim, T.~Kumagai,  J.~Wang.
		\newblock Laws of the iterated logarithm for symmetric jump processes.
		\newblock {\em Bernoulli}  23 (2017), no. 4A, 2330–2379.
		
		\bibitem{KS06}
		P.~Kim, R.~Song.
		\newblock Two-sided estimates on the density of {B}rownian motion with singular	drift.
		\newblock {\em Illinois J. Math.}  50 (2006), no. 1-4, 635–688.
		
		\bibitem{KSV20}
		P.~Kim, R.~Song,  Z.~Vondra\v{c}ek.
		\newblock On the boundary theory of subordinate killed {L}\'{e}vy processes.
		\newblock {\em Potential Anal.} 53 (2020), no. 1, 131–181.
		
		
		
		\bibitem{KS}
		V.~Knopova, R.~L. Schilling.
		\newblock On the small-time behaviour of {L}\'{e}vy-type processes.
		\newblock {\em Stochastic Process. Appl.} 124 (2014), no. 6, 2249–2265.
		
		\bibitem{Ko29}
		A.~Kolmogoroff.
		\newblock \"{U}ber das {G}esetz des iterierten {L}ogarithmus.
		\newblock {\em Math. Ann.}  101 (1929), no. 1, 126–135.
		
		\bibitem{CK16}
		T. Kumagai, C. Nakamura.
		\newblock Laws of the iterated logarithm for random walks on random conductance models.
		\newblock {\em Stochastic analysis on large scale interacting systems, RIMS Kôkyûroku Bessatsu, B59, Res. Inst. Math. Sci. (RIMS), Kyoto}, 141--156, 2016.
		
		
		\bibitem{Levy2}
		P.~L\'{e}vy.
		\newblock Propri\'{e}t\'{e}s asymptotiques de la courbe du mouvement brownien
		\`a {$N$} dimensions.
		\newblock {\em C. R. Acad. Sci. Paris}  241 (1955), 689–690.
		
		\bibitem{Me75}
		P.~A. Meyer.
		\newblock Renaissance, recollements, m\'{e}langes, ralentissement de processus
		de {M}arkov.
		\newblock {\em Ann. Inst. Fourier (Grenoble)}  25 (1975), no. 3-4, xxiii, 465–497.
		
		\bibitem{Mi16}
		A.~Mimica.
		\newblock Heat kernel estimates for subordinate {B}rownian motions.
		\newblock {\em Proc. Lond. Math. Soc.}  (3) 113 (2016), no. 5, 627–648.
		
	
		\bibitem{Ch17}
		C. Nakamura.
		\newblock Rate functions for random walks on random conductance models and related topics. 
		\newblock {\em Kodai Math. J.}  40 (2017), no. 2, 289–321.
		
		\bibitem{Pr81b}
		W.~E. Pruitt.
		\newblock General one-sided laws of the iterated logarithm.
		\newblock {\em Ann. Probab.}  9 (1981), no. 1, 1–48.
		
		
		\bibitem{Ro68}
		B.~A. Rogozin.
		\newblock On the question of the existence of exact upper sequences.
		\newblock {\em Teor. Verojatnost. i Primenen.}  13 (1968) 701–707.
		
		\bibitem{Sa92}
		L.~Saloff-Coste.
		\newblock A note on {P}oincar\'{e}, {S}obolev, and {H}arnack inequalities.
		\newblock {\em Internat. Math. Res. Notices}  1992, no. 2, 27–38.
		
		\bibitem{Sa13}
		K.-I. Sato.
		\newblock {\em L\'{e}vy processes and infinitely divisible distributions},
		volume~68 of {\em Cambridge Studies in Advanced Mathematics}.
		\newblock Cambridge University Press, Cambridge, 2013.
		
		\bibitem{Sav}
		M.~Savov.
		\newblock Small time two-sided LIL behavior for {L}\'{e}vy processes at zero.
		\newblock {\em Probab. Theory Related Fields} 144 (2009), no. 1-2, 79–98.
		
		\bibitem{SRV12}
		R.~L. Schilling, R.~Song,  Z.~Vondra\v{c}ek.
		\newblock {\em Bernstein functions}, volume~37 of {\em De Gruyter Studies in
			Mathematics}.
		\newblock Walter de Gruyter \& Co., Berlin, second edition, 2012.
		\newblock Theory and applications.
		
		\bibitem{SW19}
		Y.~Shiozawa, J.~Wang.
		\newblock Long-time heat kernel estimates and upper rate functions of
		{B}rownian motion type for symmetric jump processes.
		\newblock {\em Bernoulli}  25 (2019), no. 4B, 3796–3831.
		
		\bibitem{Str64}
		V.~Strassen.
		\newblock An invariance principle for the law of the iterated logarithm.
		\newblock {\em Z. Wahrsch, und Verw. Gebiete}  3 (1964), 211–226.
		
	
		
		
		\bibitem{Wat64}
		S.~Watanabe.
		\newblock On discontinuous additive functionals and {L}\'{e}vy measures of a
		{M}arkov process.
		\newblock {\em Jpn. J. Math.}  34 (1964), 53–70.
		
		\bibitem{WK88}
		I.~S. Wee, Y.~K. Kim.
		\newblock General laws of the iterated logarithm for {L}\'{e}vy processes.
		\newblock {\em J. Korean Statist. Soc.} 17 (1988), no. 1, 30–45.
		

		
		\bibitem{Zh97}
		Q.~S. Zhang.
		\newblock Gaussian bounds for the fundamental solutions of {$\nabla (A\nabla
			u)+B\nabla u-u_t=0$}.
		\newblock {\em Manuscripta Math.}  93 (1997), no. 3, 381–390.
		
	\end{thebibliography}
\end{document}